\documentclass[10pt,journal]{IEEEtran}

\usepackage{cite}
\usepackage{amsmath}
\usepackage{graphicx}
\usepackage{listings}
\usepackage{framed}
\usepackage{hyperref}
\usepackage{subcaption}

\begin{document}

\title{Generating Optimally Focal and Intense Current Patterns in tES via Metaheuristic L1-L1 Search: Interior-Point vs.\ Simplex Algorithms}

\author{Fernando~Galaz~Prieto,
        Maryam~Samavaki,
 and~Sampsa~Pursiainen% <-this % stops a space
\IEEEcompsocitemizethanks{\IEEEcompsocthanksitem F.\ Galaz Prieto, M.\ Samavaki, and S.\ Pursiainen are with Computing Sciences, Faculty of Information Technology, Tampere University, Korkeakoulunkatu 3, 33072 Tampere, Finland.\protect\\

E-mail: fernando.galazprieto@tuni.fi}
\thanks{Manuscript received }}

\markboth{Journal of }%
{Galaz Prieto \MakeLowercase{\textit{et al.}}: Metaheuristic L1-Norm Fitted and Regularized Current Pattern Optimization in Transcranial Electrical Stimulation: Algorithm Effect and Performance}

\IEEEtitleabstractindextext{%
\begin{abstract}

This numerical simulation study investigates solving the L1-norm fitted and regularized (L1-L1) linear programming  (LP) problem to find a well-localized volumetric current density in transcranial electrical stimulation (tES), where a current pattern is attached through contact electrodes attached on the skin to create a stimulus in a targeted brain region. We consider a metaheuristic optimization process where the problem parameters are selected so that the final solution found is optimal with respect to given metacriteria, e.g., the intensity and focality of the volumetric current density in the brain. We focus on the effect of the LP-algorithm on the solution. We examine interior-point and simplex algorithms, which constitute two major alternative ways to solve an LP-task; the interior-point algorithms are based on determining a feasible solution set to allow finding an optimizer via Newton's method, while the simplex methods take steps along the edges of a polytope, subdividing the set of candidate solutions into simplicial polygons. Interior-point methods stand out among them for their complexity and predictability in terms of convergence properties at the same time. To find the current pattern, we apply five alternative optimization toolboxes: Matlab, MOSEK, Gurobi, SDPT3, and SeDuMi. We suggest that the mutual differences between these vary based on the placement of the target as well as the resolution of the two-parameter lattice. In the numerical experiments, we investigate maximizing the focality and intensity of the L1-L1 optimized stimulation current and, in the latter regard, examine its relationship to the reciprocal current pattern, maximizing the focused current density. 

\end{abstract}

\begin{IEEEkeywords}
Trancranial Electrical Stimulation (tES), Linear Programming, L1-Norm Fitting, LP-Solvers 
\end{IEEEkeywords}}

\maketitle

\IEEEdisplaynontitleabstractindextext
\IEEEpeerreviewmaketitle

\section{Introduction}
\label{sec:introduction}

\IEEEPARstart{T}{his} study focuses on transcranial Electrical Stimulation (tES), where currents are applied non-invasively through a montage of contact electrodes attached to the skin \cite{herrmann2013transcranial}. We aim at finding an accurate, robust, and fast linear programming algorithm to optimize a current pattern via the metaheuristic L1-norm fitted and regularized (L1-L1) method recently proposed in \cite{galazprieto2022l1} as an approach to create a focused and intense current density inside the brain. In this approach, L1-norm fitting is applied to minimize a regularized L1-norm difference between a targeted field and a volumetric current distribution, and regularization is performed in order to control the number of nonzero currents in the current pattern, i.e., the sparsity of the current pattern.

L1-L1 is a linear programming (LP) approach which, on top of optimizing a regularized fit between a given target and a volumetric stimulus current density resulting from the injection, performs a two-dimensional metaheuristic search to set the regularization parameter $\alpha$ and a relative threshold $\varepsilon$ defining  the focality of the solution, respectively.  The outcome will be optimal with respect to a given metacriterion, which can be, for example, the current density $\Gamma$ at the targeted location and/or the ratio $\Theta$ between $\Gamma$ and the average field outside the target (nuisance field). In \cite{galazprieto2022l1}, L1-L1 is suggested to be beneficial in the  optimization of $\Theta$ between $\Gamma$ as compared to two convex L2-norm fitting techniques:  L1-regularized semidefinite programming (SDP) approach (L1L2) as well as a Tikhonov-regularized version of weighted least squares optimization proposed in \cite{dmochowski2011optimized,dmochowski2017optimal}. The general advantages of L1-fitting in those applications of inverse modelling and imaging, where spatial distributions need to be well localized, are widely known \cite{bertero2020introduction,kaipio2006statistical}.

Compared to the state-of-the-art in non-linear convex current pattern optimization  \cite{fernandez2020unification,wagner2016optimization,dmochowski2011optimized}, the L1-L1 constitutes a  methodological leap forward; namely, whereas the previously published methods concentrate on maximizing the focused field locally in the vicinity of the target, given an upper bound for the nuisance field, the L1-L1 performs the fitting everywhere in the brain, optimizing both the focused and nuisance field to a degree set by the focality parameter $\varepsilon$. Hence, L1-L1 finds a globally optimal volumetric fit, while the state-of-the-art solutions can be considered optimal only in the local environment of the target. A global fit found by L1-L1 can be considered to enhance nuisance field suppression compared to local maximization, where the  maximum obtainable focused field intensity might be limited {\em per se}, unless the neighborhood of the target and those areas where the nuisance field is constrained are carefully separated {\em a priori}.

As an LP technique, L1-L1 is straightforward to implement and achieves an appropriate computational performance. Being fundamental in solving convex optimization tasks, well-performing LP algorithms are available in various optimization packages, which typically provide interior-point and simplex methods and their hybrids for LP optimization. This study aims to evaluate and compare the performance of interior-point  \cite{mehrotra1992implementation} and simplex algorithms \cite{Boyd2004} in solving the  L1-L1 optimization task. These are available in the well-established and widely used commercial and open-source optimization packages, for example, in Matlab, whose LP solver \lstinline{linprog} includes both kind of algorithms  \cite{zhang1999user,mehrotra1992implementation} as well as in Mosek \cite{mosek} and Gurobi \cite{gurobi} which among other things, provide a Matlab-based alternative for  \lstinline{linprog} but also a variety of other means to access their LP routines. In this study, the commercial Matlab, Mosek and Gurobi are accompanied by the open-source SDP packages SDPT3 \cite{tutuncu2003solving} and  SeDuMi (Self-Dual Minimization) \cite{sturm1999using,sturm2000central,polik2007sedumi}, which are
applicable via the openly available CVX optimization toolbox for Matlab \cite{grant2009cvx}.

We perform our numerical analysis using a realistic head model, evaluating maps of the optimization results for the full brain and finding the computing times required by the different solvers. By performing the comparison our goal is to find a robust and fast algorithm that can approximate the actual optimizer well already with a comparably low lattice resolution for $\alpha$ and $\varepsilon$. We pay particular attention to the maximization of $\Gamma$ and $\Theta$, as well as to comparing the results obtained in the former case to those yielded by the generalized reciprocity principle \cite{fernandez2020unification}.

The numerical results obtained show significant differences between the methods investigated in this study;  the success of the metaheuristic L1-L1 search process was found to be strongly dependent on the numerical implementation and affected by the search resolution.  As a central observation, the interior-point methods outperform the simplex algorithms, while the package-based discrepancies are clear. These concern, in particular, the optimization of $\Theta$, while each solver yields  $\Gamma$-value close to its theoretical optimum, which is shown via the generalized reciprocity principle presented in  \cite{fernandez2020unification}.  The observed running time differences are notable considering the practical applicability of L1-L1.Without any restrictions, i.e., $\alpha = 0$ and $\varepsilon =1$, L1L1 can be interpreted to be equivalent to optimizing  the focused field alone. Based on the generalized reciprocity principle formulated in \cite{fernandez2020unification}, the current pattern maximizing  $\Gamma$ corresponds to a bipolar montage that can be obtained by picking the most intense pair of electrodes from the backprojected solution of the optimization task.

\section{Methods}

In tES, a real $L \times 1$ current pattern ${\bf y}$ (Ampere) is injected through contact electrodes attached on the skin into subject's head where it distributes as a volume current density (Ampere per m\textsuperscript{2}) penetrating the skull and eventually reaching the brain. The relationship between ${\bf y}$ and a real $L \times 1$ discretized volume current density $\hat{\bf x}$ is described by the governing linear system 
    $\hat{\bf L} {\bf y} = \hat{\bf x}$, 
where $\hat{\bf L}$ is a real $N\times L$ lead field matrix, which can be  obtained as shown in \cite{galazprieto2022l1}. We define the following two  volume current density components: the focused field, where the target is non-zero, and the nuisance field, where it vanishes.  The governing linear system  can be split component-wise $\hat{\bf L}_1  {\bf y} = \hat{\bf x}_1$ and  $\hat{\bf L}_2  {\bf y} =  {\bf 0}$, respectively. The optimization problem to find the best match between ${\bf y}$ and the target is approached in the following  projected form 
\begin{equation}
\label{linear_equation}
  {\bf L} {\bf y} = {\bf x},
\end{equation}
where the projection of the focused field into the direction of the target constitutes the first component, that is, 
\begin{equation}
\begin{aligned}
    {\bf L} &= \begin{pmatrix} {\bf L}_1 \\ {\bf L}_2 \end{pmatrix} = \begin{pmatrix} {\bf P} \hat{\bf L}_1 \\ \hat{\bf L}_2 \end{pmatrix}
    & \hbox{and} \quad & {\bf x} = \begin{pmatrix} {\bf x}_1 \\ {\bf 0}
    \end{pmatrix} = \begin{pmatrix} {\bf P} \hat{\bf x}_1 \\ {\bf 0} \end{pmatrix} \nonumber
\label{projected_form}
\end{aligned}
\end{equation}
with ${\bf P}$ denoting  a matrix that projects a vector into the direction of $\hat{\bf x}_1$. The target amplitude $\| {\bf x}_1 \|_2$ is assumed to be 3.85 A/m\textsuperscript{2}, which is a coarse approximation  for the excitation threshold of nerve fibers; for the upper limb area of the motor cortex, this value decreases from 6 to 2.5 A/m\textsuperscript{2} when the stimulus frequency decreases from 2.44 kHz (kilohertz) to 50 Hz (hertz) \cite{kowalski2002current}.

\subsection{L1-L1 optimization}
\label{sec:L1L1}
The L1-L1 optimization problem \cite{galazprieto2022l1} for tES is to minimize the following objective function: 
\begin{equation}
    \min_{{\bf y}} \Bigg\{\, \left\| \begin{pmatrix} {\bf L}_1 {\bf y} - {\bf x}_1 \\
    \Psi_\varepsilon [\nu^{-1} {\bf L}_2 {\bf y}] \end{pmatrix} \right\|_1 + {\alpha} \zeta \| {\bf y} \|_1 \Bigg\},
\label{objective_function}
\end{equation}
subject to ${\bf y} \preceq \gamma {\bf 1}$, $\| {\bf y}\|_1 \leq \mu$, in which $\gamma$ and $\mu$ depict the maximum applicable current and the total dose, and $\sum_{\ell = 1}^L y_\ell = 0$. Parameter ${\bf \alpha}$ sets the level of L1-regularization with respect to scaling $\zeta = \| {\bf L} \|_1$. Function $\Psi_\varepsilon[{\bf w}]_m = \max \{\, |{w_m}|, \varepsilon \,\}$ for $m = 1, 2, \ldots, M$, sets a threshold $0 \leq \varepsilon \leq 1$ for the nuisance field with respect to scaling $\nu = \| {\bf x} \|$, meaning that entries $( {\bf L}_2 {\bf y})_m$ with absolute value below $\varepsilon \nu$ do not actively contribute to the minimization process due to the threshold. The L1-L1 problem can be solved numerically in the following form \cite{Boyd2004} including a linear objective function
\begin{equation} 
\label{LP1}
    \min_{{\bf y}, {\bf t}^{(1)}, {\bf t}^{(2)}, {\bf t}^{(3)}} \!\!
    \begin{pmatrix}
    {\bf 0} \\
    {\bf 1} \\
    {\bf 1} \\
    {\bf 1}
    \end{pmatrix}^T \!\!
    \begin{pmatrix}
    {\bf y} \\
    {\bf t}^{(1)} \\
    {\bf t}^{(2)} \\
    {\bf t}^{(3)}
    \end{pmatrix} \quad {\bf t}^{(1)}, {\bf t}^{(2)},   {\bf t}^{(3)} \succeq 0
    \end{equation}
 inequality constraint
{\setlength\arraycolsep{2pt}
\begin{eqnarray}
\label{LP2}
     \begin{pmatrix} {\bf L}_1  & -{\bf I} & {\bf 0}  & {\bf 0}  \\ 
     {\bf L}_2  & {\bf 0} & -{\bf I} & {\bf 0}  \\ 
    -{\bf I} & {\bf 0} &  {\bf 0} & -{\bf I} \\ 
    -{\bf L}_1  & -{\bf I} &  {\bf 0}& {\bf 0}  \\
    -{\bf L}_2  & {\bf 0} & -{\bf I} & {\bf 0}  \\
     {\bf I} & {\bf 0} & {\bf 0}  & -{\bf I} \\
     {\bf 0} & -{\bf I} & {\bf 0}  & {\bf 0} \\     
     {\bf 0} & {\bf 0} & -{\bf I}  & {\bf 0} \\ 
     {\bf 0} & {\bf 0} & {\bf 0} & -{\bf I} \\
     {\bf 0} & {\bf 0} & {\bf 0} & {\bf I} \\ 
     {\bf 0} & {\bf 0} & {\bf 0} & {\bf 1}^T \\ 
     \end{pmatrix} \begin{pmatrix} {\bf y} \\
     {\bf t}^{(1)} \\
     {\bf t}^{(2)} \\
     {\bf t}^{(3)} \end{pmatrix} & \preceq & \begin{pmatrix} {\bf x}_1 \\
     {\bf 0} \\
     {\bf 0} \\
    -{\bf x_1} \\
     {\bf 0} \\
     {\bf 0} \\
     {\bf 0} \\
    -{\varepsilon \nu {\bf 1}} \\
     {\bf 0} \\
     \gamma {\bf 1} \\
     \mu \\ 
     \end{pmatrix} . 
     \label{eq:LP} 
\end{eqnarray}}
together with the linear equality  \begin{equation} \label{LP3} {\bf 1}^T {\bf y}  = 0 \end{equation}
To obtain the matrix and right-hand side vector of the linear system. The solution can be found in a straightforward manner using the LP optimizer functions available in optimization libraries. \label{sec:L1-L1}

\subsection{Optimal focality and intensity}
\label{sec:focality_and_intensity}

Given a target current, we consider two  meta-optimization tasks to obtain maximally focal  and intense  volumetric current density. The intensity is defined as the projection 
\begin{equation}
    \label{gamma_def}
    {\Gamma} = \frac{{\bf x}_1^T {\bf L}_1 {\bf y}}{\| {\bf x}_1 \|_2},
\end{equation} of the discretized volumetric current density ${\bf L} {\bf y}$ created by the current pattern ${\bf y}$ into the direction of the target current at the  target location. The focality is defined as the following current ratio:
\begin{equation}
\label{theta_def}
    {\Theta} = \frac{ {\Gamma}}{ \| {\bf L}_2 {\bf y} \|_2 / \sqrt{ M }}.
\end{equation}

The maximum intensity is obtained as  \begin{equation} \Gamma_{\hbox{\scriptsize{max}}} =   \arg \max_{{\bf y}, \alpha, \varepsilon} \Gamma, \end{equation} 
and maximum focality  as  \begin{equation}
\label{max_theta} \Theta_{\hbox{\scriptsize{max}}} = \arg \max_{{\bf y}, \alpha, \varepsilon} \Theta \quad \hbox{subject to} \quad \Gamma \geq \Gamma_0. \end{equation} The latter one of these tasks can be interpreted to be metaheuristic, as a subjectively selected metacriterion $\Gamma \geq \Gamma_0$  is applied to maintain an appropriate intensity at the target location; this is necessary, as without a lower bound for $\Gamma$ the intensity of the maximizer is likely to vanish. 

The meta-tasks are solved approximately by discretizing the two-dimensional plane of $\varepsilon$ and $\alpha$ using a regular lattice with steady logarithmic increments in each direction. The actual optimization task is solved in each lattice point while the approximate solution of the meta-task is found by picking the optimum out of the candidate solutions obtained at the lattice points. The metacriterion is chosen to be $\Gamma_0 = (3/4) \Gamma_{\hbox{\scriptsize{max}}}$.

\subsection{Reciprocity principle}
\label{sec:reciprocity_principle}

The generalized reciprocity principle (GRP) \cite{fernandez2020unification} is a simple approach to obtain $\Gamma_{\hbox{\scriptsize{max}}}$ based on the reciprocity of the electromagnetic field propagation. Since it finds an exact estimate, we use it as a reference technique to evalueate the performance of the optimization algorithms. 

The classical version of GRP considers the connection between the forward and reverse propagation of the electromagnetic field, which in the case of the electric field caused by  brain activity, is predicted by the lead field matrix  of Electroencephalography (EEG)  \cite{fernandez2020unification,malmivuo1996bioelectromagnetism}. GRP states that the montage yielding $\Gamma_{\hbox{\scriptsize{max}}}$ is bipolar and can be obtained by picking the two greatest EEG electrode voltages corresponding to the desired target current in the brain. Thus, these voltages can be obtained via multiplying (projecting) the target current by the EEG lead field matrix.  In this study, we define GRP analogously for the tES lead field matrix ${\bf L}$ which is formulated as a mapping between electrode currents and volume current field \cite{galazprieto2022l1}. While in the general electromagnetism, Gradient propagation is not generally reciprocal,  it can be shown  (Appendix \ref{sec:grp}) that the bipolar montage coincides with the greatest absolute backprojected currents in the vector ${\bf L^T} {\bf x}_1$. 

\subsection{LP methods}

As methods to solve the LP task of Section \ref{sec:L1-L1} we consider interior point (IP), primal simplex (PS) and dual simplex (DS) algorithms \cite{dikin1967iterative,Boyd2004}. Of these, IP methods can be considered the most advanced. Those utilize Newton's method to operate in the interior of a feasible set. IP algorithms can be subdivided into primal-dual (predictor-corrector) IP methods \cite{Boyd2004,mehrotra1992implementation,fiacco1964sequential}, and  barrier methods, which determine the feasible set via a barrier function. Simplex methods are based on seeking the solution by considering the feasible set as a convex polytope and moving along its edges. This strategy is less memory-consuming than the IP iteration, but has weaker and not fully predictable convergence properties for large problem sizes compared to IP  iterations whose computational complexity can be obtained explicitly \cite{Boyd2004}.

\subsection{Primal and dual solvers}

The concepts of primal and dual refer to the formulation of the LP problem; noting that by presenting the entries of ${\bf y}$ as  differences of non-negative variables ($y_i = s_i - p_i$, $s_i, p_i \geq 0$) and presenting the equality constraint (\ref{LP3}) via two inequalities (condition $a = 0$ is satisfied $a \leq 0$ and $-a \leq 0$)  the task (\ref{LP1})--(\ref{LP3}) can be brought back to the following standard primal formulation: \begin{equation} \max_{{\bf z}} \mathbf{c}^{T} \mathbf{z} \quad \hbox{subject to} \quad   {\bf A} \mathbf{z} \leq \mathbf{b}, \quad \mathbf{\bf z} \geq 0, \end{equation} whose dual  is given by \begin{equation} \min_{\hat{\bf z}} \mathbf{b}^{T} \hat{\mathbf{z}} \quad \hbox{subject to} \quad {\bf A}^{T} \hat{\mathbf{z}} \geq \mathbf{c}, \quad \hat{\mathbf{z}} \geq 0. \end{equation}

\subsection{Solver packages}
\label{sec:solver_packages}

\begin{table*}[h!]
    \centering
        \caption{Detailed descriptions LP solvers applied in this study. All solvers were interfaced with Matlab's version R2020b and called  by  L1-L1 optimizer of the  Zeffiro Interface (ZI) toolbox \cite{he2019zeffiro}. As computing hardware, we used Dell 5820 workstation with a 10-core Intel Core i9-10900X processor and 256GB RAM.}
    \begin{tabular}{lllll}
    \hline
    Abbreviation & Solver & Interface & Method name & Method type \\
    \hline 
      Gurobi IP &   Gurobi Optimizer 9.5.2  & Gurobi toolbox &  interior-point& parallel barrier   \\
         Gurobi PS &  Gurobi Optimizer 9.5.2   & Gurobi toolbox &  simplex & primal \\
             Gurobi DS &   Gurobi Optimizer 9.5.2   &Gurobi toolbox &  simplex & dual \\
        Matlab IPL & R2020b linprog / LIPSOL &  Optimization toolbox & interior-point legacy & primal-dual \\
                 Matlab PS &  R2020b linprog   & Optimization toolbox & simplex & primal \\
             Matlab DS &  R2020b linprog & Optimization toolbox  &  simplex & dual \\
        Mosek IP &   Mosek 9.1.0 & Mosek toolbox & interior-point & primal-dual  \\
             Mosek PS &   Mosek 9.1.0 & Mosek toolbox & simplex & primal  \\
              Mosek DS &   Mosek 9.1.0 & Mosek toolbox & simplex & dual  \\
       SDPT3 IP &  SDPT3 4.0 & CVX 2.2  &  interior-point  & primal-dual \\
        SeDuMi IP &  SeDuMi 1.1 & CVX 2.2 &  interior-point  & primal-dual \\
        \hline
    \end{tabular}
    \label{tab:solvers}
\end{table*}

The IP algorithms applied this study include  Gurobi's  parallel barrier  method and the primal-dual routine of Matlab, Mosek, SDPT3, and  SeDuMi. The simplex methods applied include Mosek's PS and DS, Gurobi's PS and DS and Matlab's DS algorithm. Matlab's Optimization toolbox has two IP solvers of which we apply the interior-point legacy (IPL) whose origin is in the Linear-programming Interior Point SOLvers (LIPSOL) package \cite{zhang1999user}. The solvers, their types and  abbreviations are described in Table \ref{tab:solvers}. 

\subsection{Computing platform}

The numerical simulations of this study were performed using Dell 5820 workstation with a 10-core Intel Core i9-10900X processor and 256GB RAM. The L1-L1 solver was implemented using the Matlab-based Zeffiro Interface (ZI) toolbox\footnote{\url{https://github.com/sampsapursiainen/zeffiro_interface}} \cite{he2019zeffiro, rezaei2021reconstructing, miinalainen2019realistic} which builds a high-resolution finite element (FE) mesh and generates a tES lead field matrix \cite{galazprieto2022l1} for a given surface-based head segmentation incorporating the Complete Electrode Model (CEM) \cite{pursiainen2012, pursiainen2017advanced}. The solver packages (Section \ref{sec:solver_packages}) were accessed by ZI via their Matlab (R2020b) interfaces (\ref{tab:solvers}).

\subsection{Domain}

A realistic 1 mm FE mesh was applied as the domain of the numerical simulations. This mesh was  generated  based on a an open T1-weighted Magnetic Resonance Imaging (MRI) dataset\footnote{\url{https://brain-development.org/ixi-dataset/}}. The data were first segmented using FreeSurfer Software Suite{\footnote{\url{https://surfer.nmr.mgh.harvard.edu/}}} which finds the complex surface boundaries between different tissue compartments, including the skin, skull, cerebrospinal fluid (CSF), gray and white matter and subcortical structures such as brain stem, thalamus, amygdala, and ventricles  \cite{fischl2012freesurfer}. The conductivity values of the different tissue layers were set according to \cite{dannhauer2010}.

\subsection{Current field and targets}
\label{sec:current_field}

To solve the inverse problem, we discretized the current field using 563 spatial degrees of freedom which were evenly spread over the brain, i.e., those compartments that are enclosed by the skull and CSF compartments in the model. Each target position contained three divergence-free Cartesian field components. To maintain the total dose, i.e., the total current flowing through all the electrodes within acceptable limits, we choose $\mu = 4$ mA which is a typical value in experimental studies  along with 2 mA. The limit $\gamma$ for the maximum injectable current is set to be half of the maximum dose, i.e., $\gamma = \mu/2$. 

\subsection{Two-stage lattice search}

The meta-optimization tasks for $\Theta_{\hbox{\scriptsize{max}}}$ and $\Gamma_{\hbox{\scriptsize{max}}}$ were calculated using a  two-dimensional lattice covering  the ranges from -160 to 0 decibels (dB) for  $\varepsilon$ and from -100 to -20 for $\alpha$. Two different lattice resolutions $15 \times 15$ and $40 \times 40$ were used to examine the sensitivity of the algorithms to parameter variation. To obtain the best the optimization outcome, the meta-optimization process was run in two stages; in the first one the full set of 128 electrodes were included in the optimization process, and in the second one, the montage of the most intense 20 electrodes was fixed and the optimization process was performed in this subset. The final result was thresholded in order to count the number of  non-zero (NNZ)  currents in the pattern: the values smaller than 0.1 \% the maximum current in the pattern were set to zero.   

\begin{figure}[h!]
\centering
    \begin{subfigure}{2.5cm}
        \begin{minipage}{2.5cm}
        \raggedright
        \includegraphics[height=2.0cm]{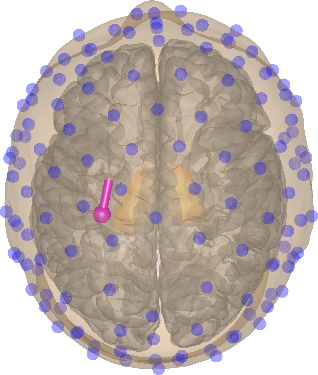}    
        \end{minipage}
        \begin{minipage}{2.5cm}
        \raggedright
        \includegraphics[height=2.0cm]{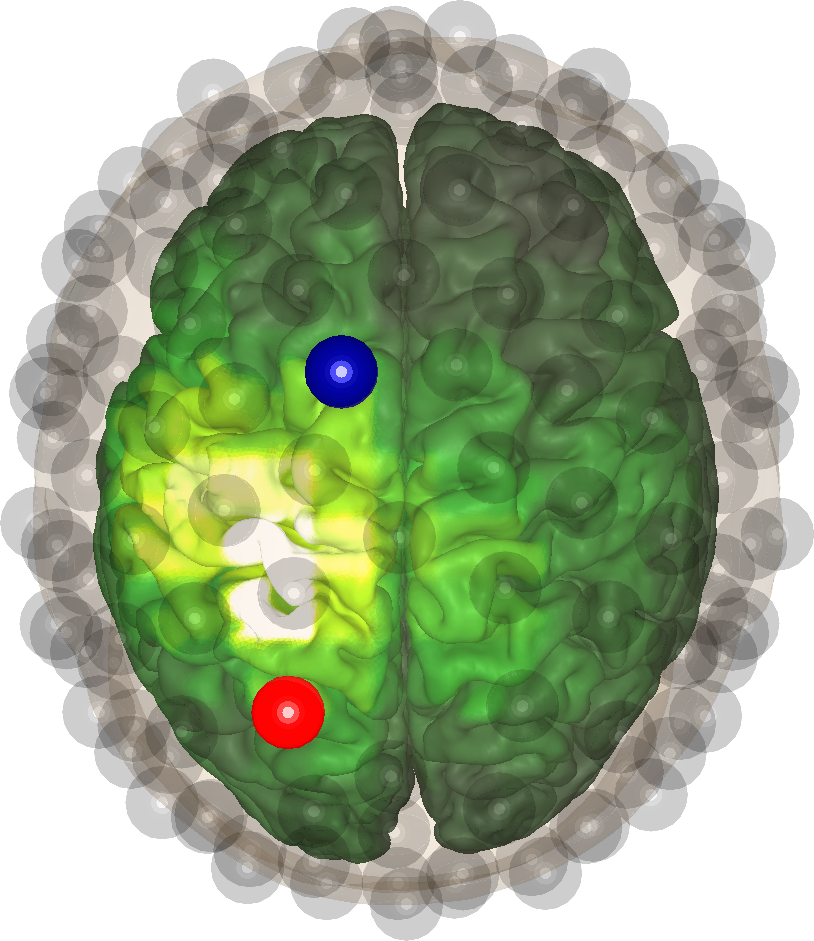}
        \includegraphics[width=0.5cm,height=2.0cm]{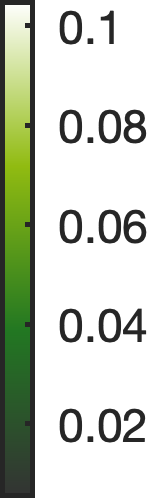}
        \end{minipage}
        \caption{Somatosensory}
        \label{fig:target_somato}
    \end{subfigure}
    \begin{subfigure}{2.5cm}
        \begin{minipage}{2.5cm}
        \raggedright
        \includegraphics[height=2.0cm]{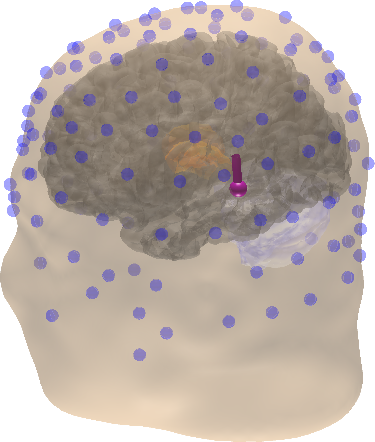}    
        \end{minipage}
        \begin{minipage}{2.5cm}
        \raggedright
        \includegraphics[height=2.0cm]{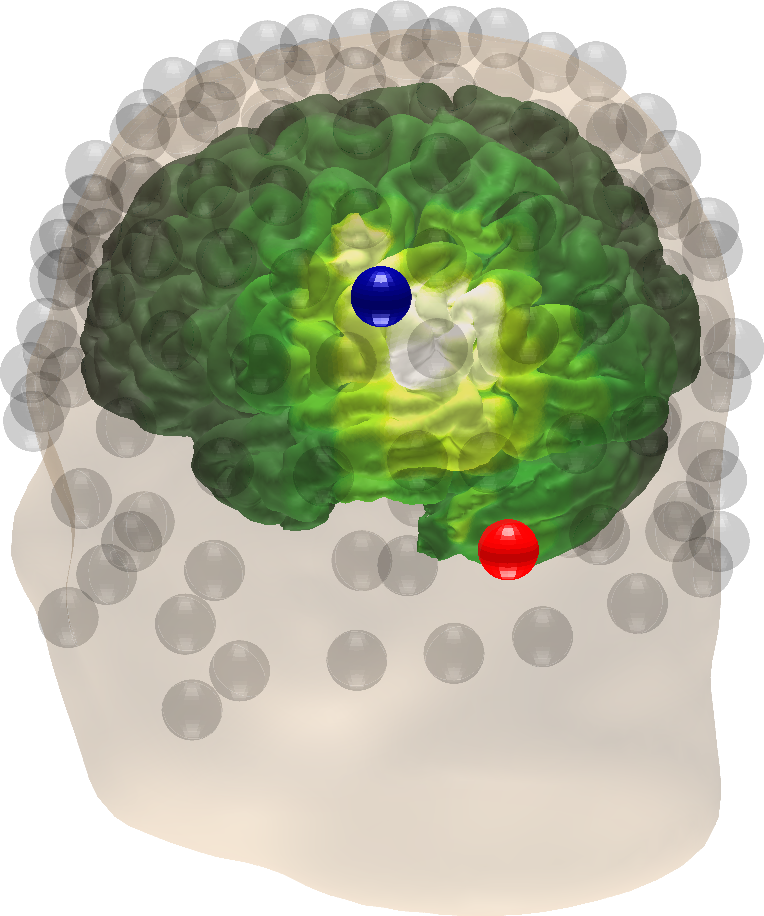}
        \includegraphics[width=0.5cm,height=2.0cm]{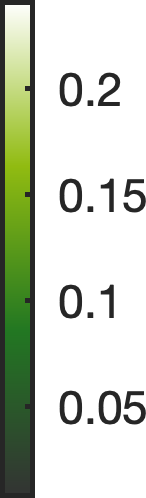}
        \end{minipage}
        \caption{Auditory}
        \label{fig:target_audio}
    \end{subfigure}
    \begin{subfigure}{2.5cm}
        \begin{minipage}{2.5cm}
        \raggedright
        \includegraphics[height=2.0cm]{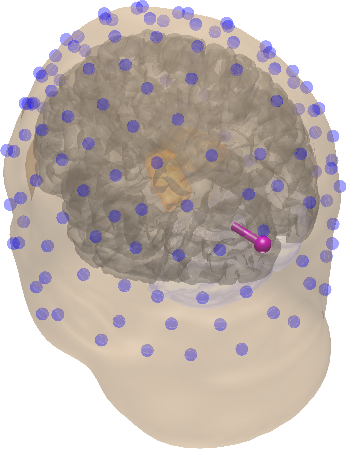}    
        \end{minipage}
        \begin{minipage}{2.5cm}
        \raggedright
        \includegraphics[height=2.0cm]{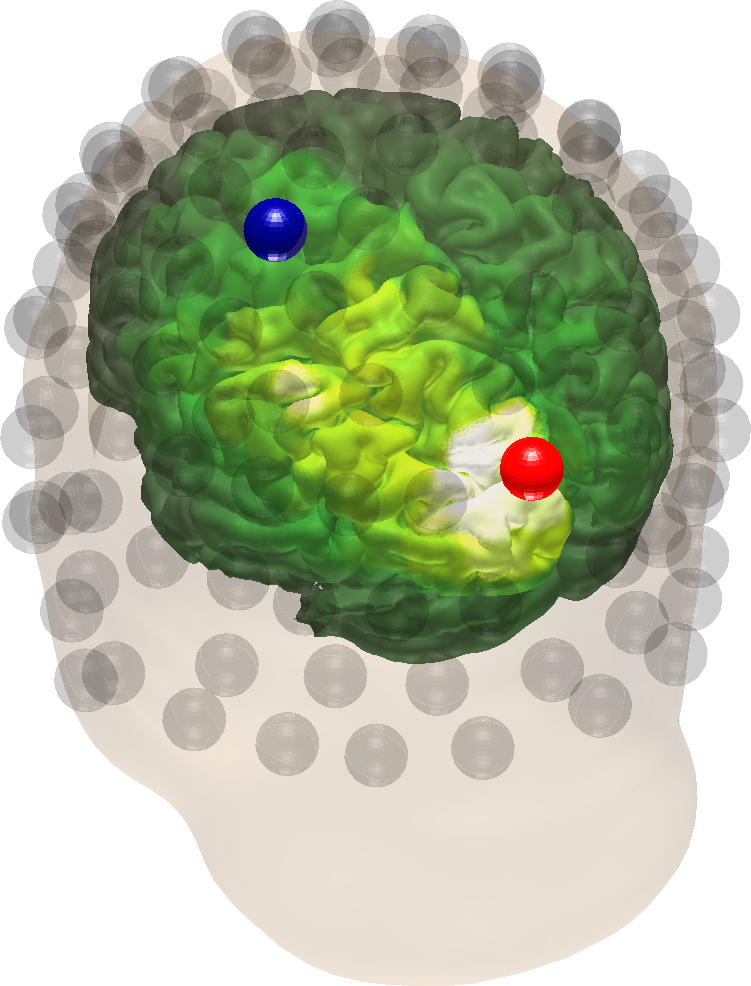}
        \includegraphics[width=0.5cm,height=2.0cm]{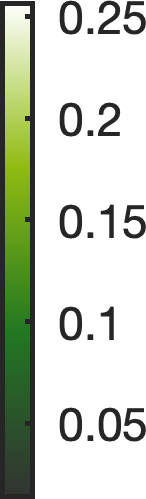}
        \end{minipage}
        \caption{Occipital}
        \label{fig:target_occipital}
    \end{subfigure}
    \caption{{\bf Top:} The head model and the distribution of the 128 electrode positions (blue) together with the placement of three targets (purple): somatosensory, auditory and occipital. {\bf Bottom:} Reciprocal bi-polar current patterns maximizing $\Gamma$ at the corresponding locations. The total dose of the current pattern is 4 mA and the volumetric current density is given in Amperes per square meter (A/m\textsuperscript{2}). }
    \label{fig:targets}
\end{figure}

% These are available in the most popular commercial and open-source packages, in particular, PDIP \cite{zhang1999user,mehrotra1992implementation} and DS algorithms implemented in Matlab's LP solver linprog as well as their alternatives available in the commercial Mosek \cite{mosek} and Gurobi \cite{gurobi} package, in particular, Mosek's PDIP and Gurobi's Parallel Barrier (PB) algorithm  together with the DS and PS implementations of these package. In addition, the PDIP algorithms of the open source semidefinite programming packages SDPT3 \cite{tutuncu2003solving} and  SeDuMi (Self-Dual Minimization) \cite{sturm1999using,sturm2000central,polik2007sedumi} are included in the comparison.

\begin{figure*}[h!]
    \centering \begin{scriptsize}
    \begin{minipage}{8.8cm}
    \centering
\rotatebox{90}{\mbox{} \hskip0.01cm Somatosensory}   \includegraphics[width = 8.5cm]{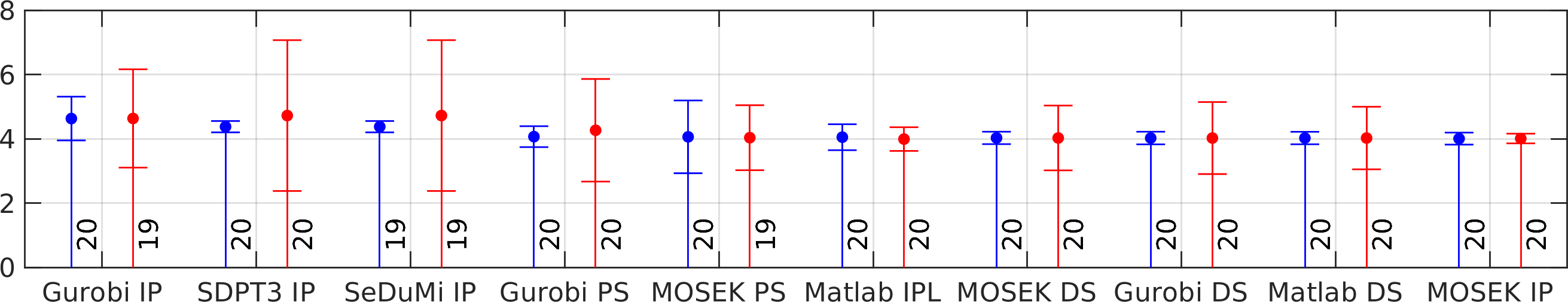}\\ \vskip0.1cm
   \rotatebox{90}{\mbox{} \hskip0.01cm Auditory}   \includegraphics[width = 8.5cm]{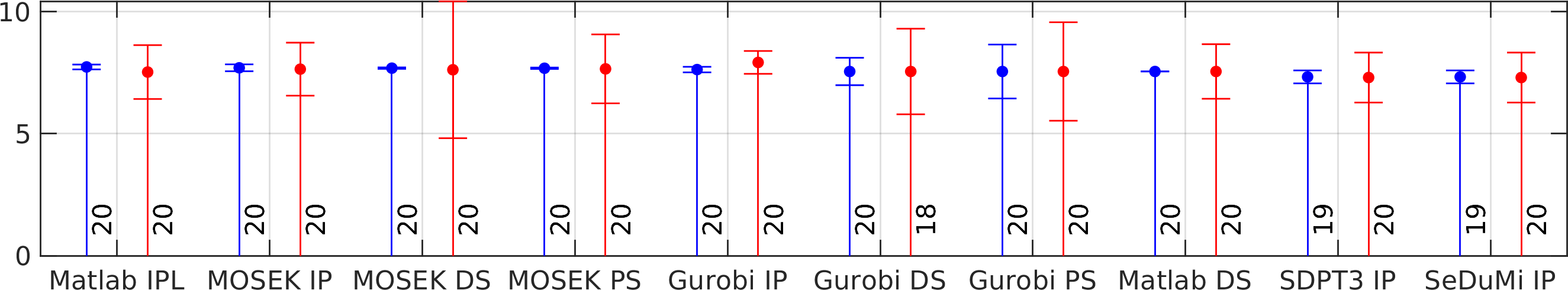} \vskip0.1cm
      \rotatebox{90}{\mbox{} \hskip0.01cm Occipital}   \includegraphics[width = 8.5cm]{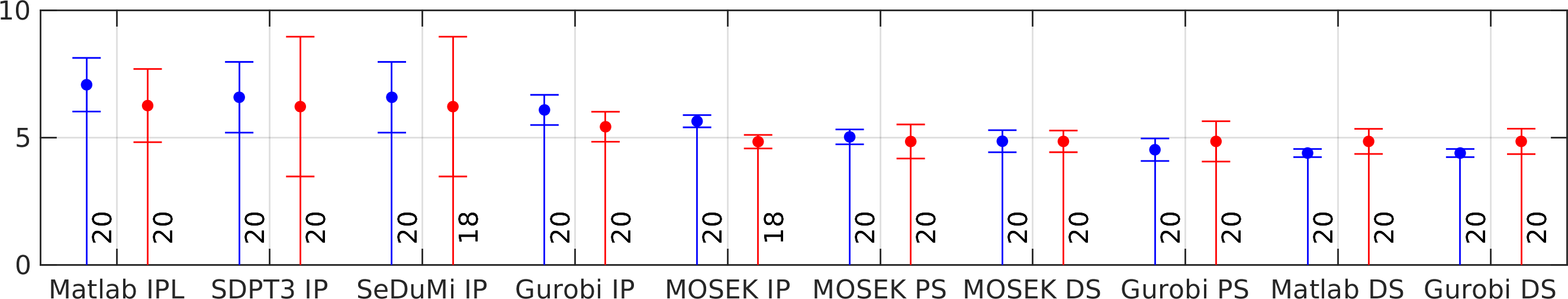} \\ \vskip0.1cm 
      $\Theta_{\hbox{\scriptsize max}}$ after the first optimization stage
%       \\ \vskip0.2cm
%       \rotatebox{90}{\mbox{} \hskip0.01cm Somatosensory}   \includegraphics[width = 8.5cm]{comparison_postcentral_1_3.png}\\ \vskip0.1cm
%   \rotatebox{90}{\mbox{} \hskip0.01cm Auditory}   \includegraphics[width = 8.5cm]{comparison_auditory_1_3.png} \vskip0.1cm
%       \rotatebox{90}{\mbox{} \hskip0.01cm Occipital}   \includegraphics[width = 8.5cm]{comparison_occipital_1_3.png} \\ \vskip0.1cm 
%       Residual in case (A) 
      \\ \vskip0.2cm
      \rotatebox{90}{\mbox{} \hskip0.01cm Somatosensory}   \includegraphics[width = 8.5cm]{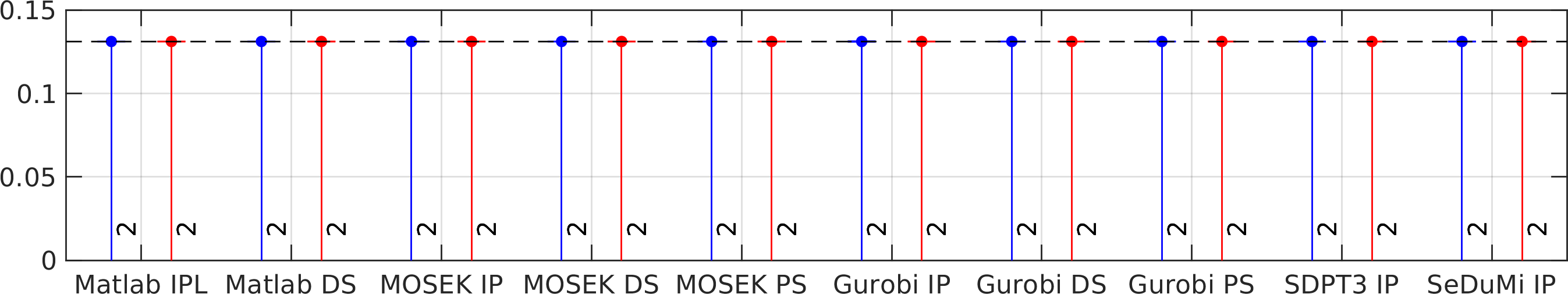}\\ \vskip0.1cm
   \rotatebox{90}{\mbox{} \hskip0.01cm Auditory}   \includegraphics[width = 8.5cm]{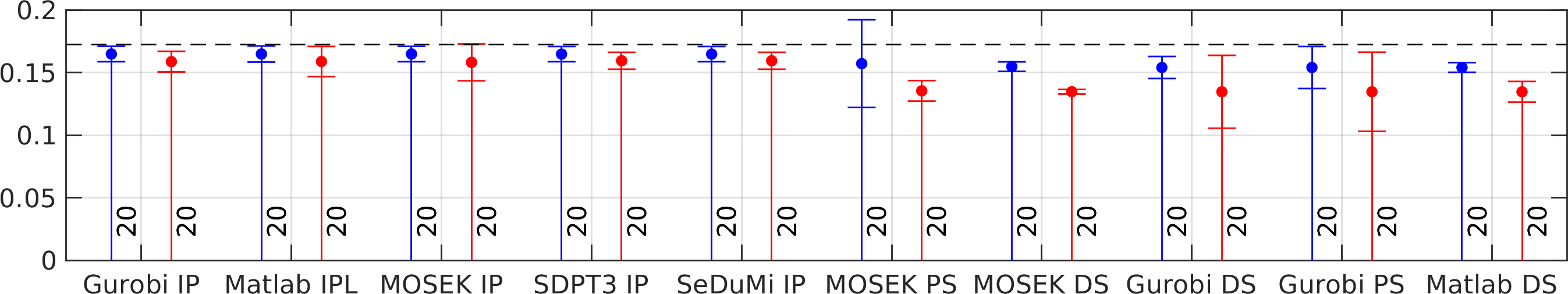} \vskip0.1cm
      \rotatebox{90}{\mbox{} \hskip0.01cm Occipital}   \includegraphics[width = 8.5cm]{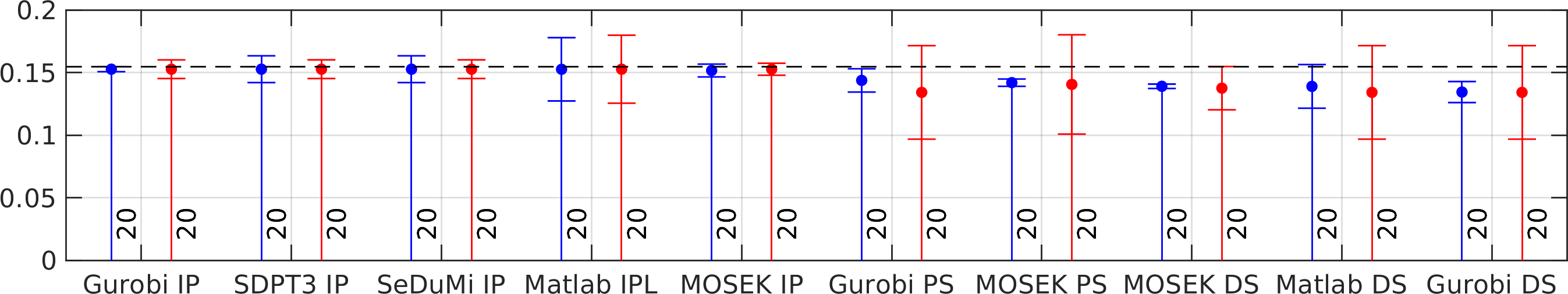}\\ \vskip0.1cm
      $\Gamma_{\hbox{\scriptsize max}}$ after the first optimization stage
       \\ \vskip0.2cm
%              \rotatebox{90}{\mbox{} \hskip0.01cm Somatosensory}   \includegraphics[width = 8.5cm]{comparison_postcentral_1_3.png}\\ \vskip0.1cm
%   \rotatebox{90}{\mbox{} \hskip0.01cm Auditory}   \includegraphics[width = 8.5cm]{comparison_auditory_1_3.png} \vskip0.1cm
%       \rotatebox{90}{\mbox{} \hskip0.01cm Occipital}   \includegraphics[width = 8.5cm]{comparison_occipital_1_3.png}\\ \vskip0.1cm
%       Residual in case (B) 
      \rotatebox{90}{\mbox{} \hskip0.01cm Somatosensory}   \includegraphics[width = 8.5cm]{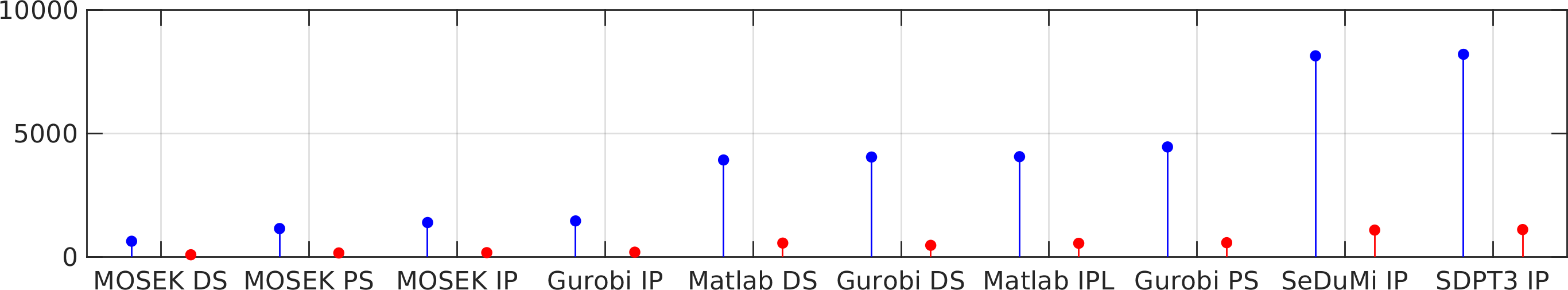}\\ \vskip0.1cm
   \rotatebox{90}{\mbox{} \hskip0.01cm Auditory}  \includegraphics[width = 8.5cm]{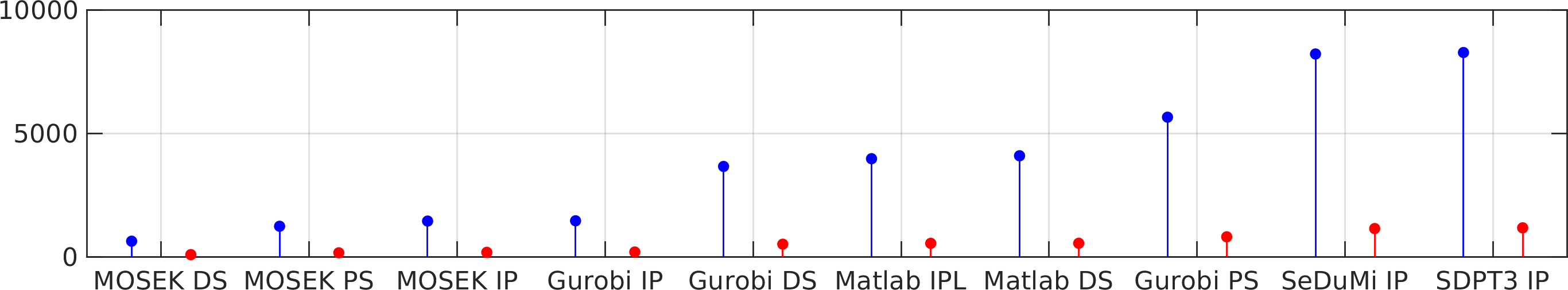} \vskip0.1cm
      \rotatebox{90}{\mbox{} \hskip0.01cm Occipital}   \includegraphics[width = 8.5cm]{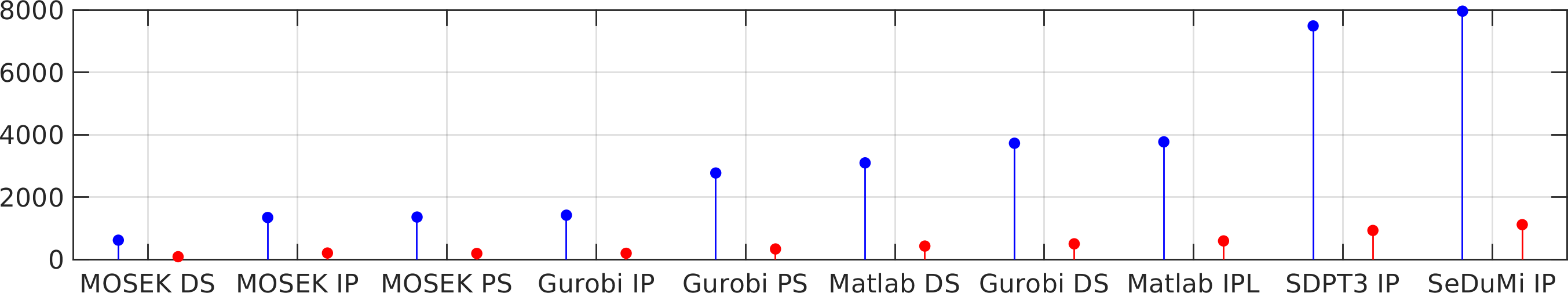} \\ \vskip0.1cm 
     Computing time for the first optimization stage
     \end{minipage}
        \end{scriptsize}
    \label{fig:my_label_computing_time_case_A}
    \centering \begin{scriptsize}
    \begin{minipage}{8.8cm}
    \centering
\rotatebox{90}{\mbox{} \hskip0.01cm Somatosensory}   \includegraphics[width = 8.5cm]{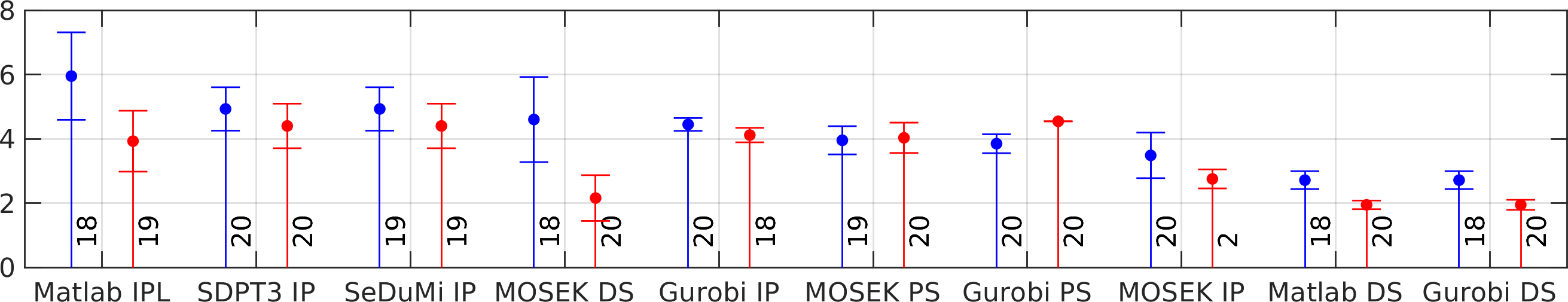}\\ \vskip0.1cm
   \rotatebox{90}{\mbox{} \hskip0.01cm Auditory}   \includegraphics[width = 8.5cm]{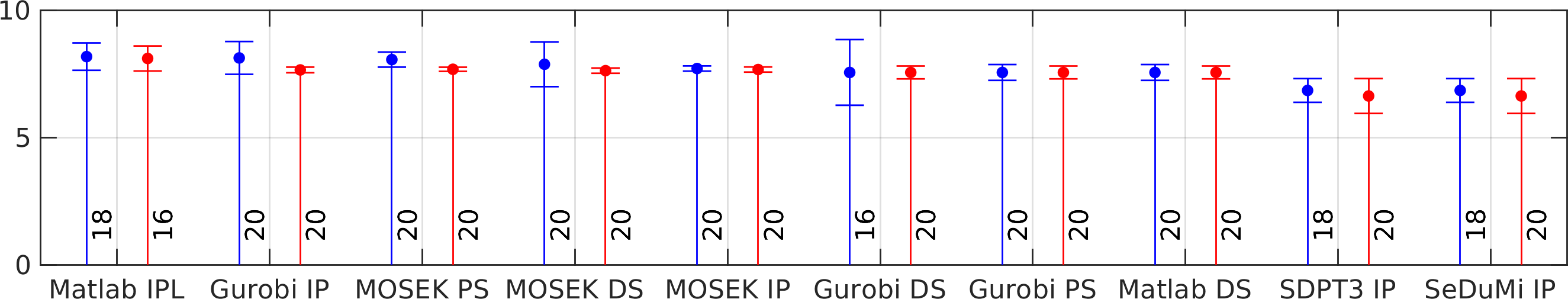} \vskip0.1cm
      \rotatebox{90}{\mbox{} \hskip0.01cm Occipital}   \includegraphics[width = 8.5cm]{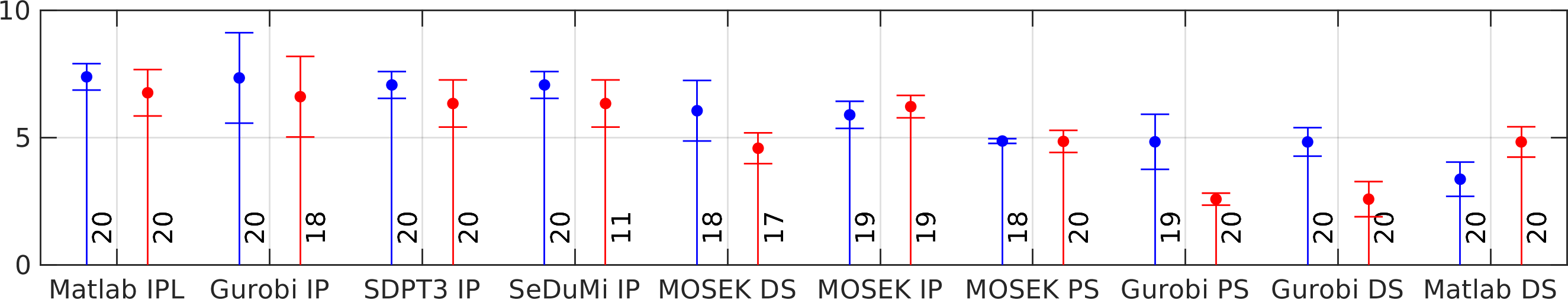} \\ \vskip0.1cm 
      $\Theta_{\hbox{\scriptsize{max}}}$ after the second optimization stage
%       \\ \vskip0.2cm
%       \rotatebox{90}{\mbox{} \hskip0.01cm Somatosensory}   \includegraphics[width = 8.5cm]{comparison_postcentral_1_3.png}\\ \vskip0.1cm
%   \rotatebox{90}{\mbox{} \hskip0.01cm Auditory}   \includegraphics[width = 8.5cm]{comparison_auditory_1_3.png} \vskip0.1cm
%       \rotatebox{90}{\mbox{} \hskip0.01cm Occipital}   \includegraphics[width = 8.5cm]{comparison_occipital_1_3.png} \\ \vskip0.1cm 
%       Residual in case (A) 
      \\ \vskip0.2cm
      \rotatebox{90}{\mbox{} \hskip0.01cm Somatosensory}   \includegraphics[width = 8.5cm]{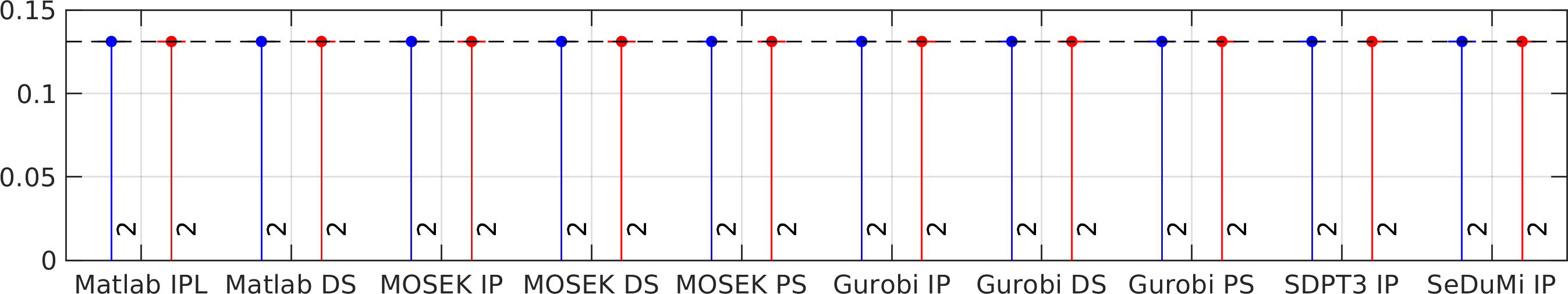}\\ \vskip0.1cm
   \rotatebox{90}{\mbox{} \hskip0.01cm Auditory}   \includegraphics[width = 8.5cm]{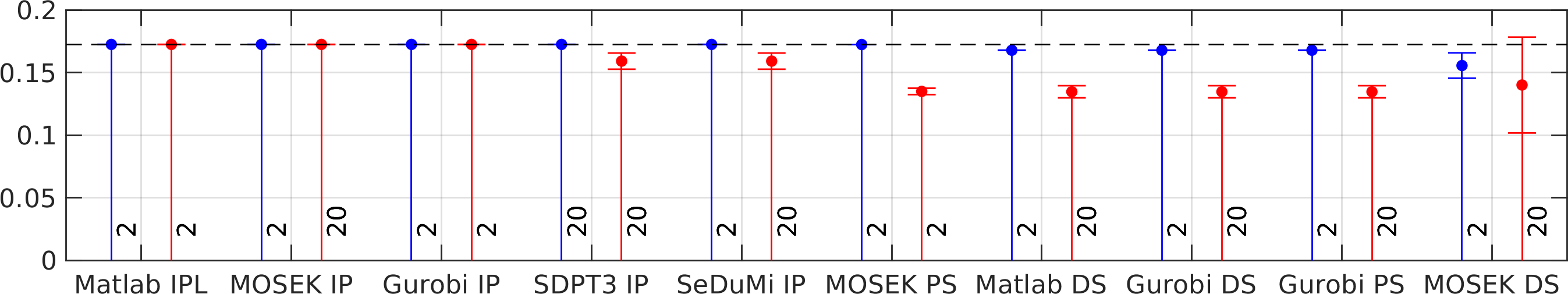} \vskip0.1cm
      \rotatebox{90}{\mbox{} \hskip0.01cm Occipital}   \includegraphics[width = 8.5cm]{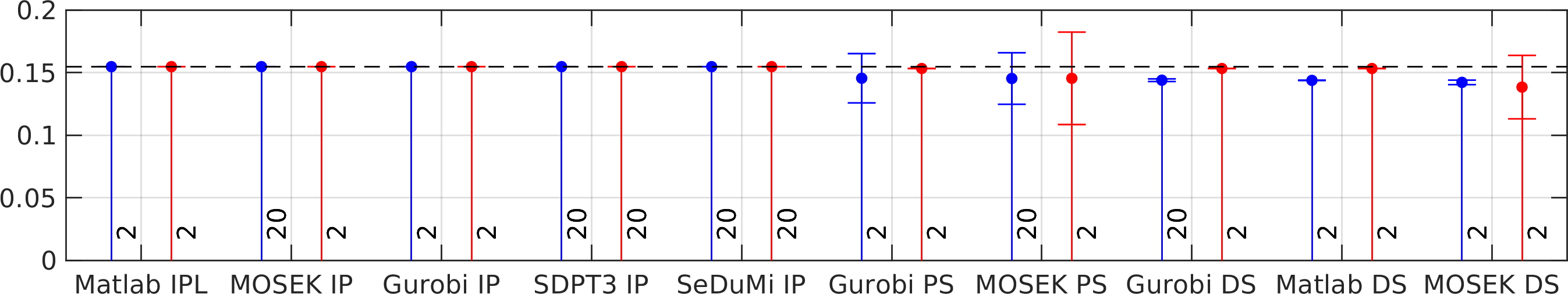}\\ \vskip0.1cm
      $\Gamma_{\hbox{\scriptsize{max}}}$  after the second optimization stage
       \\ \vskip0.2cm
%              \rotatebox{90}{\mbox{} \hskip0.01cm Somatosensory}   \includegraphics[width = 8.5cm]{comparison_postcentral_1_3.png}\\ \vskip0.1cm
%   \rotatebox{90}{\mbox{} \hskip0.01cm Auditory}   \includegraphics[width = 8.5cm]{comparison_auditory_1_3.png} \vskip0.1cm
%       \rotatebox{90}{\mbox{} \hskip0.01cm Occipital}   \includegraphics[width = 8.5cm]{comparison_occipital_1_3.png}\\ \vskip0.1cm
%       Residual in case (B) 
 \label{fig:my_label_theta_Gamma_case_B}
      \rotatebox{90}{\mbox{} \hskip0.01cm Somatosensory}   \includegraphics[width = 8.5cm]{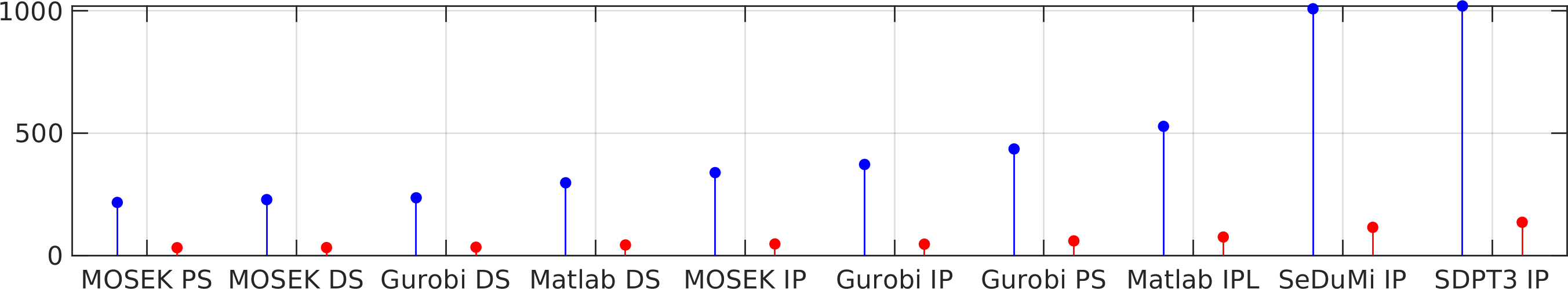}\\ \vskip0.1cm
   \rotatebox{90}{\mbox{} \hskip0.01cm Auditory}  \includegraphics[width = 8.5cm]{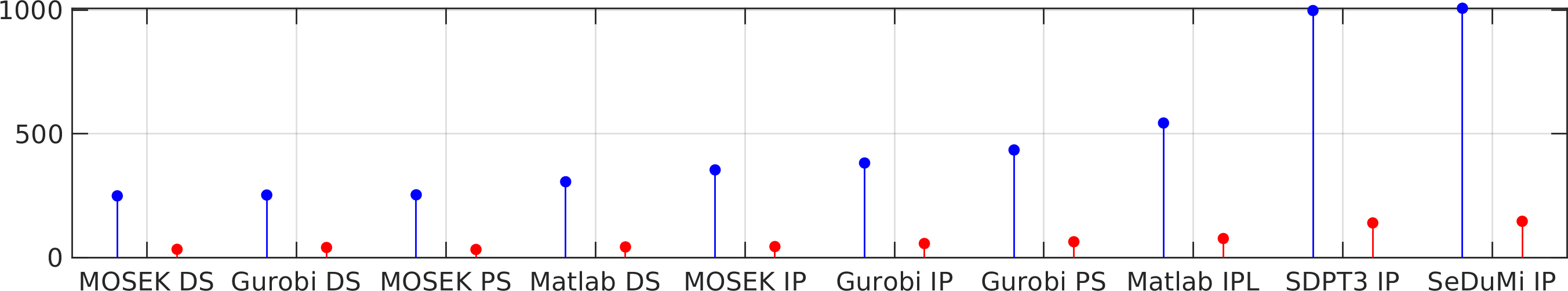} \vskip0.1cm
      \rotatebox{90}{\mbox{} \hskip0.01cm Occipital}   \includegraphics[width = 8.5cm]{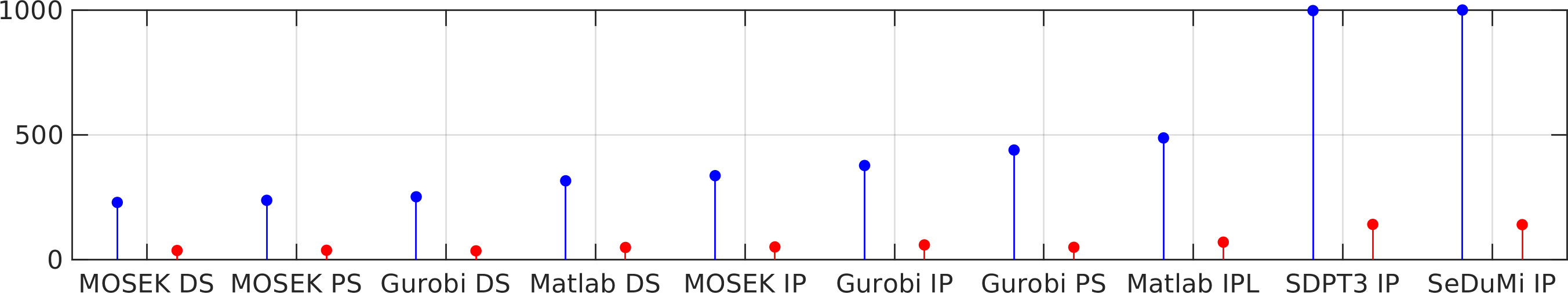} \\ \vskip0.1cm 
     Computing time in the second optimization stage
     \end{minipage}
        \end{scriptsize}
    \caption{{\bf Left:} A stem plot of the first-stage optimization results obtained for the somatosensory, auditory and occipital target. The blue and bar show the results for the denser $40 \times 40$ and coarser  $10 \times 10$ lattice, respectively. The results obtained with different optimization algorithms have been sorted into a descending order for $\Theta_{\hbox{\scriptsize max}}$ and $\Gamma_{\hbox{\scriptsize max}}$  and in an ascending order for the computing time from left to right based on the former one of these two cases. The computing time is presented as the total time required by the solvers counting together the time required by each processor core. The number of non-zero currents in the optimized pattern is shown next to the stem. The reference value obtained for $\Gamma_{\hbox{\scriptsize max}}$ using the reciprocity principle is shown as a dashed horizonal line. Approximating the accuracy of the lattice search, the whiskers for $\Theta_{\hbox{\scriptsize max}}$ and $\Gamma_{\hbox{\scriptsize max}}$ show a second order Taylor's polynomial based estimate for the maximum deviation of the result within 1/2 lattice units distance from the optimizer. The denser lattice can be observed to have a reduced deviation compared to the coarser one, demonstrating an improved reliability of the results. {\bf Right:} The second stage results can be observed to, improve the results compared to the first one which is obvious, e.g, based the active electrode counts obtained for $\Gamma_{\hbox{\scriptsize max}}$ (the bipolar exact solution has two active electrodes). The computing times are only a fraction of those required in the first stage as the limited current pattern allows a reduced representation of the lead field. }
     \label{fig:my_label_theta_Gamma_case_A}
     \label{fig:my_label_computing_time_case_B}
\end{figure*}

\begin{figure*}
    \centering
    \begin{scriptsize}
    \begin{framed}
    \begin{minipage}{4.0cm}
    \centering
    \begin{minipage}{0.20cm}
    \centering 
     \rotatebox{90}{\mbox{} \hskip0.01cm Interior-point} 
    \end{minipage}
    \begin{minipage}{1.5cm}
    \centering
    \includegraphics[height=1.1cm]{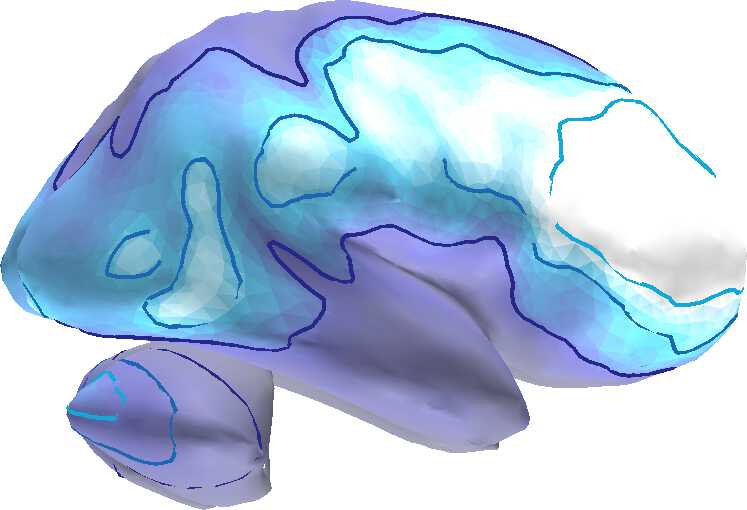}
    \end{minipage}
    \begin{minipage}{1.5cm}
    \centering
    \includegraphics[height=1.4cm]{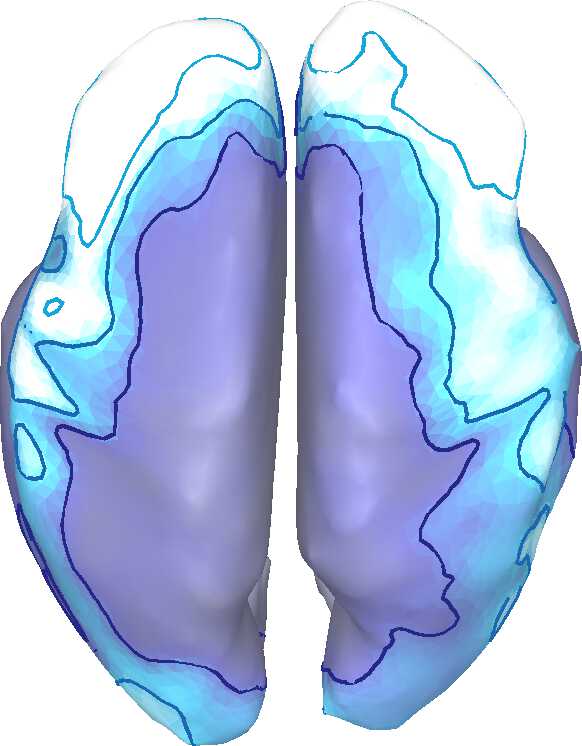}
    \end{minipage} \\ \mbox{} \vskip0.1cm
            \begin{minipage}{0.20cm}
    \centering 
     \rotatebox{90}{\mbox{} \hskip0.01cm Dual-simplex} 
    \end{minipage}
        \begin{minipage}{1.5cm}
    \centering
    \includegraphics[height=1.1cm]{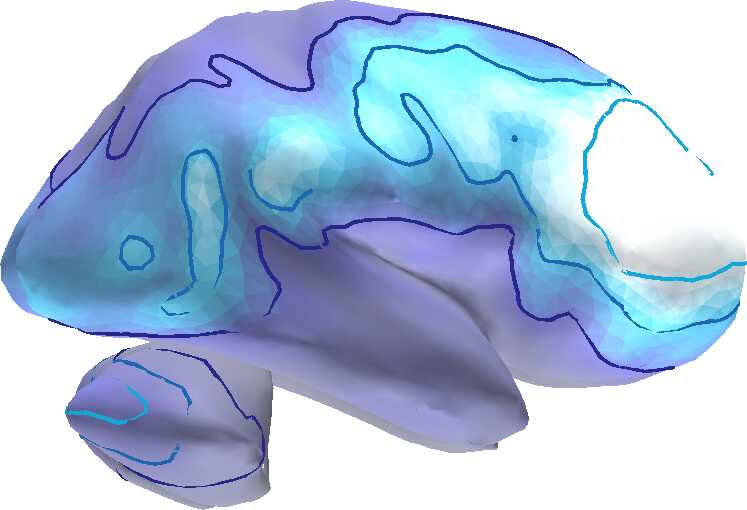}
    \end{minipage}
    \begin{minipage}{1.5cm}
    \centering
    \includegraphics[height=1.4cm]{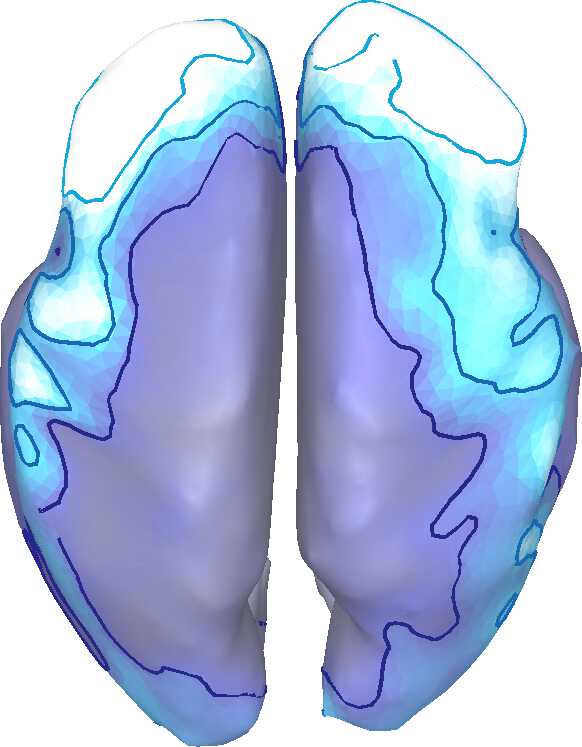}
    \end{minipage}
    \\ \mbox{} \vskip0.1cm
    \begin{minipage}{0.20cm}
    \centering 
     \rotatebox{90}{\mbox{} \hskip0.01cm Primal-simplex} 
    \end{minipage}
        \begin{minipage}{1.5cm}
    \centering
    \includegraphics[height=1.1cm]{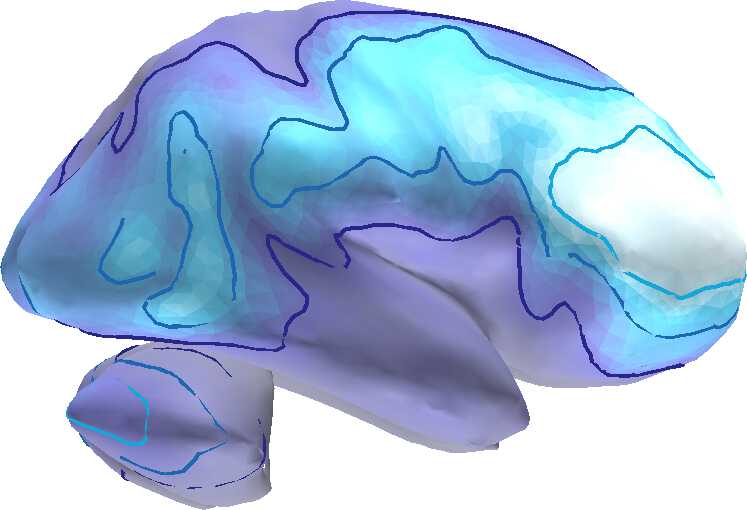}
    \end{minipage}
    \begin{minipage}{1.5cm}
    \centering
    \includegraphics[height=1.4cm]{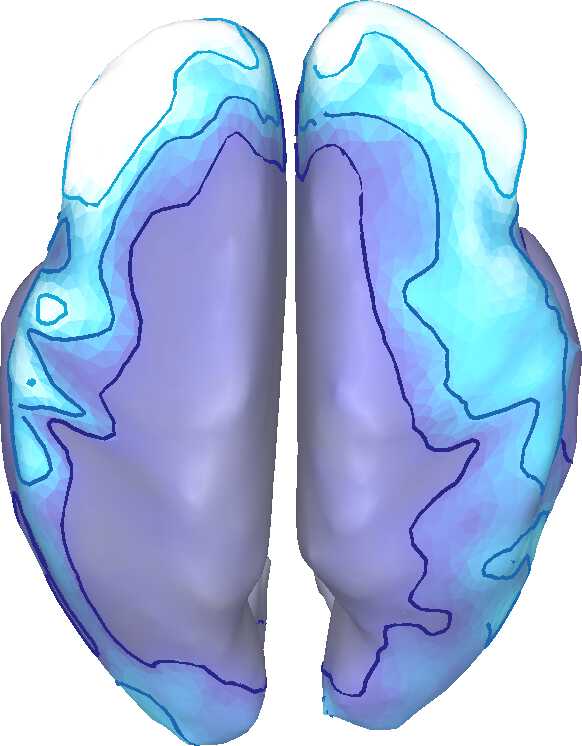}
    \end{minipage} \\ \vskip0.1cm
    $\Theta_{\hbox{\scriptsize max}}$
    \end{minipage} 
        \begin{minipage}{0.20cm} 
    \includegraphics[height=3cm]{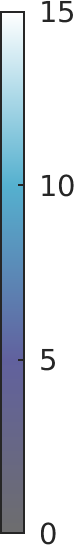}
    \end{minipage} 
            \begin{minipage}{4.0cm}
    \centering
      \begin{minipage}{0.20cm}
    \centering 
     \rotatebox{90}{\mbox{} \hskip0.01cm Interior-point} 
    \end{minipage}
    \begin{minipage}{1.5cm}
    \centering
    \includegraphics[height=1.1cm]{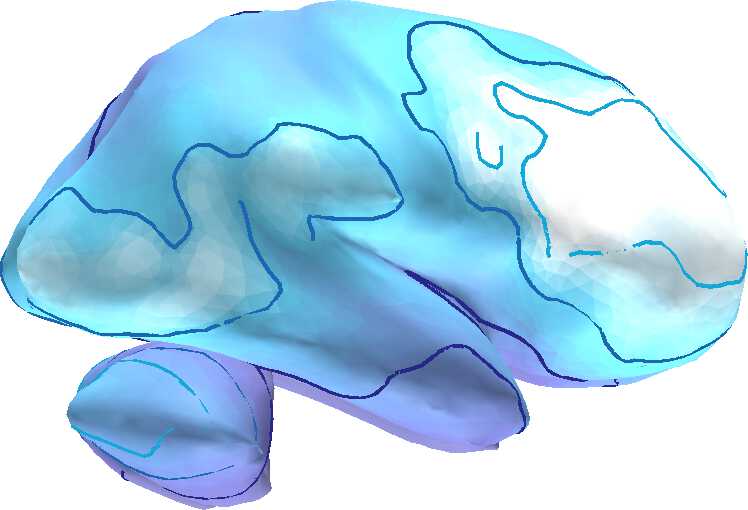}
    \end{minipage}
    \begin{minipage}{1.5cm}
    \centering
    \includegraphics[height=1.4cm]{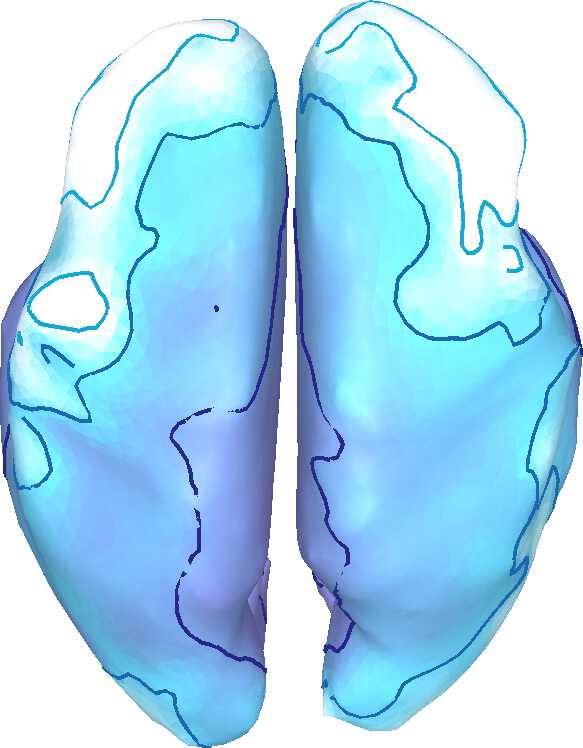}
    \end{minipage} \\ \mbox{} \vskip0.1cm
      \begin{minipage}{0.20cm}
    \centering 
     \rotatebox{90}{\mbox{} \hskip0.01cm Dual-simplex} 
    \end{minipage}
        \begin{minipage}{1.5cm}
    \centering
    \includegraphics[height=1.1cm]{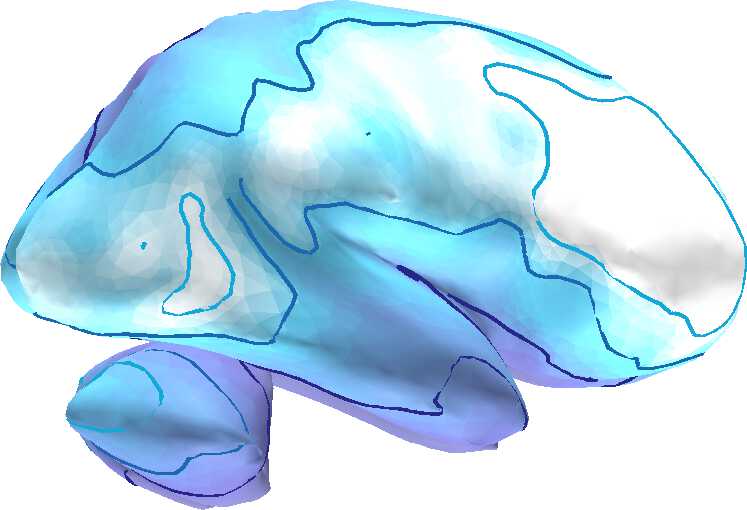}
    \end{minipage}
    \begin{minipage}{1.5cm}
    \centering
    \includegraphics[height=1.4cm]{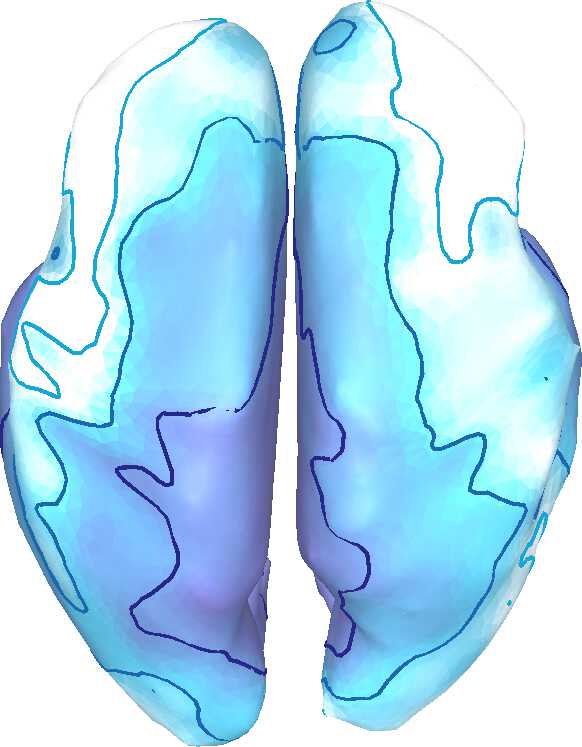}
    \end{minipage} \\
    \mbox{} \vskip0.1cm
      \begin{minipage}{0.20cm}
    \centering 
     \rotatebox{90}{\mbox{} \hskip0.01cm Primal-simplex} 
    \end{minipage}
        \begin{minipage}{1.5cm}
    \centering
    \includegraphics[height=1.1cm]{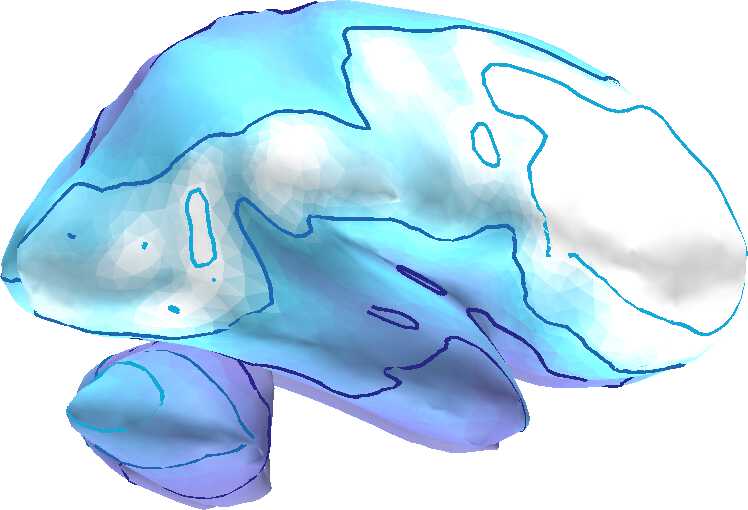}
    \end{minipage}
    \begin{minipage}{1.5cm}
    \centering
    \includegraphics[height=1.4cm]{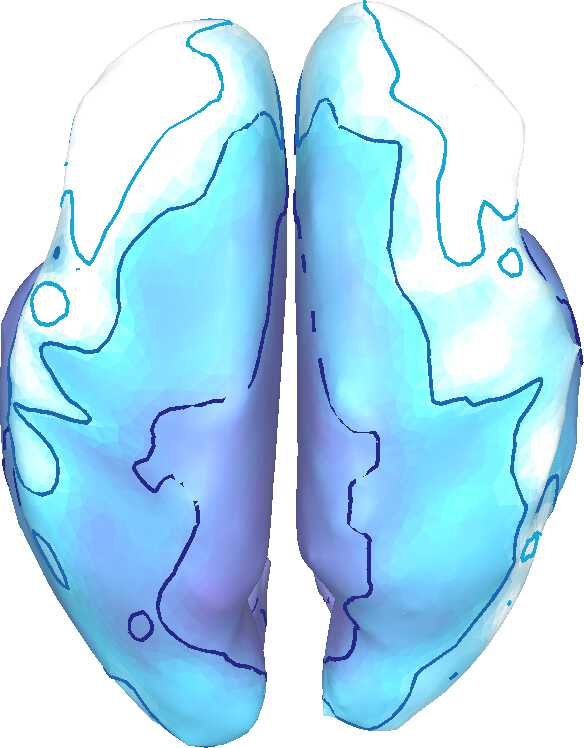}
    \end{minipage} \\ \vskip0.1cm
    $\Gamma$ 
    \end{minipage} 
    \begin{minipage}{0.20cm} 
    \includegraphics[height=3cm]{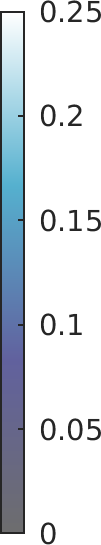}
    \end{minipage} 
    \begin{minipage}{4.0cm}
    \centering
    \begin{minipage}{0.20cm}
    \centering 
     \rotatebox{90}{\mbox{} \hskip0.01cm Interior-point} 
    \end{minipage}
    \begin{minipage}{1.5cm}
    \centering
    \includegraphics[height=1.1cm]{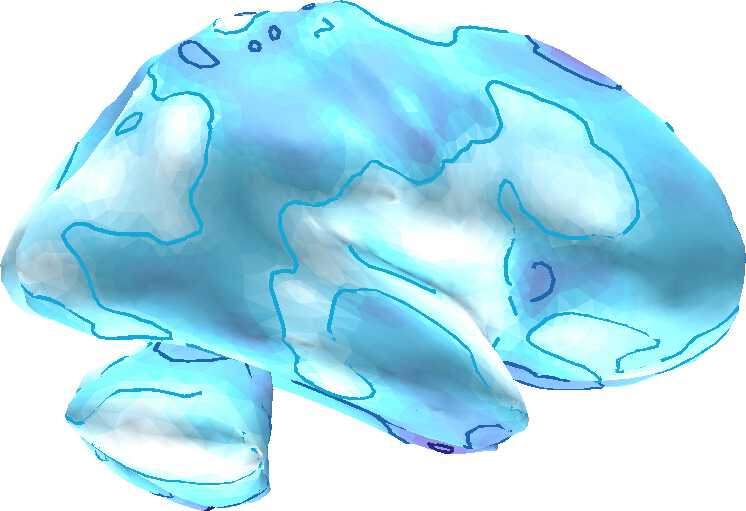}
    \end{minipage}
    \begin{minipage}{1.5cm}
    \centering
    \includegraphics[height=1.4cm]{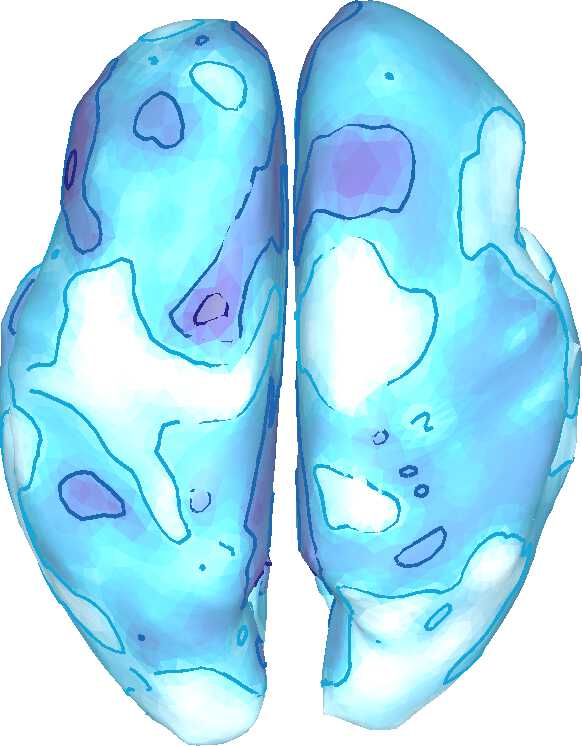}
    \end{minipage} \\ \mbox{} \vskip0.1cm
            \begin{minipage}{0.20cm}
    \centering 
     \rotatebox{90}{\mbox{} \hskip0.01cm Dual-simplex} 
    \end{minipage}
        \begin{minipage}{1.5cm}
    \centering
    \includegraphics[height=1.1cm]{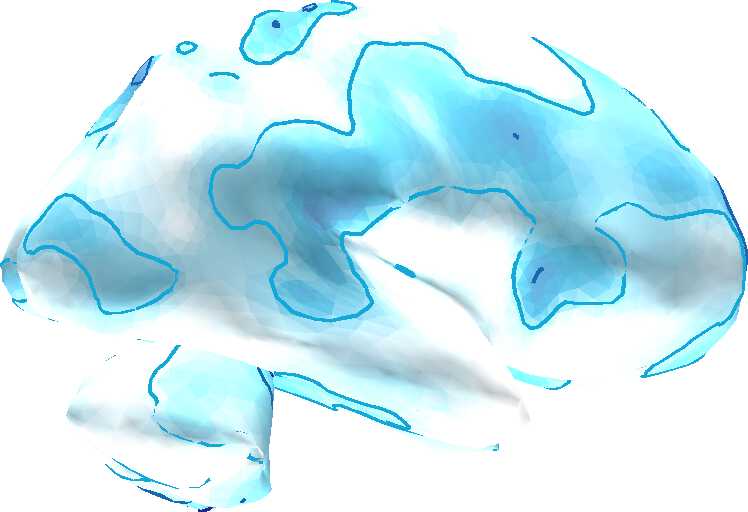}
    \end{minipage}
    \begin{minipage}{1.5cm}
    \centering
    \includegraphics[height=1.4cm]{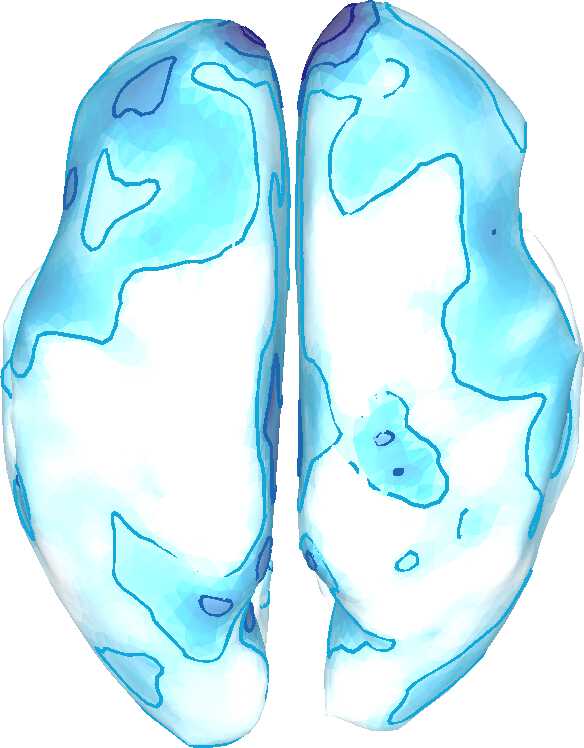}
    \end{minipage} \\
    \mbox{} \vskip0.1cm
    \begin{minipage}{0.20cm}
    \centering 
    \rotatebox{90}{\mbox{} \hskip0.01cm Primal-simplex} 
    \end{minipage}
    \begin{minipage}{1.5cm}
    \centering
    \includegraphics[height=1.1cm]{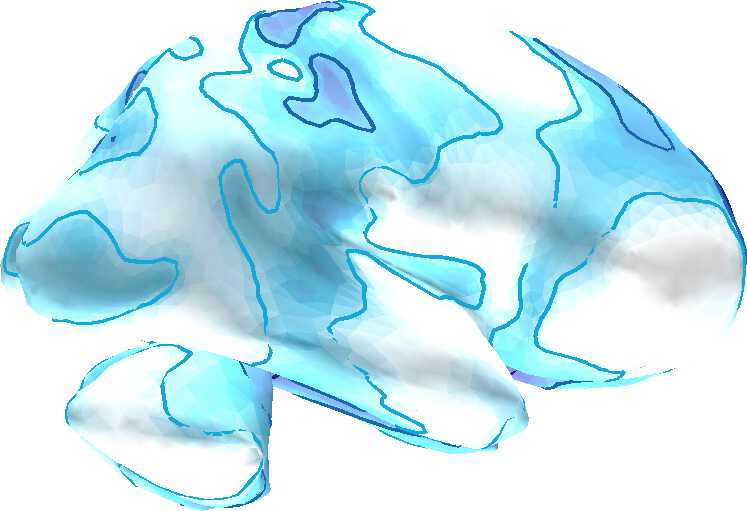}
    \end{minipage}
    \begin{minipage}{1.5cm}
    \centering
    \includegraphics[height=1.4cm]{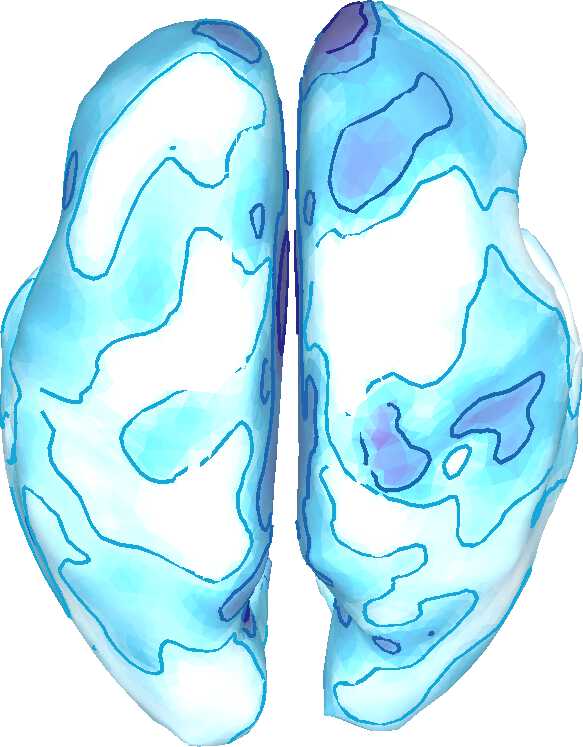}
    \end{minipage} \\ \vskip0.1cm
    Maximum current
    \end{minipage} 
    \begin{minipage}{0.20cm} 
    \includegraphics[height=3cm]{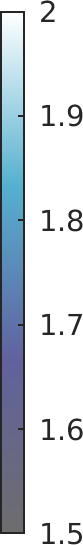}
    \end{minipage} 
    \begin{minipage}{4.0cm}
    \centering
    \begin{minipage}{0.20cm}
    \centering 
    \rotatebox{90}{\mbox{} \hskip0.01cm Interior-point} 
    \end{minipage}
    \begin{minipage}{1.5cm}
    \centering
    \includegraphics[height=1.1cm]{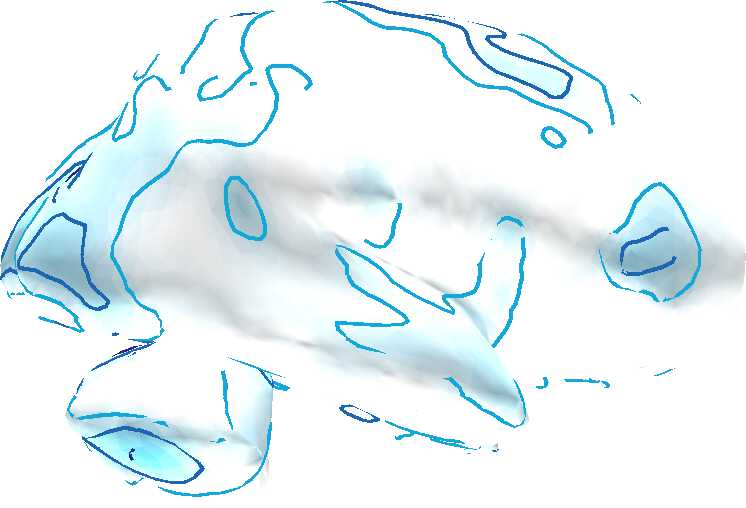}
    \end{minipage}
    \begin{minipage}{1.5cm}
    \centering
    \includegraphics[height=1.4cm]{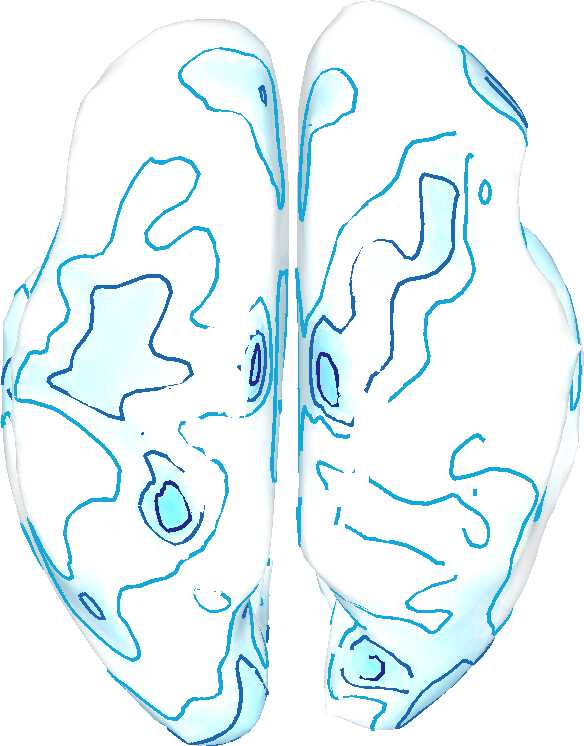}
    \end{minipage} \\ \mbox{} \vskip0.1cm
    \begin{minipage}{0.20cm}
    \centering 
    \rotatebox{90}{\mbox{} \hskip0.01cm Dual-simplex} 
    \end{minipage}
    \begin{minipage}{1.5cm}
    \centering
    \includegraphics[height=1.1cm]{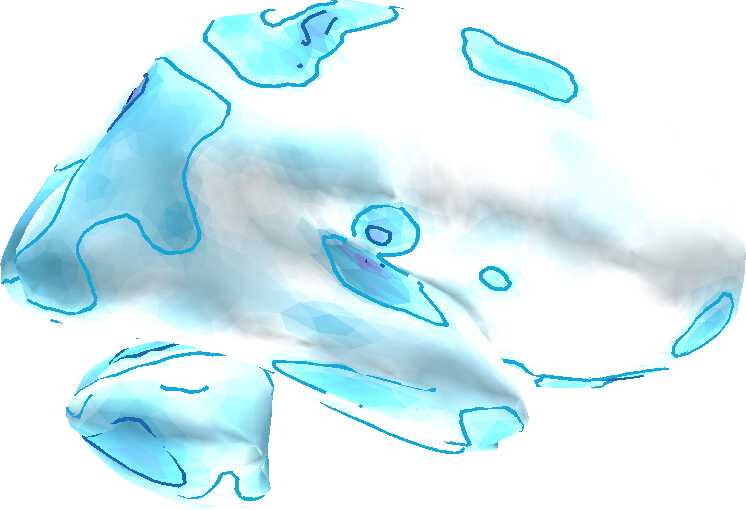}
    \end{minipage}
    \begin{minipage}{1.5cm}
    \centering
    \includegraphics[height=1.4cm]{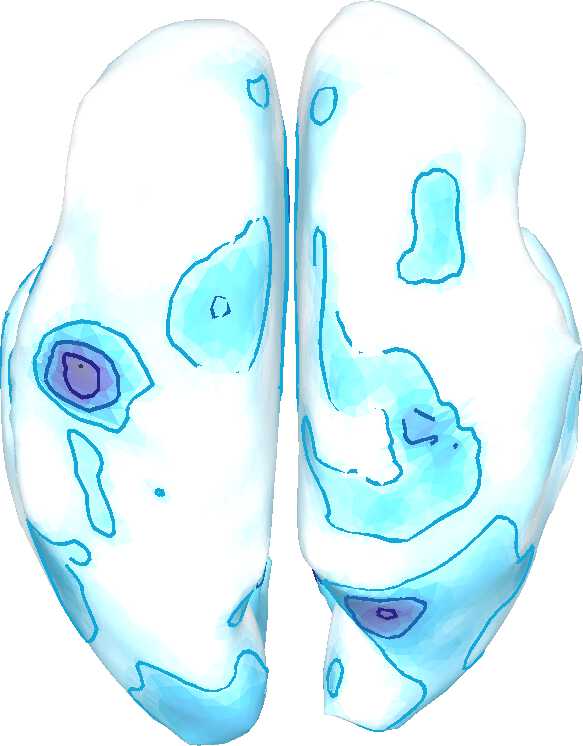}
    \end{minipage} \\
    \mbox{} \vskip0.1cm
    \begin{minipage}{0.20cm}
    \centering 
    \rotatebox{90}{\mbox{} \hskip0.01cm Primal-simplex} 
    \end{minipage}
    \begin{minipage}{1.5cm}
    \centering
    \includegraphics[height=1.1cm]{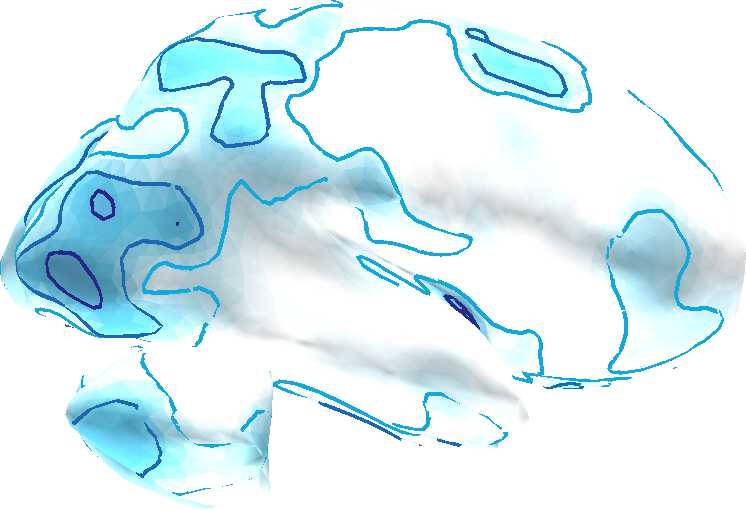}
    \end{minipage}
    \begin{minipage}{1.5cm}
    \centering
    \includegraphics[height=1.4cm]{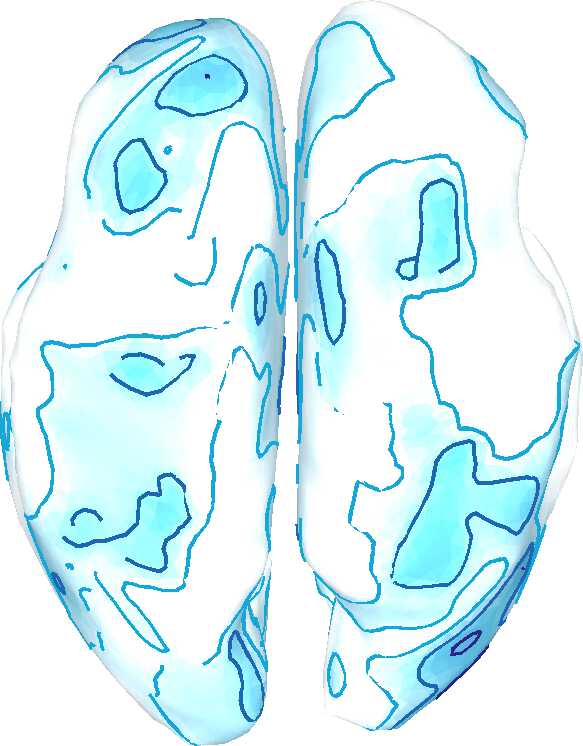}
    \end{minipage} \\ \vskip0.1cm
    NNZ
    \end{minipage} 
    \begin{minipage}{0.20cm} 
    \includegraphics[height=3cm]{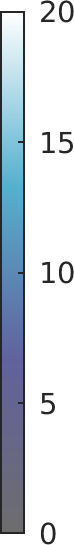}
    \end{minipage} \\ \vskip0.5cm
    The process to estimate $\Theta_{\hbox{\scriptsize max}}$
    \end{framed}
    \vskip0.3cm
    \label{fig:my_label_max_sup_nnz}
    \begin{framed}
    \begin{minipage}{4.0cm}
    \centering
    \begin{minipage}{0.20cm}
    \centering 
    \rotatebox{90}{\mbox{} \hskip0.01cm Interior-point} 
    \end{minipage}
    \begin{minipage}{1.5cm}
    \centering
    \includegraphics[height=1.1cm]{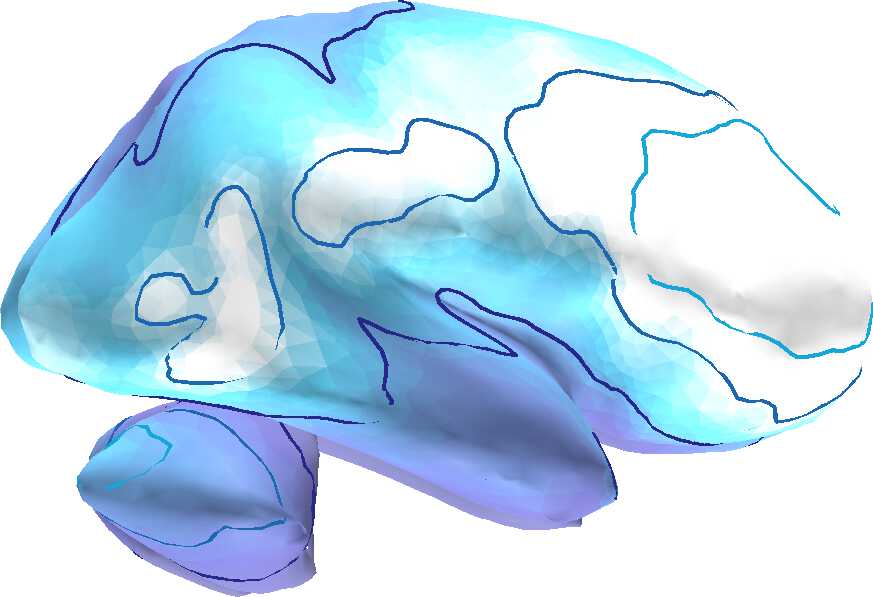}
    \end{minipage}
    \begin{minipage}{1.5cm}
    \centering
    \includegraphics[height=1.4cm]{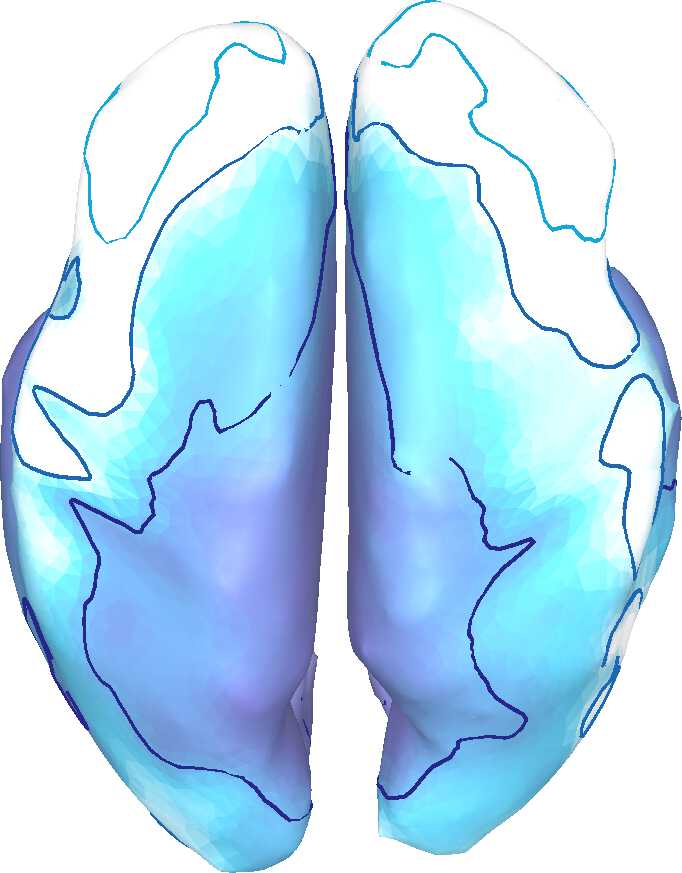}
    \end{minipage} \\ \mbox{} \vskip0.1cm
          \begin{minipage}{0.20cm}
    \centering 
    \rotatebox{90}{\mbox{} \hskip0.01cm Dual-simplex} 
    \end{minipage}
    \begin{minipage}{1.5cm}
    \centering
    \includegraphics[height=1.1cm]{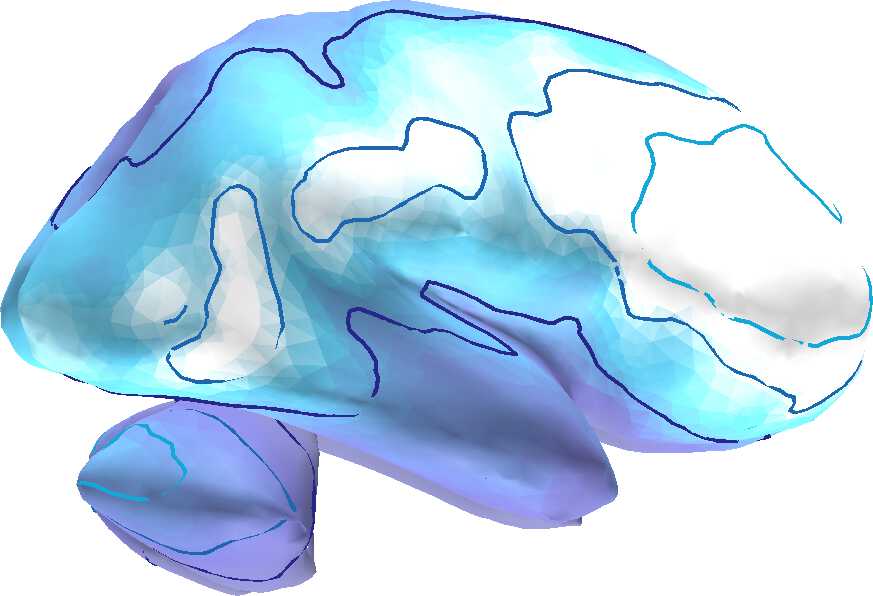}
    \end{minipage}
    \begin{minipage}{1.5cm}
    \centering
    \includegraphics[height=1.4cm]{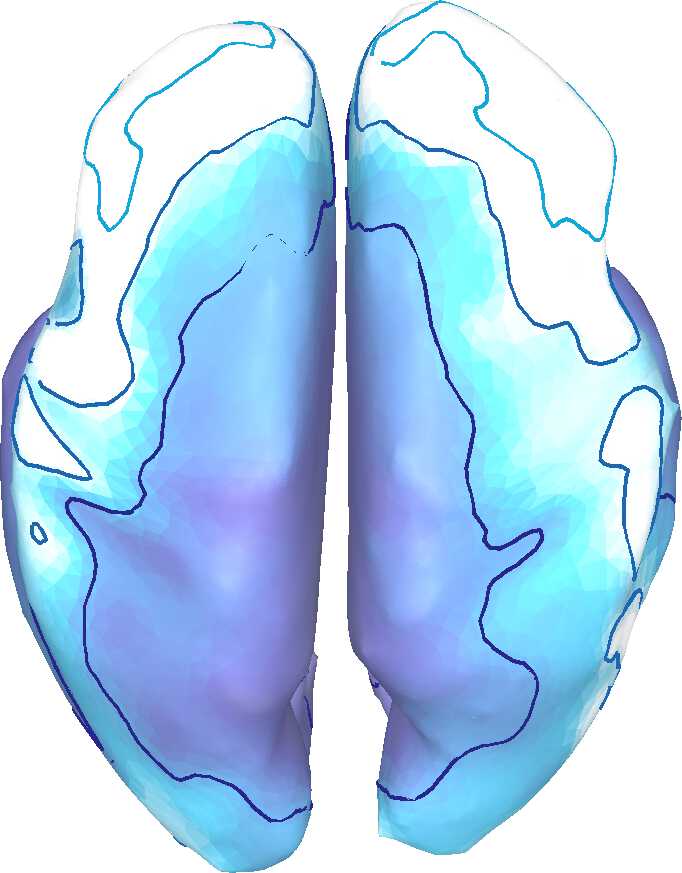}
    \end{minipage} \\
    \mbox{} \vskip0.1cm
    \begin{minipage}{0.20cm}
    \centering 
    \rotatebox{90}{\mbox{} \hskip0.01cm Primal-simplex} 
    \end{minipage}
    \begin{minipage}{1.5cm}
    \centering
    \includegraphics[height=1.1cm]{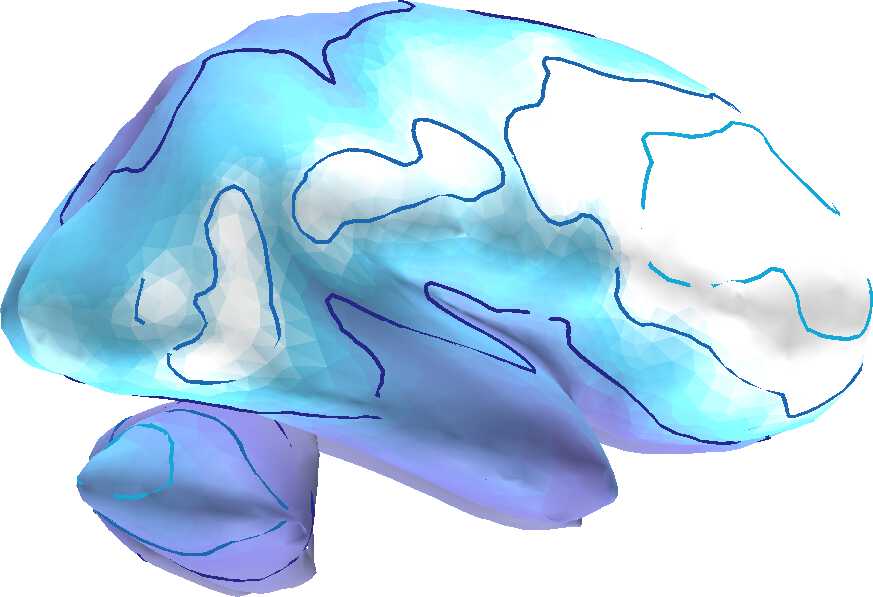}
    \end{minipage}
    \begin{minipage}{1.5cm}
    \centering
    \includegraphics[height=1.4cm]{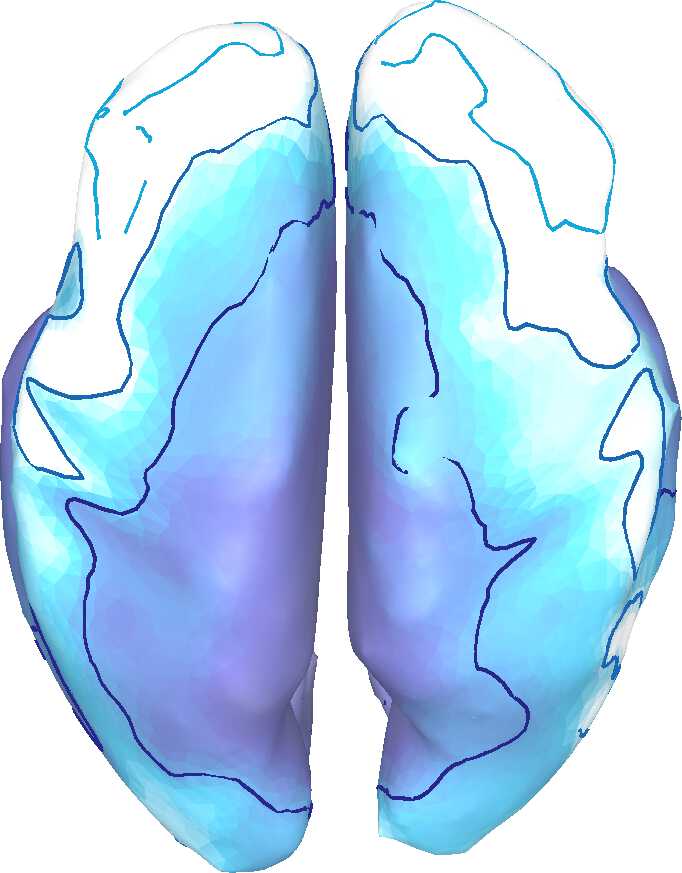}
    \end{minipage} \\ \vskip0.1cm
    $\Gamma_{\hbox{\scriptsize max}}$
    \end{minipage} 
    \begin{minipage}{0.20cm} 
    \includegraphics[height=3cm]{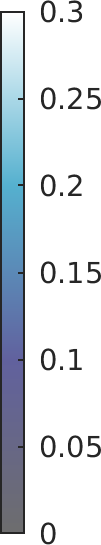}
    \end{minipage} 
    \label{fig:my_label_B_Gamma}
    \centering
    \begin{minipage}{4.0cm}
    \centering
    \begin{minipage}{0.20cm}
    \centering 
     \rotatebox{90}{\mbox{} \hskip0.01cm Interior-point} 
    \end{minipage}
    \begin{minipage}{1.5cm}
    \centering
    \includegraphics[height=1.1cm]{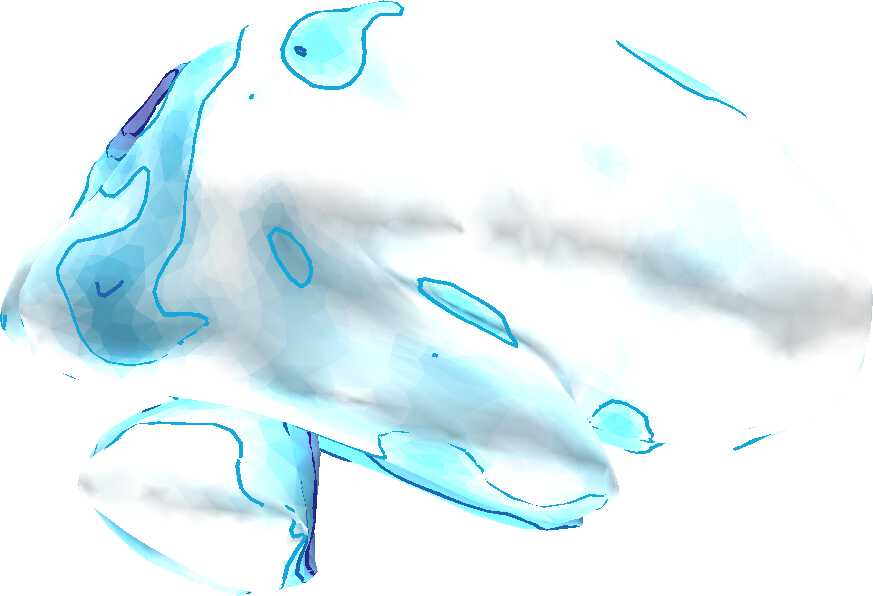}
    \end{minipage}
    \begin{minipage}{1.5cm}
    \centering
    \includegraphics[height=1.4cm]{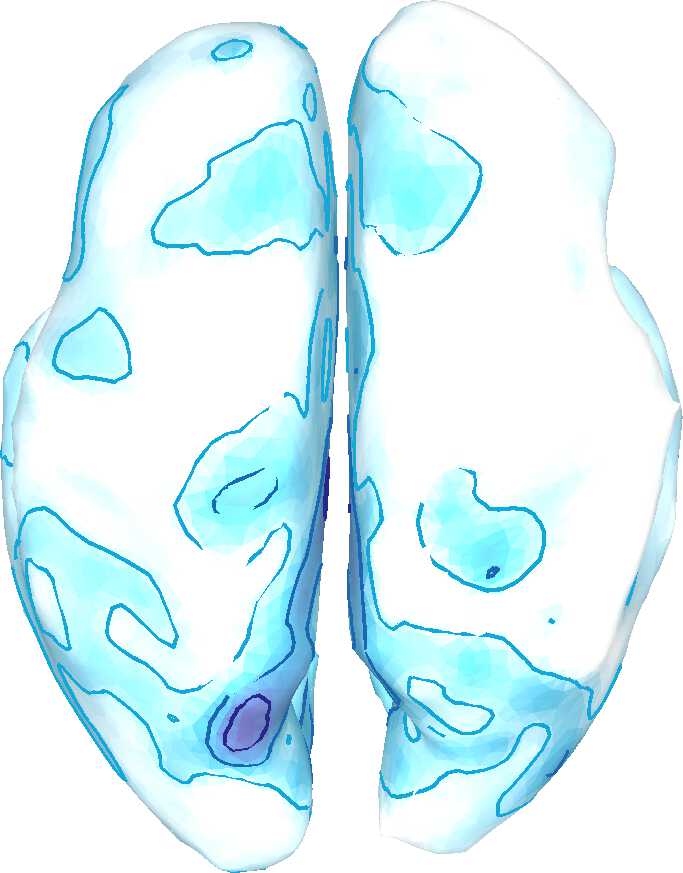}
    \end{minipage} \\ \mbox{} \vskip0.1cm
    \begin{minipage}{0.20cm}
    \centering 
     \rotatebox{90}{\mbox{} \hskip0.01cm Dual-simplex} 
    \end{minipage}
    \begin{minipage}{1.5cm}
    \centering
    \includegraphics[height=1.1cm]{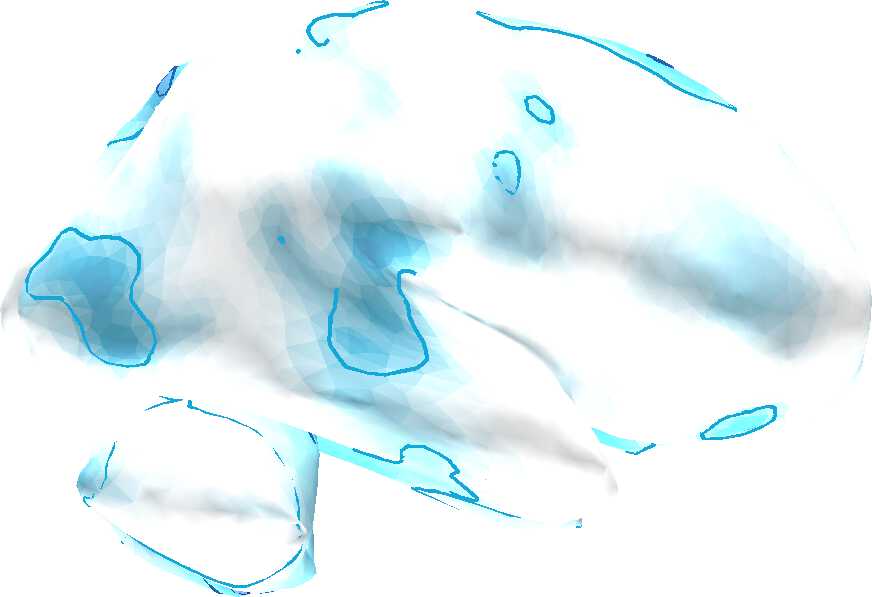}
    \end{minipage}
    \begin{minipage}{1.5cm}
    \centering
    \includegraphics[height=1.4cm]{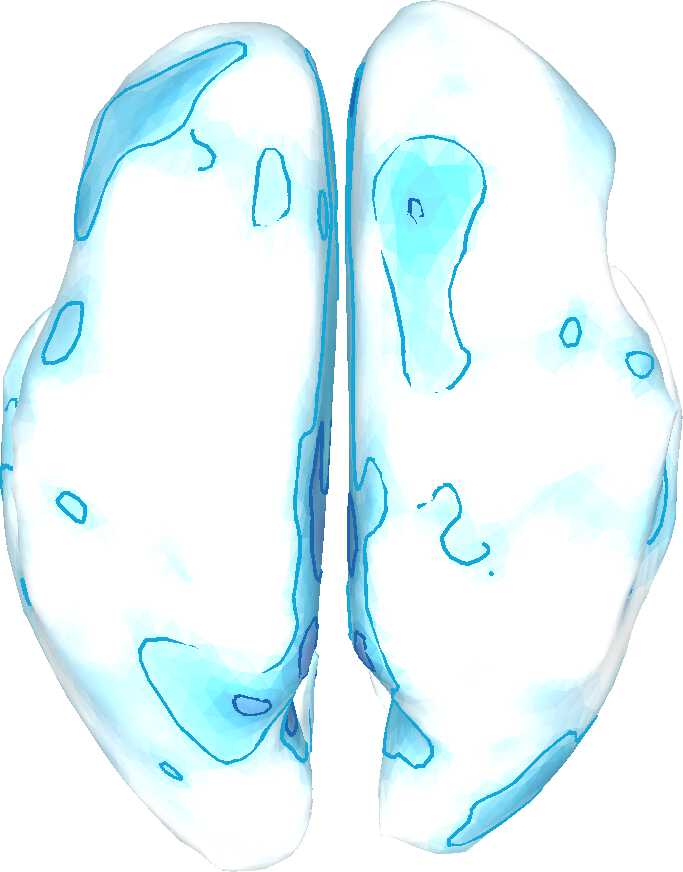}
    \end{minipage} \\
    \mbox{} \vskip0.1cm
    \begin{minipage}{0.20cm}
    \centering 
    \rotatebox{90}{\mbox{} \hskip0.01cm Primal-simplex} 
    \end{minipage}
    \begin{minipage}{1.5cm}
    \centering
    \includegraphics[height=1.1cm]{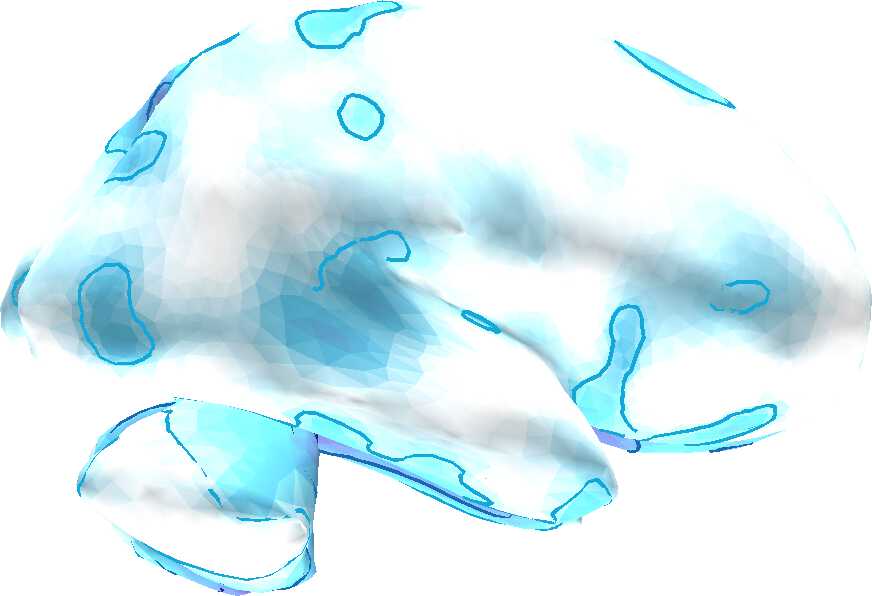}
    \end{minipage}
    \begin{minipage}{1.5cm}
    \centering
    \includegraphics[height=1.4cm]{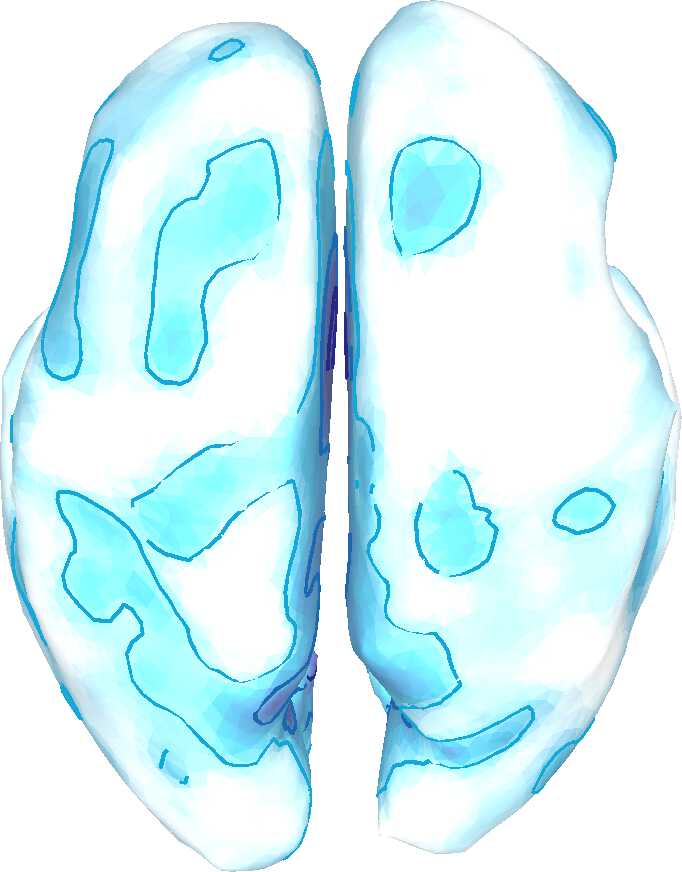}
    \end{minipage} \\ \vskip0.1cm
    Maximum current
    \end{minipage} 
    \begin{minipage}{0.20cm} 
    \includegraphics[height=3cm]{bar_3.png}
    \end{minipage} 
    \begin{minipage}{4.0cm}
    \centering
    \begin{minipage}{0.20cm}
    \centering 
    \rotatebox{90}{\mbox{} \hskip0.01cm Interior-point} 
    \end{minipage}
    \begin{minipage}{1.5cm}
    \centering
    \includegraphics[height=1.1cm]{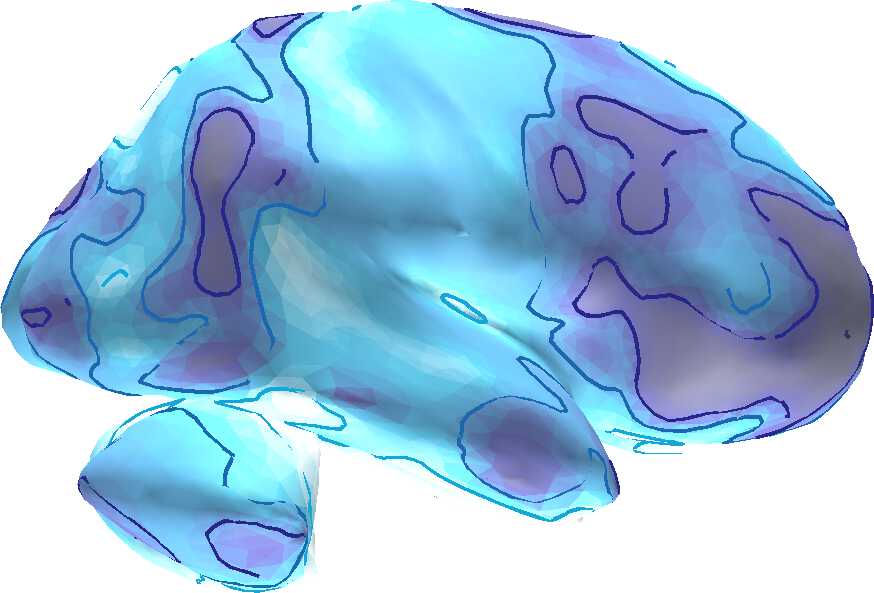}
    \end{minipage}
    \begin{minipage}{1.5cm}
    \centering
    \includegraphics[height=1.4cm]{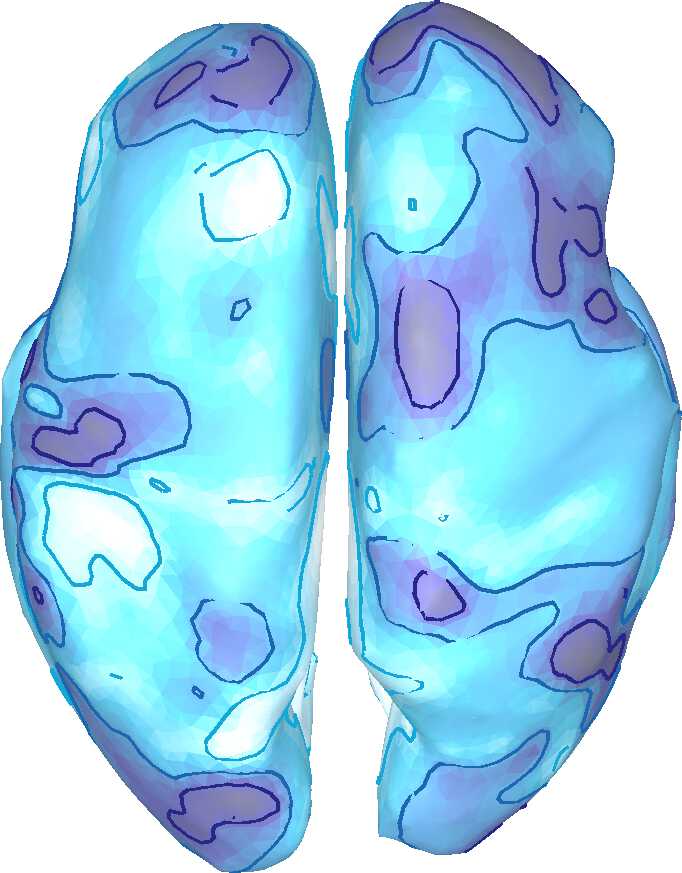}
    \end{minipage} \\ \mbox{} \vskip0.1cm
    \begin{minipage}{0.20cm}
    \centering 
    \rotatebox{90}{\mbox{} \hskip0.01cm Dual-simplex} 
    \end{minipage}
        \begin{minipage}{1.5cm}
    \centering
    \includegraphics[height=1.1cm]{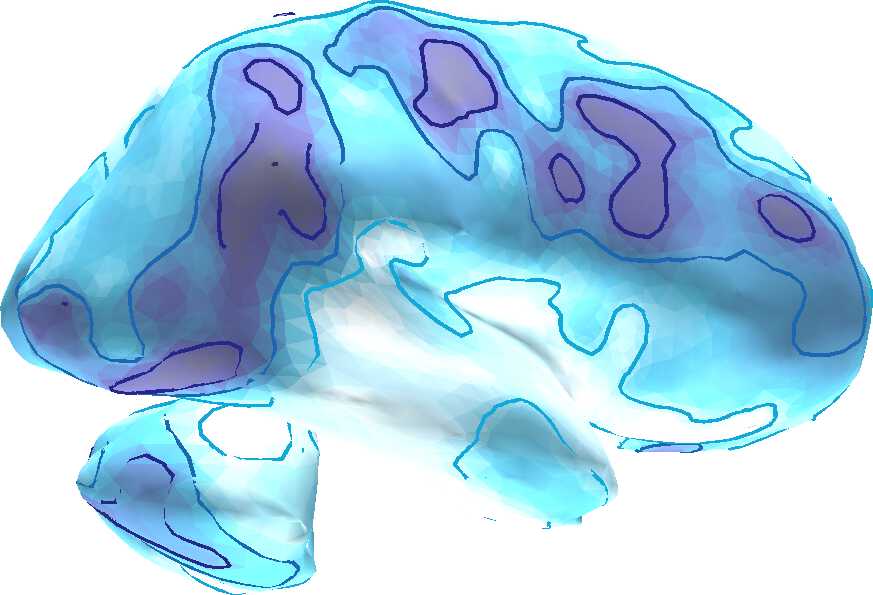}
    \end{minipage}
    \begin{minipage}{1.5cm}
    \centering
    \includegraphics[height=1.4cm]{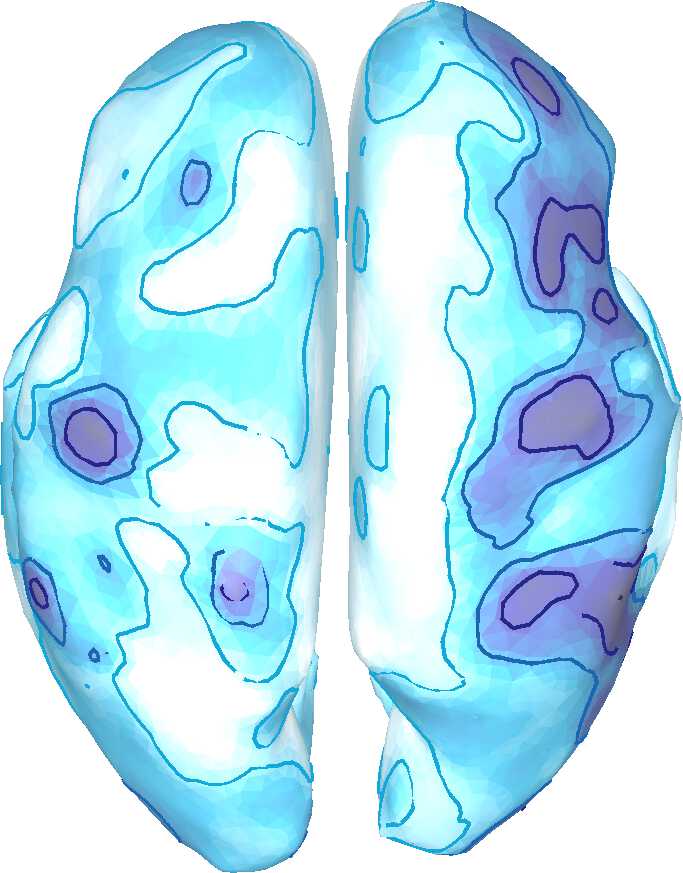}
    \end{minipage} \\
    \mbox{} \vskip0.1cm
    \begin{minipage}{0.20cm}
    \centering 
    \rotatebox{90}{\mbox{} \hskip0.01cm Primal-simplex} 
    \end{minipage}
    \begin{minipage}{1.5cm}
    \centering
    \includegraphics[height=1.1cm]{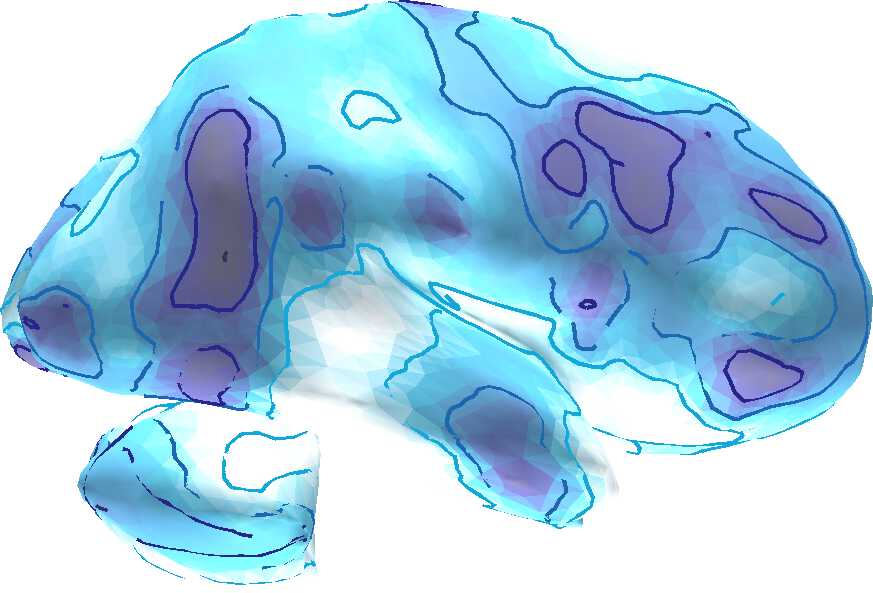}
    \end{minipage}
    \begin{minipage}{1.5cm}
    \centering
    \includegraphics[height=1.4cm]{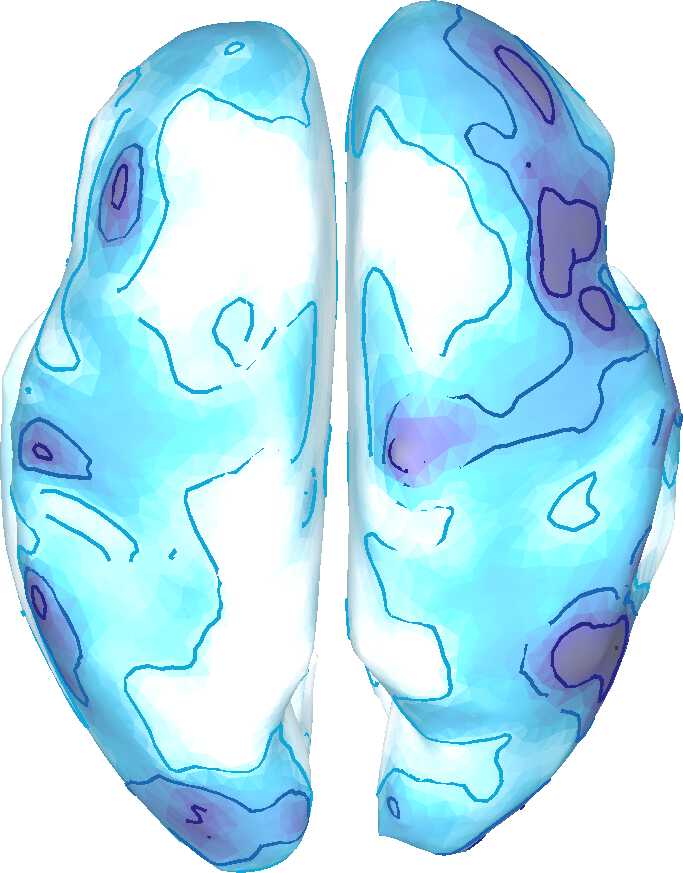}
    \end{minipage} \\ \vskip0.1cm
    NNZ
    \end{minipage} 
    \begin{minipage}{0.20cm} 
    \includegraphics[height=3cm]{bar_5.png}
    \end{minipage} \\ \vskip0.5cm 
    The process to estimate $\Gamma_{\hbox{\scriptsize max}}$
    \end{framed}
    \vskip0.1cm
    \end{scriptsize}
    \caption{{\bf Top:} Maps of the maximum focality  $\Theta_{\hbox{\scriptsize max}}$ and the values of intensity $\Gamma$, the maximum current and NNZ corresponding to the maximizer for coarser $15 \times 15$ lattice and the interior-point (IP), dual-simplex (DS) and primal-simplex (PS) algorithm of the MOSEK package. Each map includes countours for 20, 40 and 75 \% of the maximum entry. While the maps obtained with different algorithms are topographically similar, IP can be observed to enhance $\Theta_{\hbox{\scriptsize max}}$, especially, in the areas, where the solutions are likely to be focal. In these areas, DS and PS yield a greater intensity which comes with the cost of suppressed focality. {\bf Bottom:} Maps of the maximum intensity $\Gamma_{\hbox{\scriptsize max}}$. The mutual differences between the maps are minor, both topographically and amplitude-wise. Of  $\Theta_{\hbox{\scriptsize max}}$ and $\Gamma_{\hbox{\scriptsize max}}$,   the former one corresponds to a greater NNZ and a smaller maximum current.}
    \label{fig:my_label_ratio_theta}
    \label{fig:my_label_NNZ}
\end{figure*}

\subsection{Solver packages}
%We are dealing with two different classes and even two sub classes. Under the simplex class there has two sub classes which is DS and PS, and then Interior point class includes the parallel barer by Gurobi’s PS. So we take 11 different methods to measure the computing time on three different sources, occipital , postcentral and auditory.
\section{Results}

%***************************************************
%***************************************************

\begin{figure*}
    \centering
    \begin{scriptsize}
    \begin{minipage}{8.8cm}
    \centering
\begin{minipage}{0.5cm}
\centering
\mbox{}
\end{minipage}
    \begin{minipage}{2.4cm}
    \centering
   Interior-point
    \end{minipage}
       \begin{minipage}{2.4cm}
    \centering
   Dual-sumplex
    \end{minipage}
           \begin{minipage}{2.4cm}
    \centering
   Primal-simplex
    \end{minipage}
               \\ \vskip0.1cm
\begin{minipage}{0.5cm}
\centering
 \rotatebox{90}{\mbox{} \hskip0.01cm Matlab} 
\end{minipage}
    \begin{minipage}{2.4cm}
    \centering
    \includegraphics[width=2.3cm]{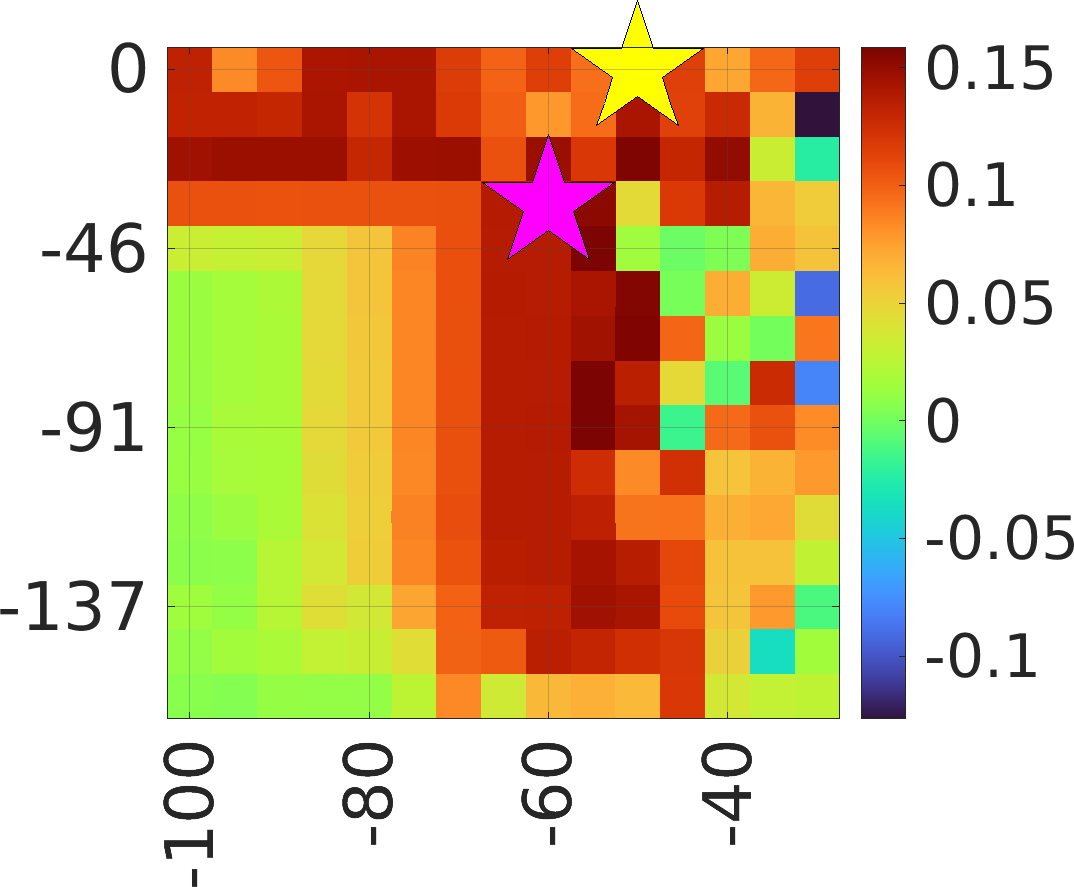}
    \end{minipage}
       \begin{minipage}{2.4cm}
    \centering
        \includegraphics[width=2.3cm]{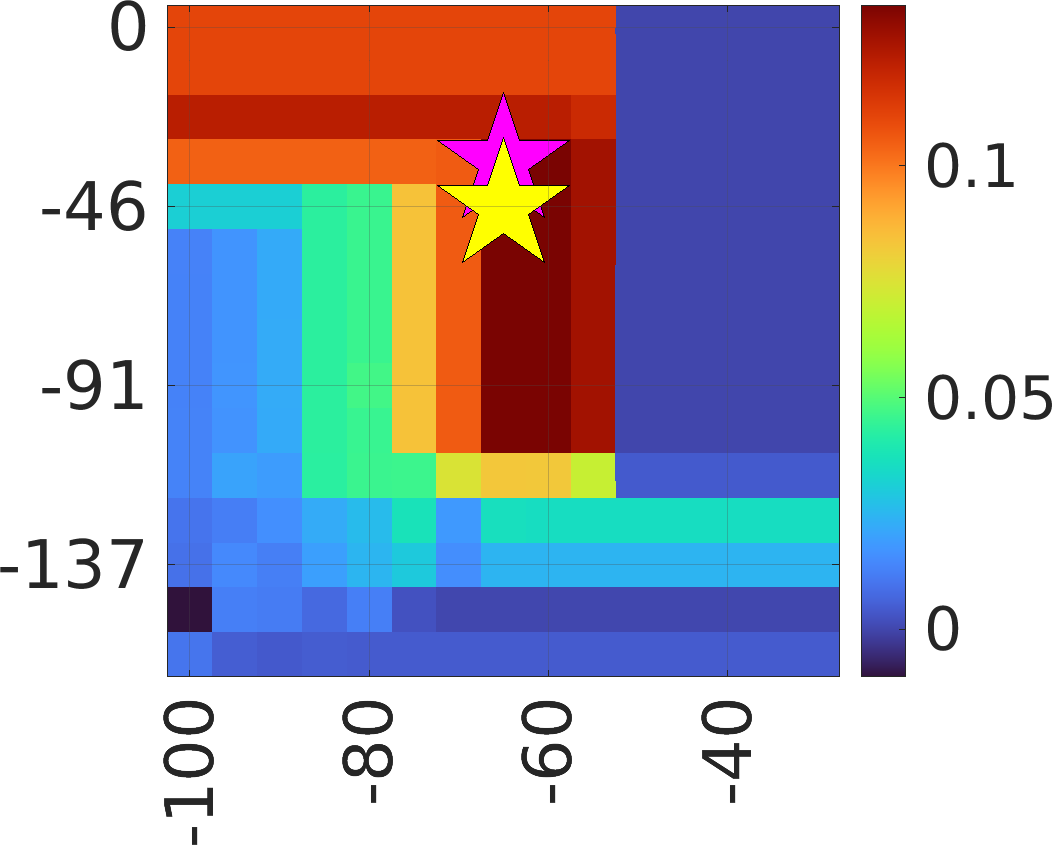}
    \end{minipage}
           \begin{minipage}{2.4cm}
    \centering
\mbox{}
    \end{minipage}
    \\ \vskip0.1cm
    \begin{minipage}{0.5cm}
\centering
 \rotatebox{90}{\mbox{} \hskip0.01cm MOSEK} 
\end{minipage}
               \begin{minipage}{2.4cm}
    \centering
    \includegraphics[width=2.3cm]{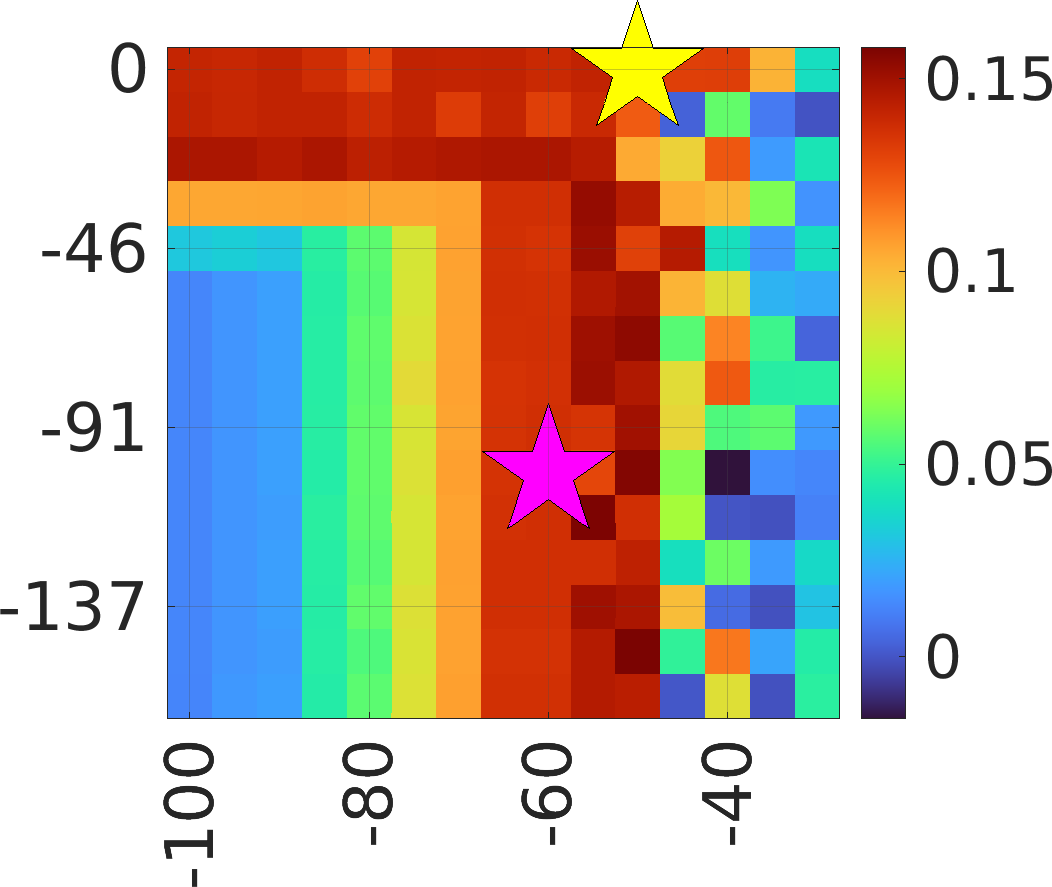}
    \end{minipage}
               \begin{minipage}{2.4cm}
    \centering
    \includegraphics[width=2.3cm]{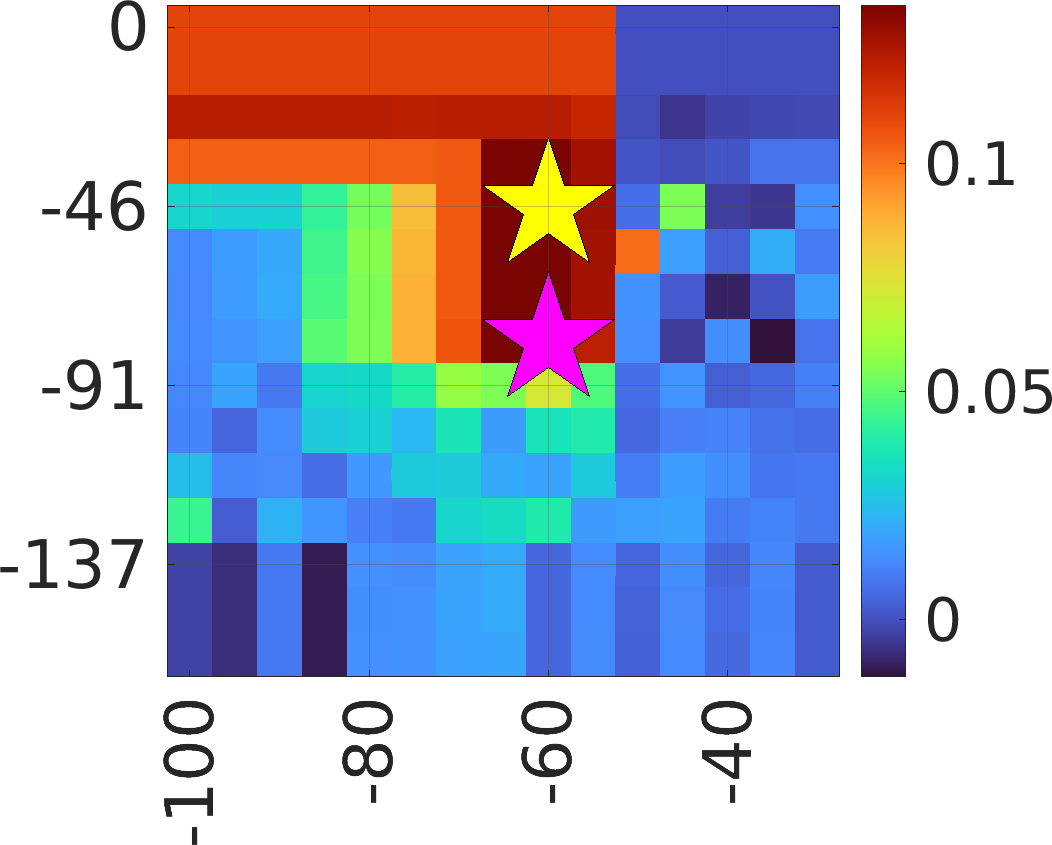}
    \end{minipage} 
        \begin{minipage}{2.4cm}
    \centering
    \includegraphics[width=2.3cm]{chart_3_1_2_1_2_2_1.png}
    \end{minipage}
    \\ \vskip0.1cm
    \begin{minipage}{0.5cm}
\centering
 \rotatebox{90}{\mbox{} \hskip0.01cm Gurobi} 
\end{minipage}
       \begin{minipage}{2.4cm}
    \centering
    \includegraphics[width=2.3cm]{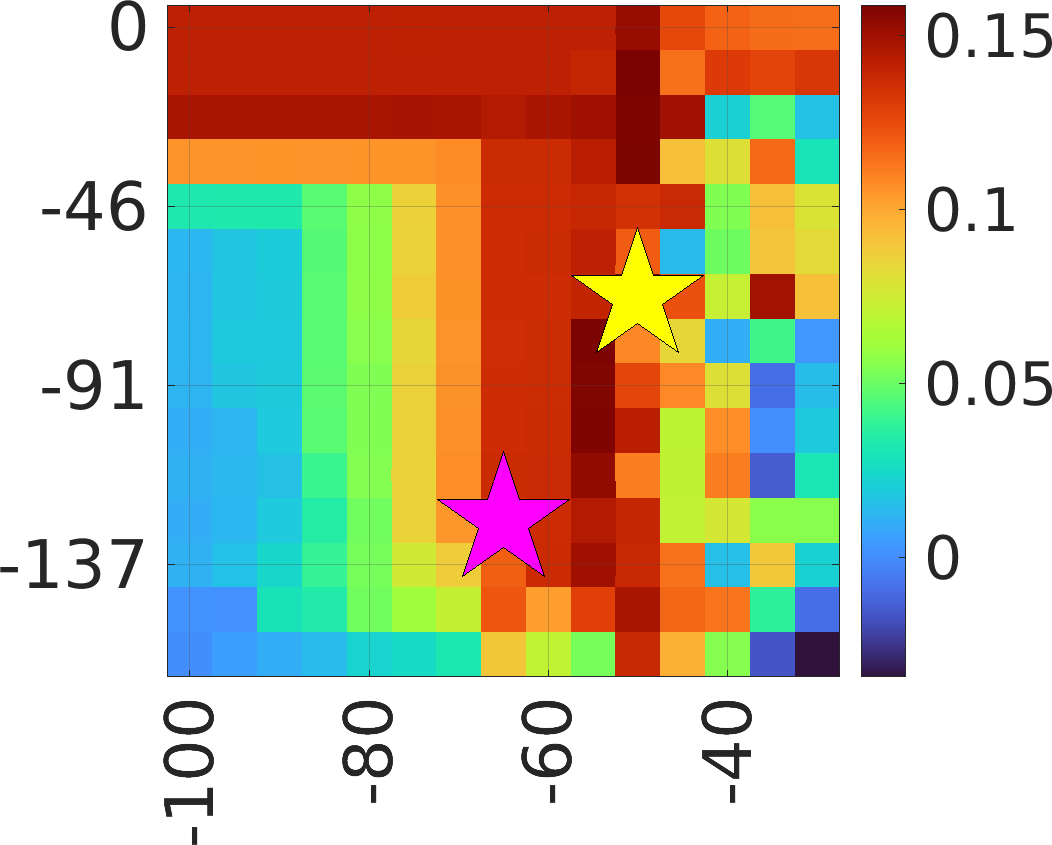}
    \end{minipage}
           \begin{minipage}{2.4cm}
    \centering
    \includegraphics[width=2.3cm]{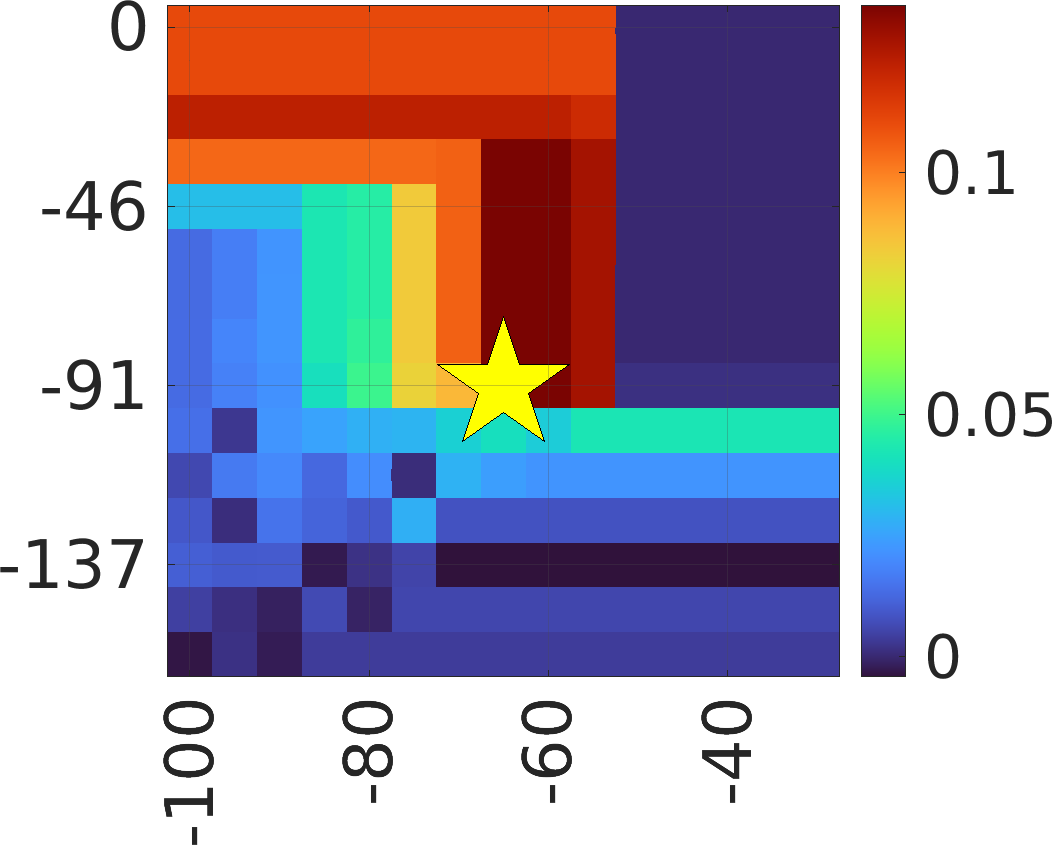}
    \end{minipage}
               \begin{minipage}{2.4cm}
    \centering
 \includegraphics[width=2.3cm]{chart_3_1_3_1_2_2_1.png}
    \end{minipage}
    \\ \vskip0.1cm
    \begin{minipage}{0.5cm}
\centering
 \rotatebox{90}{\mbox{} \hskip0.01cm SDPT3} 
\end{minipage}
               \begin{minipage}{2.4cm}
    \centering
   \includegraphics[width=2.3cm]{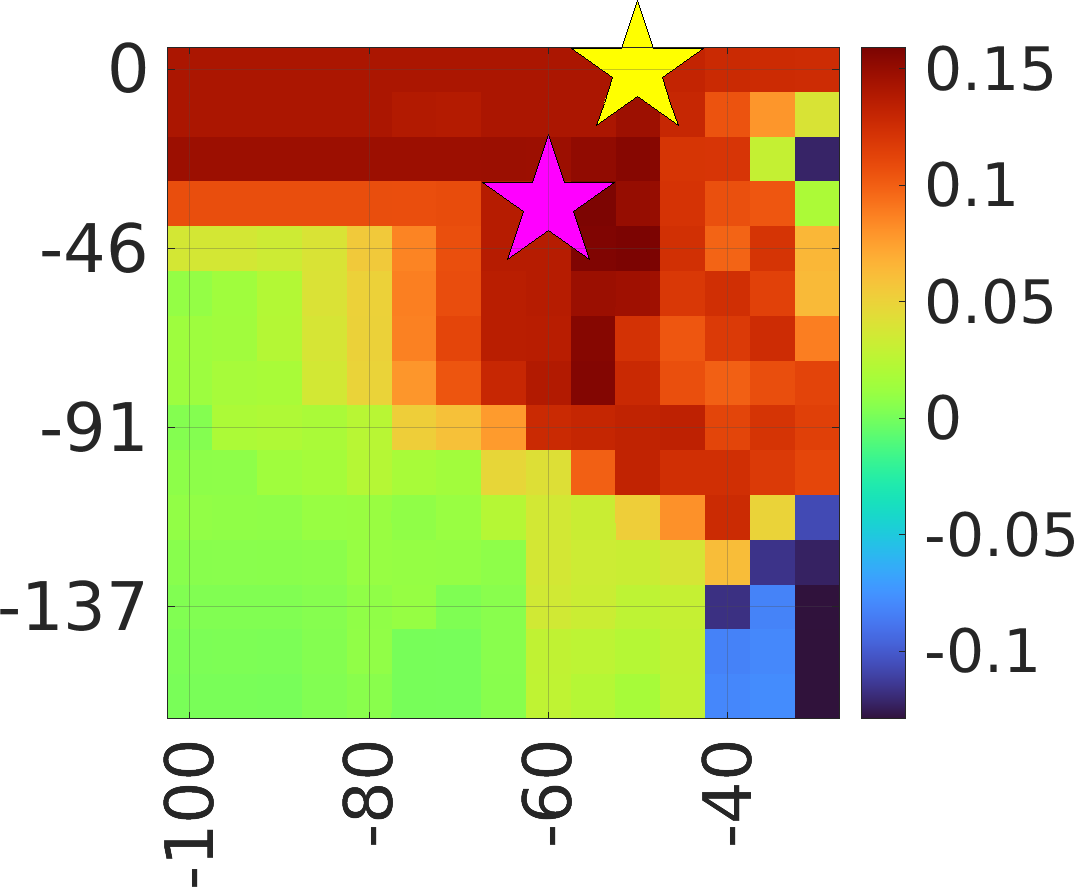}
    \end{minipage}
            \begin{minipage}{2.4cm}
    \centering
    \mbox{}
    \end{minipage}
       \begin{minipage}{2.4cm}
    \centering
\mbox{}
    \end{minipage}
    \\ \vskip0.1cm 
    \begin{minipage}{0.5cm}
\centering
 \rotatebox{90}{\mbox{} \hskip0.01cm SeDuMi} 
\end{minipage}
           \begin{minipage}{2.4cm}
    \centering
     \includegraphics[width=2.3cm]{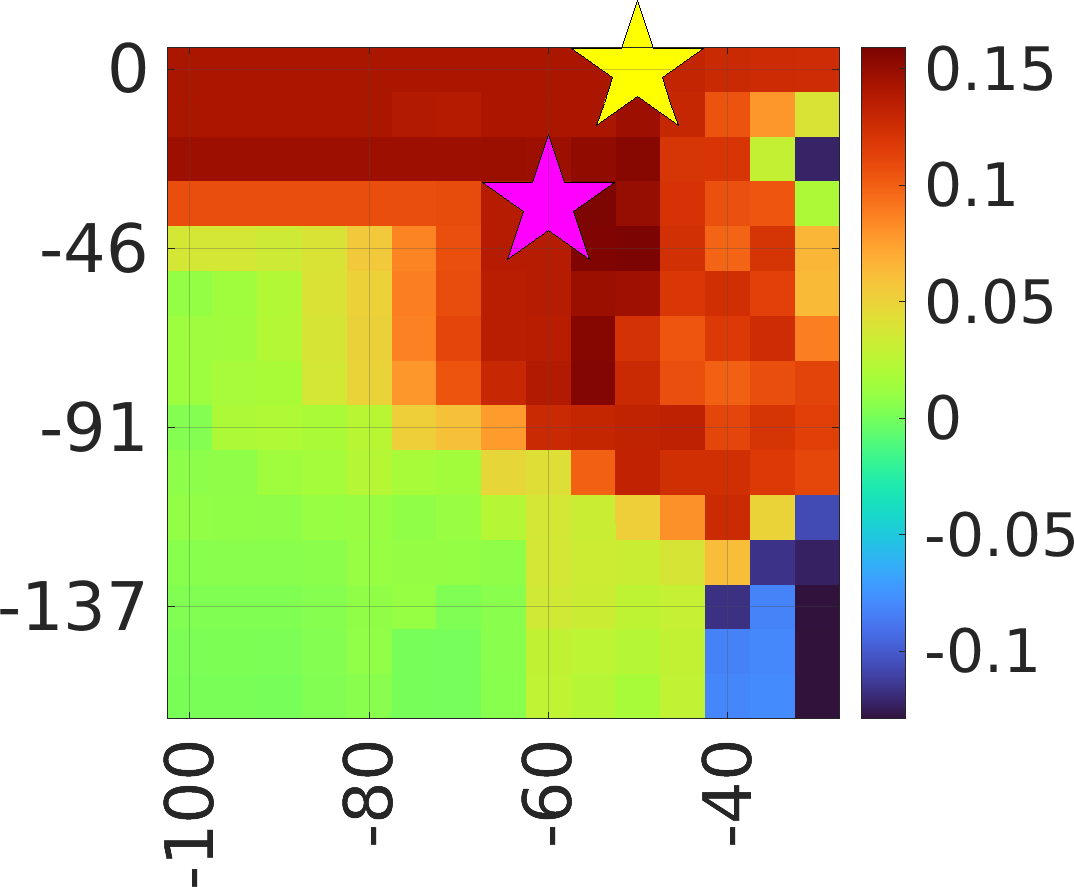}
    \end{minipage}
               \begin{minipage}{2.4cm}
    \centering
 \mbox{}
    \end{minipage}
               \begin{minipage}{2.4cm}
    \centering
    \mbox{}
    \end{minipage}
    \\
    \vskip0.1cm
    Coarse $15 \times 15$ lattice
    \end{minipage}
    \end{scriptsize}
    \centering
    \begin{scriptsize}
    \begin{minipage}{8.8cm}
    \centering
    \begin{minipage}{0.5cm}
    \centering
    \mbox{}
    \end{minipage}
    \begin{minipage}{2.4cm}
    \centering
    Interior-point
    \end{minipage}
    \begin{minipage}{2.4cm}
    \centering
    Dual-sumplex
    \end{minipage}
    \begin{minipage}{2.4cm}
    \centering
    Primal-simplex
    \end{minipage}
    \\
    \vskip0.1cm
    \begin{minipage}{0.5cm}
    \centering
    \rotatebox{90}{\mbox{} \hskip0.01cm Matlab} 
    \end{minipage}
    \begin{minipage}{2.4cm}
    \centering
    \includegraphics[width=2.3cm]{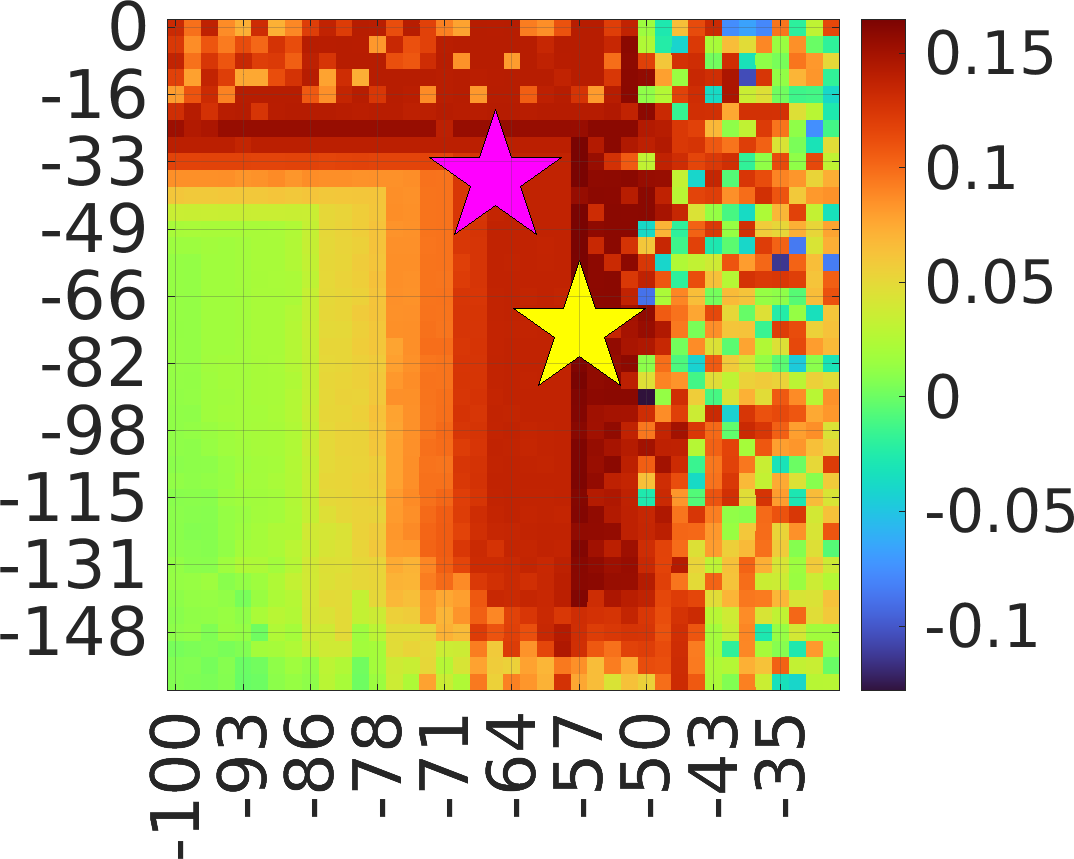}
    \end{minipage}
    \begin{minipage}{2.4cm}
    \centering
    \includegraphics[width=2.3cm]{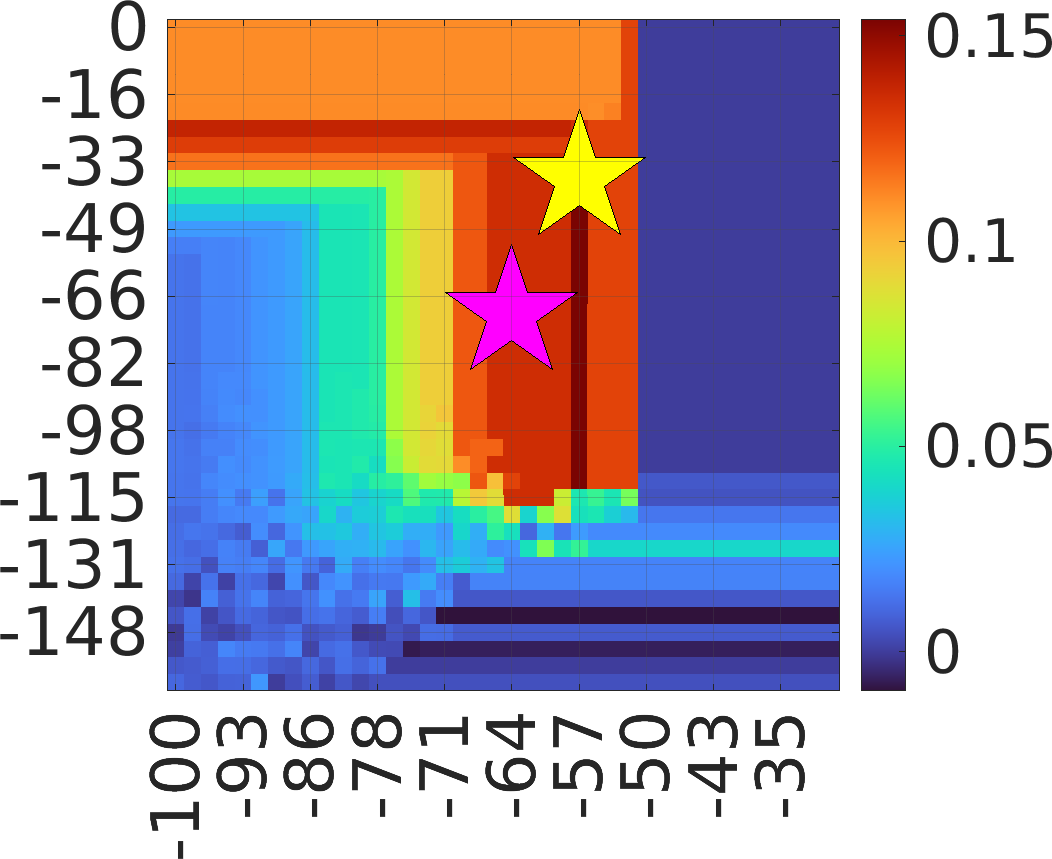}
    \end{minipage}
    \begin{minipage}{2.4cm}
    \centering
    \mbox{}
    \end{minipage}
    \\
    \vskip0.1cm
    \begin{minipage}{0.5cm}
    \centering
    \rotatebox{90}{\mbox{} \hskip0.01cm MOSEK} 
    \end{minipage}
    \begin{minipage}{2.4cm}
    \centering
    \includegraphics[width=2.3cm]{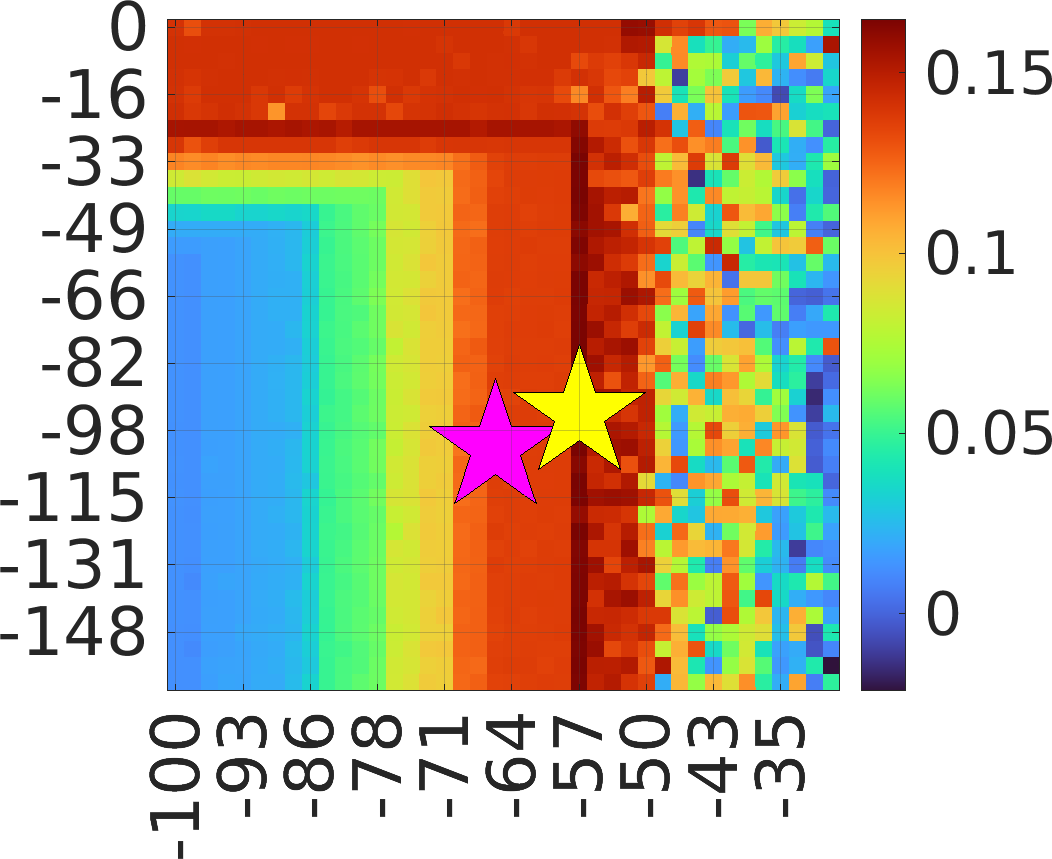}
    \end{minipage}
    \begin{minipage}{2.4cm}
    \centering
    \includegraphics[width=2.3cm]{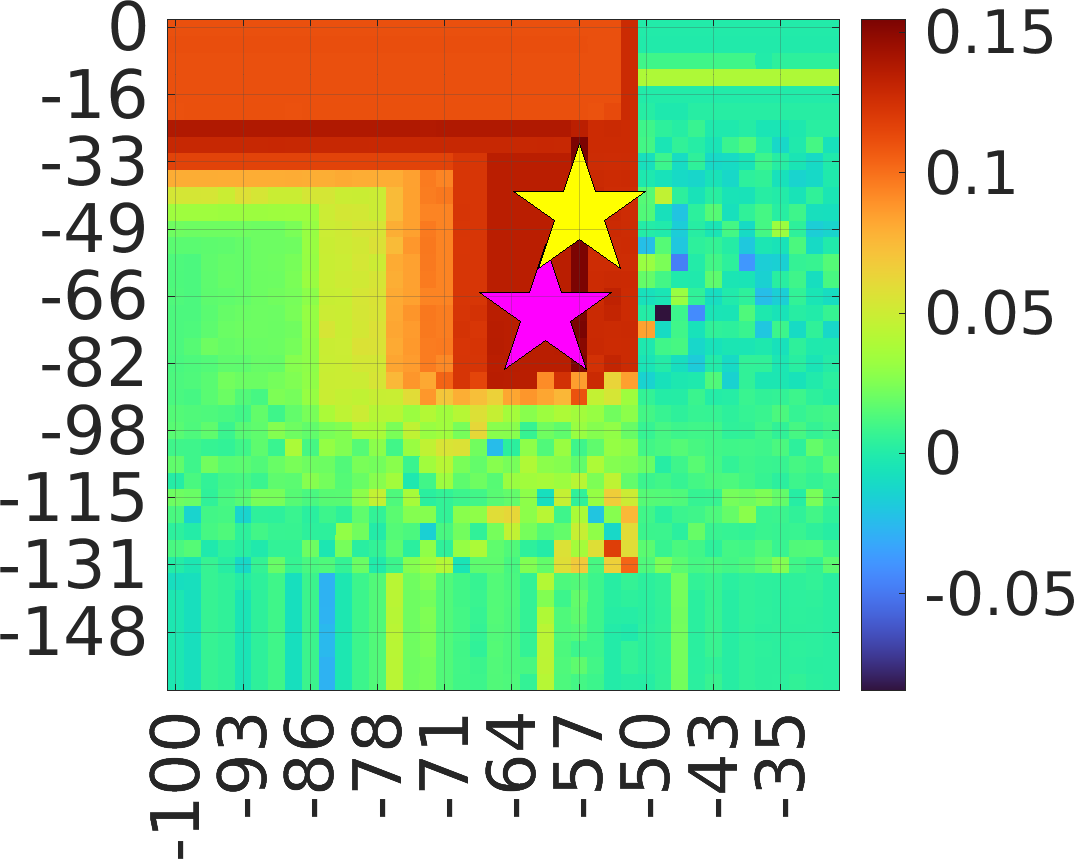}
    \end{minipage} 
    \begin{minipage}{2.4cm}
    \centering
    \includegraphics[width=2.3cm]{chart_3_2_2_1_2_2_1.png}
    \end{minipage}
    \\
    \vskip0.1cm
    \begin{minipage}{0.5cm}
    \centering
    \rotatebox{90}{\mbox{} \hskip0.01cm Gurobi} 
    \end{minipage}
    \begin{minipage}{2.4cm}
    \centering
    \includegraphics[width=2.3cm]{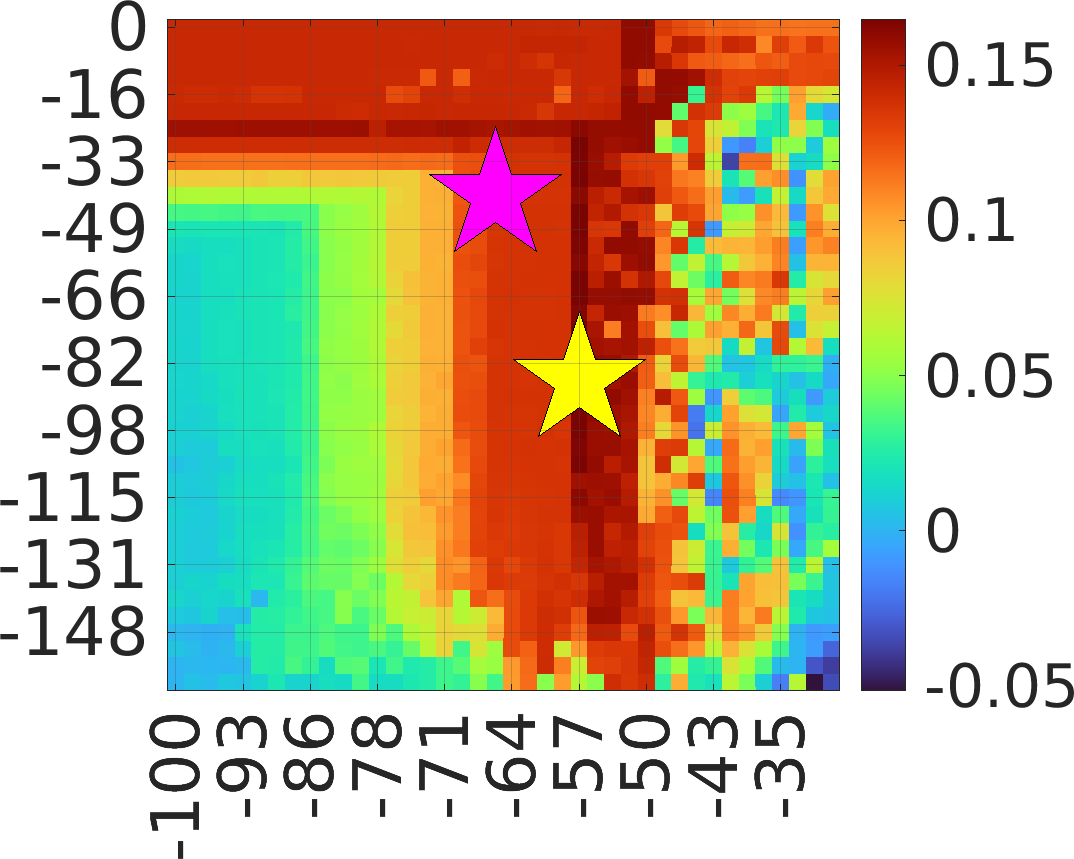}
    \end{minipage}
    \begin{minipage}{2.4cm}
    \centering
    \includegraphics[width=2.3cm]{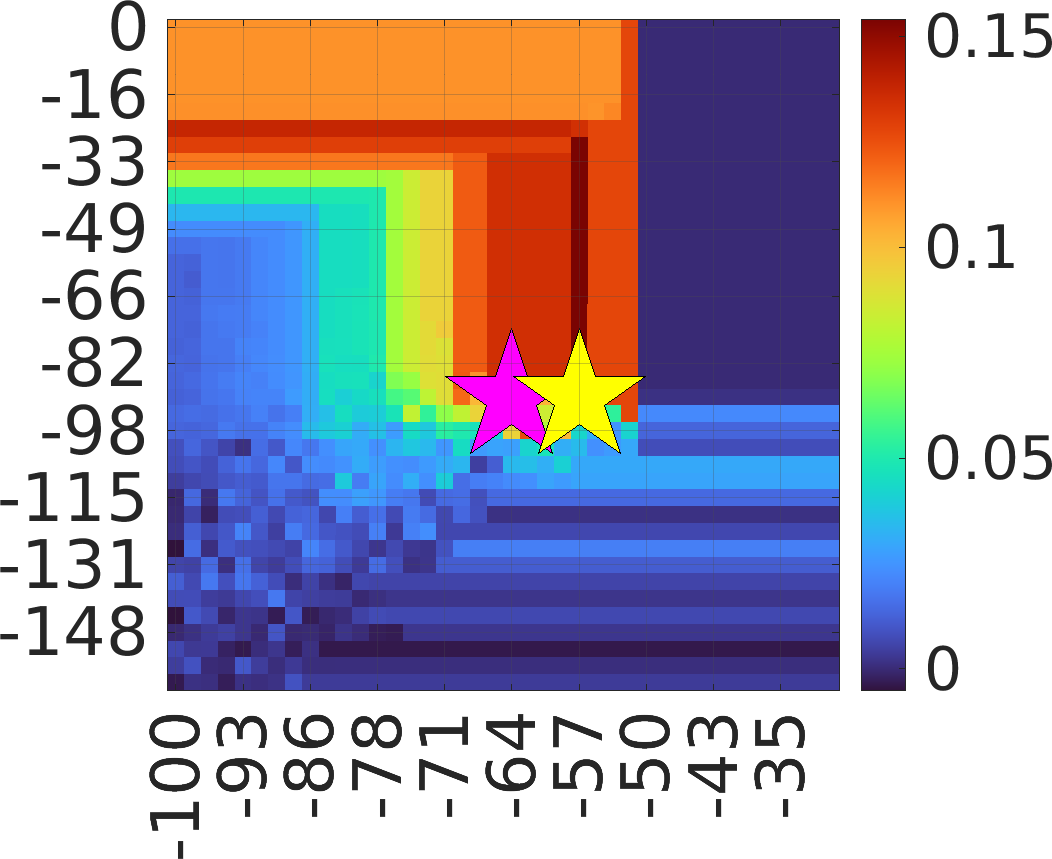}
    \end{minipage}
    \begin{minipage}{2.4cm}
    \centering
    \includegraphics[width=2.3cm]{chart_3_2_3_1_2_2_1.png}
    \end{minipage}
    \\
    \vskip0.1cm
    \begin{minipage}{0.5cm}
    \centering
    \rotatebox{90}{\mbox{} \hskip0.01cm SDPT3} 
    \end{minipage}
    \begin{minipage}{2.4cm}
    \centering
    \includegraphics[width=2.3cm]{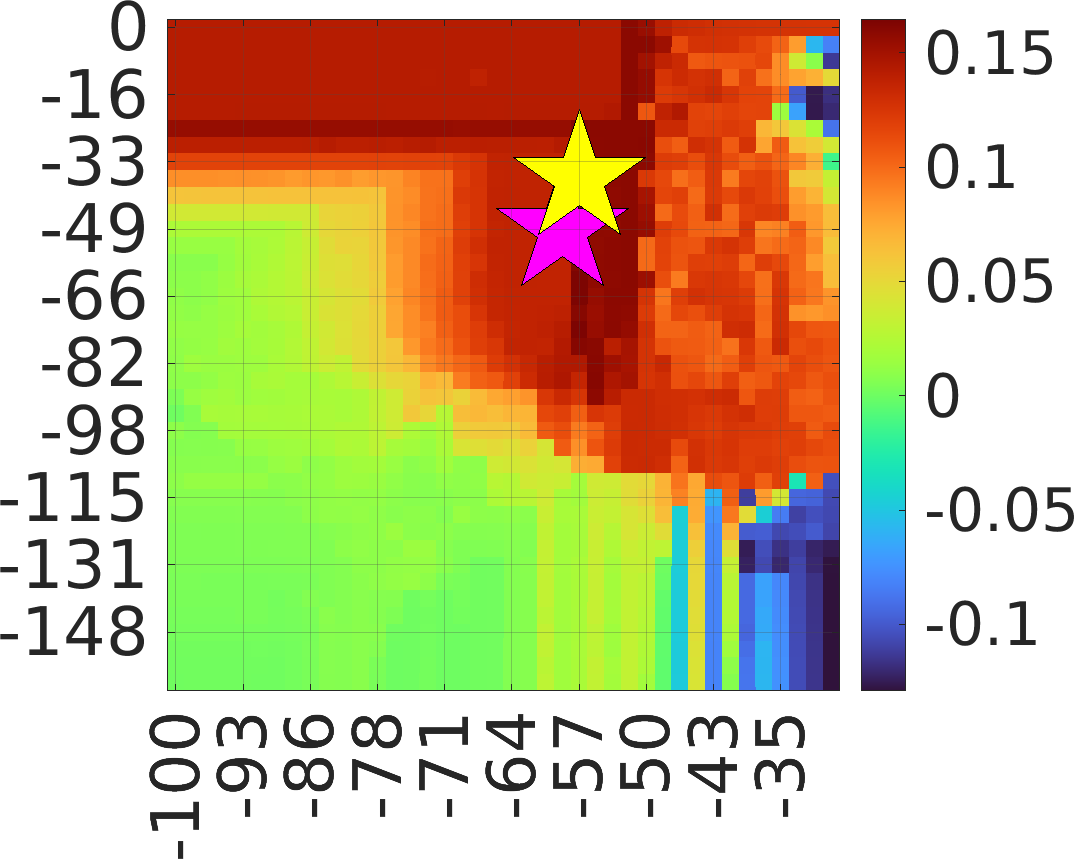}
    \end{minipage}
    \begin{minipage}{2.4cm}
    \centering
    \mbox{}
    \end{minipage}
    \begin{minipage}{2.4cm}
    \centering
    \mbox{}
    \end{minipage}
    \\
    \vskip0.1cm 
    \begin{minipage}{0.5cm}
    \centering
    \rotatebox{90}{\mbox{} \hskip0.01cm SeDuMi} 
    \end{minipage}
    \begin{minipage}{2.4cm}
    \centering
    \includegraphics[width=2.3cm]{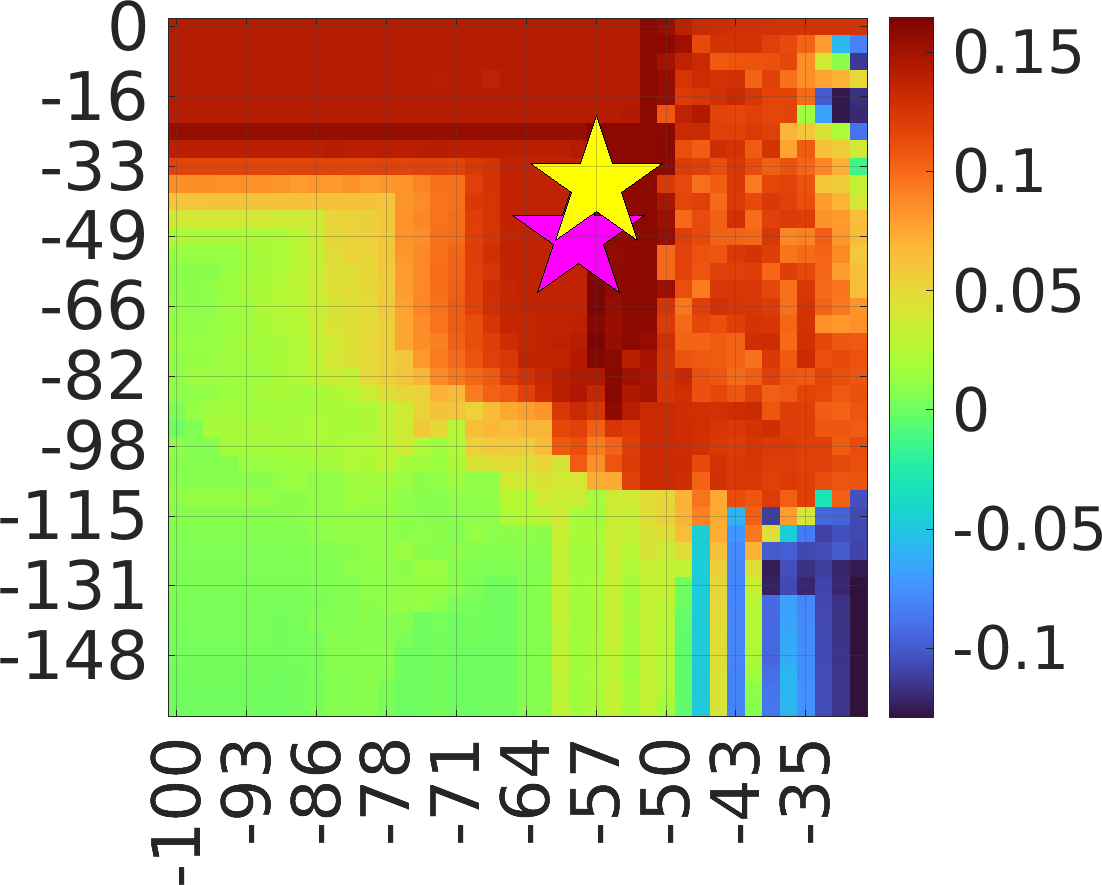}
    \end{minipage}
    \begin{minipage}{2.4cm}
    \centering
    \mbox{}
    \end{minipage}
    \begin{minipage}{2.4cm}
    \centering
    \mbox{}
    \end{minipage}
    \\
    \vskip0.1cm
    Dense $40 \times 40$ lattice
    \end{minipage}
    \end{scriptsize}
    \caption{{\bf Left:} Charts showing value of $\Gamma$ in the $15 \times 15$ coarser $\alpha \epsilon$-lattice of the two-stage meta-optimization process for different optimization processes. The yellow star shows the maximizer of the intensity $\Gamma_{\hbox{\scriptsize{max}}}$ and purple one that of the focality  $\Theta_{\hbox{\scriptsize{max}}}$. The maximizers are rather generally found with  $\varepsilon$-values (thershold levels) smaller than -30 dB, reflecting the mutual scale difference between the reference amplitude of ${\bf x}_1$ and the volumetric current intensities that can be generated within the given limits for the maximum dose. The regularization parameter $\alpha$ can be observed to have a sweet spot between -60 and -50 dB, where the intensity obtains its maximum. The maximum focality is found with a slightly smaller $\alpha$ compared to the intensity. Significant mutual differences can be observed between how the different optimization methods localize the sweet spot; in particular, the IP algorithms tend to find a smoother area with an elevated $\Gamma$ than the simplex methods. {\bf Right:} Charts showing the denser $40 \times 40$  $\alpha \epsilon$-lattice for $\Gamma$. Compared to the coarser sweet spot between -60 and -50 dB can be observed to be sharper and to have an enhanced amplitude, in particular, in the case of the DS and PS algorithm, suggesting that they might be less robust than IP with respect to the selection of $\alpha$ and $\epsilon$ parameter.}
    \label{fig:my_label_low_resolution_lattice} 
    \label{fig:my_label_high_resolution_lattice} 
\end{figure*}

Figure \ref{fig:my_label_theta_Gamma_case_A} shows the results obtained in the first and second optimization stage for the somatosensory, auditory and occipital target, demonstrating the performance differences between the optimization methods. As expected, DS and PS tend to be faster but overall less accurate compared to the IP algorithms, which is is clealy visible in the results obtained for $\Theta_{\hbox{\scriptsize{max}}}$, whereas the differences observed in the estimates of  $\Gamma_{\hbox{\scriptsize{max}}}$ are minor, especially when the denser lattice is used. Somewhat surprisingly, Matlab IPL, the legacy routine of  the LIPSOL package, is the winner method considering the amplitudes of both focality $\Theta_{\hbox{\scriptsize{max}}}$ and intensity  $\Gamma_{\hbox{\scriptsize{max}}}$ after two optimization runs. For the somatosensory target, each algorithm finds the bipolar montage which maximizes the intensity $\Gamma_{\hbox{\scriptsize{max}}}$, while this is not the case for the other two targets. Solutions with more than two active electrodes are found by Mosek IP, SeDuMi IP, SDPT3 IP, Mosek PS, and Gurobi DS with the denser $40 \times 40$ lattice and Mosek IP, SeDuMi IP, SDPT3 IP, Gurobi DS, Gurobi PS, Matlab DS, and Mosek DS with the coarser $15 \times 15$ lattice.  

Whereas the best mathematical perfomance is obtained with Matlab IPL, the open solver of the LIPSOL package, the commercial solvers are fast. Looking at the computing times in the first and second stages (Figures \ref{fig:my_label_computing_time_case_B}), the first stage takes up to 10000 seconds vs the second stage takes up to 1000 seconds for the denser lattice. It is obvious that the second stage must be faster to be evaluate, since it involves a limited montage and, thereby, a reduced lead field size which is roughly 15 \% of the original.  MOSEK IP, DS and PS turned out to be by far superior regarding the total computing time, MOSEK DS being the fastest of these three. The computing time of Gurobi IP was close to that of MOSEK IP, while Gurobi DS and PS a well as Matlab IPL and DS required approximately three times that.  The slowest performing SDPT3 and SeDuMi took as much as six  times the run time of MOSEK IP. 

For their fast performance, we applied MOSEK IP, DS and PS together with the $15 \times 15$ search lattice to evaluate maps of $\Theta_{\hbox{\scriptsize{max}}}$ and $\Gamma_{\hbox{\scriptsize{max}}}$ by solving the optimization task for each spatial and directional degree of freedom in the discretized source distribution (Figure \ref{fig:my_label_ratio_theta}). These maps were observed to be mutually similar between the three solvers. In the case of the focality  $\Theta_{\hbox{\scriptsize{max}}}$ (first row), the IP optimizer achieves enhanced amplitudes compared to DS and PS in the area, where finding a  focal stimulus can be considered feasible. In these areas, DS and PS yield a greater intensity which comes with the cost of suppressed focality and is less interesting  in the process to  find  $\Theta_{\hbox{\scriptsize{max}}}$ as it does not primarily aim at  maximizing the intensity but the focality instead. The maps obtained for $\Gamma_{\hbox{\scriptsize{max}}}$ (second row ) have  minor differences both topographically and amplitude-wise. Of  $\Theta_{\hbox{\scriptsize max}}$ and $\Gamma_{\hbox{\scriptsize max}}$, the former one corresponds to a greater NNZ and a smaller maximum current. A smaller number of NNZ in finding $\Gamma_{\hbox{\scriptsize max}}$ is expected, since the  optimizers often find the bipolar pattern.  Hence, the   overall current pattern has a greater coverage but a smaller maximal intensity in the case $\Theta_{\hbox{\scriptsize max}}$ compared to the case of $\Gamma_{\hbox{\scriptsize max}}$. 

Figure \ref{fig:my_label_high_resolution_lattice} shows the first-stage lattice search of the meta-optimization processes: $\Theta_{\hbox{\scriptsize max}}$ (purple star) is found using as a metacriterion the threshold condition corresponding to 75\% of the maximum amplitude obtainable with the bipolar two-channel montage, and $\Gamma_{\hbox{\scriptsize max}}$ is the global maximizer of $\Gamma$ (yellow star). The search results are shown for each optimization algorithm for two different lattice resolutions $15 \times 15$ and $40 \times 40$. In the denser lattice, a thin vertical line maximizing current amplitude is well visible. This sweet spot appears in the regularization parameter range from -60 to -50 dB (horizontal axis) and the nuisance field threshold range below -30 dB. Emphasizing the importance of a sufficient search resolution, the denser lattice contains some details that are absent in the coarser lattice which affect the optimization process. For the IP algorithms, the lattice structure is overall smoother and the differences between the two resolutions are less drastic than for DS and PS. Hence, it seems that DS and PS methods are less robust regarding the variation of the search parameters as compared to the IP techniques, partially explaining the mutual performance differences between IP, DS and PS observed in the numerical experiments of this study. 

The denser lattice not only increases the optimization accuracy but also the numerical stability of the results; the denser lattice deviates less in the vicinity of the optimum than the coarser one which can be observed in Figure \ref{fig:my_label_theta_Gamma_case_A} based on the whiskers showing  second order Taylor's polynomial  estimates for the maximal deviation within a 1/2 lattice size distance from the optimizer position. 

\section{Discussion}

In this study, we compared the performances of several IP \cite{mehrotra1992implementation}, DS and PS algorithms \cite{Boyd2004} to solve the recently introduced L1-L1 (L1 fitted and regularized) optimization task of tES \cite{herrmann2013transcranial}. These techniques are currently the predominating ones in solving Linear Programming (LP) problems, and therefore they are available in most open and commercial optimization toolboxes. The comparison was motivated by our earlier results  \cite{galazprieto2022l1}  which suggest that L1-L1 provides a theoretically attractive approach to obtain a focal stimulus with a focal current pattern and that it performs appropriately compared to complex L2-fitting and regularized least squares techniques \cite{galazprieto2022l1}. In particular, the optimization process itself is completely parameter free and does not necessitate any prior knowledge about where and how much the nuisance field needs to be suppressed, which is the case in the earlier LP formulations of the tES curent optimization problem \cite{wagner2016optimization,fernandez2020unification,dmochowski2011optimized}. In L1-L1, the nuisance field is settled optimally by a meta-optimization process that is performed in a two-dimensional lattice. 

We examined solvers of the commercial Matlab, Mosek \cite{mosek} and Gurobi  \cite{gurobi} toolboxes as well as  SDPT3 \cite{tutuncu2003solving},  SeDuMi (Self-Dual Minimization) \cite{sturm1999using,sturm2000central,polik2007sedumi}  which are openly available and accessible through the open CVX toolbox for Matlab \cite{grant2009cvx}. Even though commercial, Matlab IPL, i.e., Matlab's legacy IP solver, can also be regarded as an open, since it originates from the open LIPSOL  \cite{zhang1999user} toolbox. It was selected to this study, since Matlab's main IP algorithm was observed to stall in the middle of a two-stage meta-optimization process and not return an appropriate result at all. 

As a reference technique in the comparison, we used the generalized reciprocity principle (GRP) \cite{fernandez2020unification}, which states that, given the target, the current pattern corresponding to the maximum intensity can be obtained by picking the two electrodes with the largest absolute backprojected currents. The validity of GRP was shown for the present lead field matrix ${\bf L}$ (Appendix \ref{sec:grp}), since the reciprocity of propagation concerns an  electromagnetic field itself, measured for instance in EEG, but not necessarily the gradient of the field evaluated by ${\bf L}$ which projects the current flux through the boundary electrodes to a volumetric current density.

In the numerical experiments, we observed significant mutual differences between the solvers. These differences concern mainly the greatest obtainable focality  amplitudes and the number of active electrodes, whereas rather similar intensities and topography maps could be obtained regarless of the method.  A somewhat surprising finding was that Matlab IPL demonstrated superior mathematical performance in the mutual comparison between the different solvers; it was the winner after two optimization runs considering both $\Theta_{\hbox{\scriptsize max}}$ and $\Gamma_{\hbox{\scriptsize max}}$ for each target location. The commercial solvers were observed to be fast compared to the other methods, especially MOSEK IP, DS, and PS, which were the fastest of their kind in this study. While the best performing solvers show that the L1-L1 method is suitable for maximizing both focality and intensity, a few of them did not find the bipolar current pattern that maximizes $\Gamma_{\hbox{\scriptsize max}}$. Notably, SDTP3 did not find a bipolar pattern at all,  verifying our earlier hypothesis \cite{galazprieto2022l1} that the performance of L1-L1 might be highly solver-based. Part of the discrepancies between the optimization methods can obviously be explained by a different sensitivity with respect to parameter variation or the resolution of the meta-optimization lattice. 

The IP, DS and PS methods were natural choices for this study, since they are  available in  well-known optimization packages. As DS and PS require fewer machine resources than the IP methods and  provide an appropriate estimate  $\Gamma_{\hbox{\scriptsize max}}$, they might be beneficial in some applications, where the hardware performance is limited, for example in a potential FPGA implementation of the L1-L1 optimizer \cite{bayliss2006fpga,gensheimer2014simplex}. Otherwise, IP algorithms seem preferable due to their better mathematical performance for the L1-L1 task. Alternative solver algorithms might include Alternating Direction Method of Multipliers (ADMM) \cite{lin2021admm} which is increasingly applied to LP problems. ADMM was not included in this investigation, as achieving an appropriate convergence with ADMM seemed more difficult, e.g., due to its dependence on a step-length parameter. 

As for the future work, we currently plan to use L1-L1 in applications. One interesting application would be deep brain stimulation (DBS), where an advanced optimization technique is needed to target the subcortical nuclei of the brain. See, for example, a recent study \cite{anderson2018optimized}, where an IP algorithm has been applied.

  \section*{Acknowledgment}
  
  FGP, MS, and SP were supported by the Academy of Finland Center of Excellence in Inverse Modelling and Imaging 2018--2025, DAAD project (334465) and by the ERA PerMed project PerEpi (344712). 
  
  \begin{appendix}
  
  \section{Generalized reciprocity principle for a tES lead field matrix}
  
  \label{sec:grp}
  
GRP can be formulated for a restricted system ${\bf L} {\bf R}_K {\bf y}_K = {\bf x}$, where ${\bf R}_K$ denotes a real $N\times K$ ($K \leq N$) restriction matrix, whose nonzero entries $r_{i_j,j} = 1$ corresponds to an ordered subset of electrodes $\mathcal{S} = \{i_j \, : \, j = 1, 2, \ldots, K, \hbox{ with } |({{\bf L}_1^T {{\bf x}_1})_{i_1}| \geq |({{\bf L}_1^T {\bf x}_1})_{i_2}| \geq \cdots \geq |({\bf L}_1^T {\bf x}_1})_{i_K} \}$.  GRP follows by writing the intensity as $\Gamma = \sigma_K \, {\bf y}_K^T \,  {\bf s}_K$ with $\sigma_K = \| {\bf R} {\bf L}^T_1 {\bf x}_1 \|_1 / \| {\bf x}_1 \|_2 $ and ${\bf s}_K = {\bf R}  {\bf L}^T_1 {\bf x}_1 / \| {\bf R} {\bf L}^T_1 {\bf x}_1\|_1$ which shows that $\Gamma$ can be interpreted as a projection of ${\bf y}_K$ on  $\sigma_K \, {\bf s}_K$. 
 
 Thus, the  maximum of $\Gamma$ is obtained, when ${\bf y}_K$ is  parallel to ${\bf s}_K$. Since the maximizer also needs to match the maximum applicable dose $\mu$, as otherwise it can be upscaled, it follows that   ${\bf y}_K = \mu \, {\bf s}_K$. The corresponding maximum intensity is $\Gamma = \mu\, \sigma_K \,  \|  {\bf s}_{K} \|_2^2$.  The optimal montage is the maximizer of  
\begin{equation}
\max_{K} \mu \, \sigma_K \|  {\bf s}_K \|_2^2, 
\end{equation}
where by definition  $\| {\bf s}_K  \|_1 = 1$ and the entries of ${\bf s}_K$ are ordered in a descending order with respect to their absolute value. Assuming that these entries are given by $\sum_{j = 1}^K \lambda_j$ with $\lambda_1 \geq \lambda_2 \geq \cdots \geq \lambda_K \geq 0$ it holds that $\|\sum_{j = 1}^{K}  \lambda_{i_j}^2  \|$ 
{\setlength\arraycolsep{2pt} \begin{eqnarray}
\| {\bf s}_K \|^2_2 - \| {\bf s}_{K-1} \|^2_2  & = & \! \frac{\lambda_K}{1 - \lambda_K} \!  \left( \lambda_K^2 \! - \! \lambda_K \! + \! 2 \sum_{j = 1}^{K-1}  \lambda_{i_j}^2  \right) \nonumber\\
& \geq & \! \frac{\lambda_K}{1 - \lambda_K} \! \left( K \lambda_K^2 \! - \! \lambda_K \!+\!  \sum_{j = 1}^{K-1}  \lambda_{i_j}^2  \right).
\label{reciprocal}
\end{eqnarray}}
Here the equality follows from a straightforward substitution and the inequality is obtained as $(K-1) \lambda_K^2 \leq \sum_{j=1}{K-1} \lambda_{i_j}^2$. Following from the discriminant together with the Arithmetic Mean -- Quadratic Mean inequality \begin{equation} \frac{1}{K-1} \sum_{j = 1}^{K-1} \lambda_{i_j}^2 \geq \left( \frac{1}{1-K} \sum_{j=1}^{K-1} \lambda_{i_j} \right)^2, \end{equation} the second factor in (\ref{reciprocal}) does not have roots if 
\begin{equation}
    K \sum_{j = 1}^{K-1}  \lambda_{i_j}^2  \geq \left(  \sum_{j = 1}^{K-1}  \lambda_{i_j} \right)^2 \geq \frac{1}{4}.
\end{equation}
That is, if $\sum_{j=1}^{K-1} \lambda_{i_j}  \geq \frac{1}{2}$. This assumption is valid, since a montage can have  minimally two active channels, none of which can contain more than half of the total dose, since the sum of currents for any montage is assumed to be zero. Hence, it follows that $\| {\bf s}_K \|^2_2 - \| {\bf s}_{K-1} \| \geq 0 $ for any montage and, recursively, that the maximum of $\Gamma$ is obtained with the bipolar pattern that corresponds to the first two entries $i_1$ and $i_2$ in the set $\mathcal{S}$.

---
\end{appendix}
  
\bibliographystyle{IEEEtran}
\bibliography{references}

\end{document}